\documentclass[11pt]{article}
\usepackage[a4paper,top=1.05in,bottom=1.05in,left=1in,right=1in]{geometry}
\usepackage{amsmath,amssymb,amsthm,mathtools}
\allowdisplaybreaks[4]
\usepackage[T1]{fontenc}\usepackage{newtxtext,newtxmath,microtype}
\usepackage{titlesec}
\titleformat{\section}{\normalfont\Large\bfseries}{\thesection}{0.75em}{}
\titleformat{\subsection}{\normalfont\large\bfseries}{\thesubsection}{0.75em}{}
\titlespacing*{\section}{0pt}{2.6ex plus 0.8ex minus 0.2ex}{1.2ex plus 0.2ex}
\titlespacing*{\subsection}{0pt}{2.0ex plus 0.6ex minus 0.2ex}{0.9ex plus 0.2ex}
\usepackage{fancyhdr}
\fancypagestyle{plain}{\fancyhf{}\fancyfoot[C]{\thepage}}
\usepackage[dvipsnames]{xcolor}
\definecolor{linkblue}{RGB}{20,67,128}
\definecolor{citegreen}{RGB}{25,105,80}
\definecolor{urlpurple}{RGB}{100,55,125}
\usepackage[colorlinks=true,linkcolor=linkblue,citecolor=citegreen,urlcolor=urlpurple,plainpages=false,hypertexnames=false]{hyperref}
\usepackage{bookmark}
\numberwithin{equation}{section}
\newtheorem{thm}{Theorem}[section]\newtheorem{lem}[thm]{Lemma}
\newtheorem{prop}[thm]{Proposition}\newtheorem{cor}[thm]{Corollary}
\theoremstyle{remark}
\newcommand{\R}{\mathbb R}\newcommand{\C}{\mathcal C}\newcommand{\V}{\mathcal V}
\newcommand{\dd}{\,\mathrm d}\newcommand{\norm}[1]{\left\|#1\right\|}
\DeclareMathOperator{\dist}{dist}
\title{Uniqueness of tangent flows in mean curvature flow (I):\\ smooth minimal cone}
\author{Hanbing Fang\qquad and \qquad Yongheng Han}
\date{\today}
\begin{document}
	\maketitle
	\begin{abstract}
		We prove uniqueness of multiplicity-one minimal-cone tangent flows in
		mean curvature flow when the cone has a smooth embedded two-sided link and satisfies strict Hardy
		stability, radial Gaussian nonresonance, and weighted stationary rigidity.
		We also establish a discrete Lojasiewicz inequality, polynomial convergence
		rates, and rigidity among smooth self-shrinkers. In particular, these conclusions hold for all Simons cones $\mathcal C_{p,q}$ with $p,q \geq 1, p+q\geq 6$.
	\end{abstract}
	\begingroup\small\tableofcontents\endgroup
	
	\section{Introduction}
	\subsection{Main theorem}
	Mean curvature flow is the evolution of an immersion $X(\cdot,t):M^m\to\R^{m+1}$ by
	\[
	\partial_tX=\mathbf H,
	\]
	where $\mathbf H$ is its mean curvature vector. A smooth flow starting from
	a closed hypersurface develops a singularity in finite time: if
	$M_0\subset B_R(x)$, the avoidance principle compares $M_t$ with the shrinking
	sphere $\partial B_{\sqrt{R^2-2mt}}(x)$ and gives $T\le R^2/(2m)$.
	The analysis of such singularities is therefore a central problem in mean
	curvature flow. Huisken's monotonicity formula shows that subsequential
	limits of parabolic rescalings at a singular point are self-similarly
	shrinking flows \cite{Huisken}. In the normalization $t=-1$, their slices
	satisfy
	\[
	\mathbf H+\tfrac12X^\perp=0.
	\]
	A basic question is whether this tangent flow is independent of the
	rescaling sequence. Uniqueness determines a single asymptotic model for the
	singularity and permits quantitative estimates of convergence to that model.

	For smooth tangent slices, Schulze \cite{Schulze} proves uniqueness when
	one tangent flow is a closed, smoothly embedded multiplicity-one
	self-shrinker, in arbitrary codimension. Colding--Minicozzi \cite{CM}
	prove uniqueness for multiplicity-one round cylindrical tangent flows of
	hypersurface mean curvature flow, including uniqueness of the axis. 
	Recently, Bamler--Lai \cite{BamlerLai} use their PDE--ODI principle
	to give an alternative proof of the uniqueness of round cylindrical tangent
	flows. Zhu \cite{Zhu} establishes Lojasiewicz inequalities and uniqueness for a
	class of cylindrical shrinkers in arbitrary codimension, including round
	cylinders and cylinders over Abresch--Langer curves. For noncompact tangent
	slices with conical ends, Chodosh--Schulze \cite{CS} prove uniqueness of
	multiplicity-one asymptotically conical tangent flows in codimension one;
	Lee--Zhao \cite{LZ} extend this conclusion to arbitrary codimension.
	These results develop the relation between analytic gradient inequalities
	and convergence established by Simon \cite{Simon}.
	
	When the tangent slice itself is singular, smooth graphical convergence
	fails near its singular set. For Brakke flows satisfying their structural hypothesis,
	Stuvard--Tonegawa \cite{StuvardTonegawa} prove that if one
	tangent flow is a static multiplicity-one triple junction, then it is
	the unique tangent flow. Lotay--Schulze--Sz\'ekelyhidi \cite{LSS} prove
	uniqueness at the first singular time of a zero-Maslov, rational Lagrangian
	mean curvature flow in a compact Calabi--Yau surface when one tangent flow
	is the static union of two transverse planes. For exact Lagrangian mean curvature flow in $\mathbb C^n$ with uniformly
	bounded Lagrangian angle and area ratios, Li--Sz\'ekelyhidi
	\cite{LiSzekelyhidiLMCF} prove uniqueness when one tangent flow
	is a multiplicity-one integrable special Lagrangian cone with a smooth
	connected link. They also prove uniqueness for tangent flows given by
	special Lagrangian unions of pairwise transverse planes in $\mathbb C^2$. Stolarski \cite{Stolarski}
	proves uniqueness of multiplicity-one stationary conical tangent flows
	with smooth links under locally uniform mean-curvature bounds, for integral
	Brakke flows in arbitrary codimension. In that setting, the curvature bound
	allows a reduction to uniqueness of tangent cones for varifolds with
	bounded mean curvature.

	There are parallel convergence and uniqueness results for Ricci flow.
	Sun--Wang \cite{SunWang} prove polynomial convergence, modulo
	diffeomorphisms, of normalized K\"ahler--Ricci flows starting sufficiently
	close to a positive K\"ahler--Einstein metric, using a Lojasiewicz inequality
	for Perelman's $\mu$ functional. Chan--Ma--Zhang \cite{CMZ} prove uniqueness
	when a tangent flow is closed and smooth. For noncompact models,
	Li--Wang \cite{LiWang} prove that a Ricci shrinker sufficiently close to the
	round cylinder $S^{n-1}\times\R$ in pointed Gromov--Hausdorff topology is
	isometric to it. Colding--Minicozzi \cite{CMRicci} establish strong rigidity
	of cylindrical Ricci shrinkers and show that if one tangent flow at a future
	singular point is a cylinder, then all tangent flows there are cylinders.
	Fang--Li \cite{FL} prove a pointed entropy Lojasiewicz inequality at
	cylindrical singularities and obtain strong uniqueness, including convergence
	in a fixed parametrization. These results connect rigidity of stationary
	models with estimates controlling an evolving flow between successive
	rescalings. More recently, Li--Yan--Zhang \cite{LYZAC}
	prove strong uniqueness at the first singular time of a closed Ricci
	flow when one tangent flow is generated by a complete smooth
	asymptotically conical Ricci shrinker. Their modified flow converges
	in a fixed gauge, with quantitative estimates on expanding domains.
	
	In this paper, we prove uniqueness when one tangent flow is a static
	multiplicity-one minimal cone, under strict Hardy stability, radial
	Gaussian nonresonance, and weighted stationary rigidity; see
	Theorem~\ref{thm:main}. We also prove a discrete Lojasiewicz inequality,
	quantitative convergence estimates, and rigidity among smooth self-shrinkers.
	The corresponding stationary problem is uniqueness of tangent
	cones. Allard--Almgren \cite{AllardAlmgren} proved uniqueness and a power rate
	of convergence for multiplicity-one tangent cones with smooth links
	under an integrability assumption. Simon \cite{Simon} proves uniqueness for multiplicity-one tangent
	cones with smooth links by analytic convergence methods; his later work
	\cite{SimonCyl} treats classes of cylindrical cones with nonisolated
	singularities. For more work related to uniqueness of tangent cones and the geometry
	of minimal hypersurfaces near singular cones, see
	\cite{HS,ES,Szekelyhidi,Wang}.

	\paragraph{Main result.}
	Let $m\ge3$ and let $\Sigma^{m-1}\subset S^m$ be a smooth, closed, connected,
	embedded, two-sided minimal hypersurface. Put
	\[
	\C:=\{r\omega:r\ge0,\ \omega\in\Sigma\}\subset\R^{m+1}.
	\]
	Every cone in this paper has multiplicity one. We use $|\C|$ for the
	associated integral varifold and $\mathcal H^m$ for $m$-dimensional
	Hausdorff measure. Choose a unit normal
	$n_\Sigma$. We denote the unit normals to $\C$ and $N$ by
	$n_\C$ and $n_N$, respectively. We write $A_\Sigma$ for the second fundamental form and $H_\Sigma$
	for scalar mean curvature. With $\Delta_\Sigma=\operatorname{div}_\Sigma\nabla_\Sigma$
	nonpositive, we define the Jacobi operator and its kernel by
	\[
	L_\Sigma:=\Delta_\Sigma+|A_\Sigma|^2+(m-1),\qquad
	K:=\ker L_\Sigma,\qquad d_K:=\dim K.
	\]
	The orthogonal complement $K^\perp$ and projections are taken with respect
	to the fixed $L^2(\Sigma,\dd\mu_\Sigma)$ inner product. For any orthonormal
	basis $\{\psi_l\}_{l=1}^{d_K}$ of $K$,  define
	\begin{equation}
		P_Kv:=\sum_{l=1}^{d_K}\left(\int_\Sigma v\psi_l\dd\mu_\Sigma\right)\psi_l,
		\qquad P_{K^\perp}v:=v-P_Kv.
	\end{equation}
	Let $(\mu_j,\phi_j)$ be an $L^2(\Sigma)$ orthonormal eigensystem of
	$-L_\Sigma+(m-1)=-\Delta_\Sigma-|A_\Sigma|^2$, with
	$\mu_0$ the lowest eigenvalue and $\phi_0>0$. We use the following cone
	conditions.
	\begin{itemize}
		\item \textbf{Strict Hardy stability.}
		\begin{equation}\tag{S1}
			\mu_0>-(m-2)^2/4.
		\end{equation}
		Under (S1),  define
		\[
		\nu_j:=\sqrt{(m-2)^2/4+\mu_j},\qquad
		\gamma_j^\pm:=-(m-2)/2\pm\nu_j.
		\]
		The numbers $\gamma_j^\pm$ are the homogeneities of the ordinary Jacobi
		equation in the $j$th angular eigenspace.
		
		\item \textbf{Radial Gaussian nonresonance.}
		\begin{equation}\tag{S2}
			(1-2k)(m-1-2k)\notin\operatorname{Spec}(-L_\Sigma+(m-1))
			\quad(1\le k<m/4,\ k\in\mathbb Z).
		\end{equation}
		This excludes additional radial kernel functions of the Gaussian Jacobi
		operator; see Proposition~\ref{prop:spectrum}.
		
		\item \textbf{Weighted stationary rigidity.}
		\begin{equation}\tag{S3}
			\left.
			\begin{gathered}
				V\ne0\text{ is stationary and integral, with finite area ratios},\\
				(D_{1/R})_\#V\longrightarrow|\C|\quad(R\longrightarrow\infty),\\
				\displaystyle\int_{\R^{m+1}}\frac{\dist(X,\C)^2}{1+|X|^2}\dd\|V\|(X)<\infty
			\end{gathered}\right\}\quad\Longrightarrow\quad V=|\C|.
		\end{equation}
		Here $D_{1/R}(X)=X/R$, and the convergence is full multiplicity-one
		varifold blowdown convergence. The condition applies to every stationary
		integral varifold satisfying these bounds, including those that are
		neither stable nor area-minimizing.
		
	\end{itemize}
	
	Conditions (S1) and (S2) identify the Gaussian kernel at the base
	cone. Condition (S3) excludes nontrivial stationary limits arising
	from concentration near the vertex. Now we can state our main theorem as follows.
	
	\begin{thm}[Uniqueness of tangent flows]\label{thm:main}
		Suppose $\C$ satisfies \textup{(S1), (S2), and (S3)}.
		Let $M_t$, $t<T$, be a smooth embedded boundaryless mean curvature flow
		which is spatially and spacetime proper and satisfies
		\begin{equation}\label{equ:areabd}
			\mathcal H^m(M_t\cap B_R(x))\le\Lambda R^m
			\quad\text{for every }t<T,\ x\in\R^{m+1},\ R>0.
		\end{equation}
		If one tangent flow at $(x_0,T)$ is the static multiplicity-one cone $\C$,
		then every tangent flow there is the same cone in the original coordinates.
	\end{thm}
	
	We verify the cone hypotheses for two classes of examples. Strict
	minimization gives (S3) by a stationary asymptotic-rate bound, while
	quadratic cones can be treated by their explicit link spectra and
	stationary rigidity; see Section~\ref{sec:examples}. Here strict
	minimization is the Hardt--Simon strengthening of area minimization:
	the canonical minimizing leaves on both sides of the cone have a
	nonzero leading term at the larger principal indicial root
	$\gamma_0^+$. We give the precise definition in
	Subsection~\ref{sec:strict-minimization}.
	
	\begin{cor}\label{cor:examples}
		The conclusion in Theorem~\ref{thm:main} holds for a strictly minimizing
		cone satisfying \textup{(S1) and (S2)}. It also holds for
		\[
		\C_{p,q}=\{(x,y)\in\R^{p+1}\oplus\R^{q+1}:q|x|^2=p|y|^2\},
		\quad p,q\ge1,\quad m=p+q+1\ge7.
		\]
	\end{cor}

	We next give quantitative consequences of the estimates used to prove
	uniqueness. To state them, at the base point $(x_0,T)$ we define the \textbf{rescaled mean curvature flow} and the \textbf{shrinker quantity} by
	\begin{equation}\label{eq:flow}
		N_s:=e^{s/2}(M_{T-e^{-s}}-x_0),\quad
		\rho:=(4\pi)^{-m/2}e^{-|X|^2/4},\quad
		\Phi:=\mathbf H+\frac{X^\perp}{2}.
	\end{equation}
	We use $c_m:=(4\pi)^{-m/2}$ for the Gaussian normalization constant.
	For a hypersurface $N$, we write
	\begin{equation}
		F(N):=\int_N\rho\,\dd\mu,\quad G(N)^2:=\int_{N}|\Phi|^2\rho\, \dd\mu.
	\end{equation}
	Along the rescaled flow, we define
	\begin{align}
		E(s):=F(N_s)-F(\C).
	\end{align}

	Huisken's monotonicity formula and the prescribed tangent flow give
	$E\ge0$, $E\to0$, and $E'(s)=-G(N_s)^2$. The parameter $a(s)$ fits the analytic reduced family
	$\C_a$ constructed in Proposition~\ref{prop:link-reduction}.
We choose $a(s)$ so that the graph height of $N_s$ over $\C_{a(s)}$
is orthogonal to the normal variations of the cone family on a fixed
annulus, with the cutoff Gaussian weight of Lemma~\ref{lem:gauge}. The reference links
	need not be minimal.
	The Gaussian distance and the decrease of Gaussian area on a unit interval are
	\[
	D_a(N)^2:=\int_N\dist(X,\C_a)^2\rho\dd\mu,\quad
	d_j:=D_{a(j)}(N_j),\quad \delta_j^2:=E(j)-E(j+1).
	\]
	
	Theorem~\ref{thm:energy-loja} bounds the Gaussian area excess by its decrease on
	each sufficiently late interval of length two. It also gives a bound for
	the remaining length of the rescaled flow and a rate of decay. For related discrete Lojasiewicz inequalities at cylindrical
	singularities in mean curvature flow and Ricci flow, see
	Colding--Minicozzi \cite[Theorem~6.1]{CM} and
	Fang--Li \cite[Theorem~1.3]{FL}, respectively.
	
	\begin{thm}[Discrete Lojasiewicz inequality]\label{thm:energy-loja}
		Under the hypotheses of Theorem~\ref{thm:main}, there are
		$s_0$, $C<\infty$, and $0<\beta<1/2$ such that, for every real $s\ge s_0$, we have
		\begin{equation}\label{eq:energy-loja}
			E(s)^{2(1-\beta)}\le C\left(E(s)-E(s+2)\right).
		\end{equation}
		Moreover, we have
		\begin{equation}\label{eq:energy-length}
			\begin{gathered}
				\int_s^\infty G(N_{\sigma})\dd\sigma\le CE(s)^\beta,\\
				E(s)\le C(1+s-s_0)^{-1/(1-2\beta)}
				\quad\text{for }0<\beta<1/2.
			\end{gathered}
		\end{equation}
	If $E(s_*)=0$ at a finite time $s_*$, then the rescaled flow has zero
	length on $[s_*,\infty)$.
	\end{thm}
	The inequality holds along the rescaled flow on every sufficiently late
	time interval of length two. We first prove that the flow remains in the
	neighborhood where the local estimates apply, using a continuation argument
	on finite intervals. We then obtain \eqref{eq:energy-loja} and the
	finite-length estimate \eqref{eq:energy-length}.
	
	The local estimates also apply to a constant rescaled flow. Since its
	shrinker quantity vanishes identically, they yield the following rigidity
	result.
	
	\begin{cor}\label{cor:shrinker-rigidity}
		Let $\C$ be nonflat and satisfy the cone hypotheses of
		Theorem~\ref{thm:main}. There is no sequence of smooth, properly embedded,
		boundaryless hypersurfaces $S_i$ satisfying
		$\mathbf H_{S_i}+X^\perp/2=0$ that converges locally as varifolds to
		$|\C|$ with multiplicity one. Equivalently, a local varifold neighborhood
		of $|\C|$ contains no such smooth normalized shrinker.
	\end{cor}
	We prove Corollary~\ref{cor:shrinker-rigidity} by applying the local
	residual contraction estimate to a constant rescaled flow,
	for which the shrinker quantity vanishes.
	
	In the forthcoming Part~II, we prove uniqueness at singular points
	admitting a static multiplicity-one tangent flow of the form
	$\C\times\R^k$, where $\C$ is a minimal cone with smooth link and $k\ge1$.
	The result assumes strict Hardy stability, spectral conditions on the
	link, and a weighted stationary rigidity condition for the product cone.
	In particular, it applies to $\C_{p,q}\times\R^k$ for $p,q\ge1$ and
	$p+q\ge7$.
	
	\subsection{Strategy and organization}
	We need to control both the deformation of the link and possible concentration
	near the vertex. We first separate these two issues,
	then explain how the resulting distance estimates imply convergence and
	the quantitative conclusions.
	
	\paragraph{(A) Analytic reduction and fitting the reference family.}
	The fixed cone $\C$ has degree-one Gaussian Jacobi fields
	$r\psi(\omega)$, $\psi\in K$. These fields lie in the kernel of the
	linearized shrinker equation, so an inverse estimate cannot control
	the corresponding component of the height. We therefore choose a
	reference cone $\C_{a(s)}$ that varies with the flow, and
	estimate the remaining height relative to that cone. This choice
	separates the finite-dimensional link parameter from the part to which
	an inverse estimate can apply.
	
	In Proposition~\ref{prop:link-reduction}, we construct the reference links
	by solving the complementary spherical equation. More precisely, {we have}
	\[
	k(a)=a+\eta(a),\qquad
	\eta(a)\in K^\perp,\qquad
	P_{K^\perp}\mathcal G(k(a))=0.
	\]
	The reduced area $f(a)=\mathcal A(k(a))$ satisfies
	\[
	q(a):=\nabla f(a)=\mathcal G(k(a))\in K,
	\qquad
	|f(a)-f(0)|^{1-\theta}\le C|q(a)|,
	\quad 0<\theta\le1/2.
	\]
	Their mean curvature and Gaussian area must therefore be included in the
	estimates:
	\[
	H_{\C_a}=-\frac{q(a)}{r\mathfrak w_a},\qquad
	|F(\C_a)-F(\C)|^{1-\theta}\le C|q(a)|;
	\]
	see \eqref{eq:reduced-cone-curvature} and
	\eqref{eq:reduced-cone-gradient-comparison}. The smooth positive
	factor $\mathfrak w_a$ is uniformly bounded above and below.
	
	To determine the parameter, we use only a fixed regular annulus: the
	flow is graphical there even when it fails to be graphical near the
	vertex. Let $u_s$ be its height over $\C_{a(s)}$, and let
	$\Psi_{a,l}$ be the normal velocities of this conical family. The
	compact observations in Lemma~\ref{lem:gauge} determine $a(s)$ by
	\[
	\int_{\C_{a(s)}}\chi\Psi_{a(s),l}u_s\rho\dd\mu=0,
	\qquad 1\le l\le d_K.
	\]
	At $a=0$ these velocities are $r\psi_l$. Differentiating the system
	and using \eqref{eq:fitting-matrix-positive}, we obtain
	\[
	|a'(s)|\le CG(N_s),\qquad
	|a(s_2)-a(s_1)|\le C\int_{s_1}^{s_2}G(N_\sigma)\dd\sigma.
	\]
	Thus it follows that the observations fix the degree-one directions in the normalized
	limit, while the speed estimate controls motion of the reference cone
	by the same quantity used to estimate the flow.
	
	\paragraph{(B) Normalized compactness.}
	After fitting, we seek a contraction estimate for distance together
	with the link residual. On a unit time interval, fix the terminal
	parameter $a_i=a(N_1^i)$ and put
	\[
	d_i:=D_{a_i}(N_0^i),\qquad e_i:=D_{a_i}(N_1^i),\qquad
	q_i:=q(a_i),\qquad
	\delta_i^2:=\int_0^1G(N_s^i)^2\dd s.
	\]
	The estimate to be proved is
	\[
	e_i^2+|q_i|^2\le\eta d_i^2+C_\eta\delta_i^2
	\qquad(\eta>0).
	\]
	A failure of this estimate gives a sequence for which
	$b_i^2:=d_i^2+|q_i|^2\to0$ and $\delta_i/b_i\to0$.
	By dividing the graph height by $b_i$, we retain its possible nonzero limit
	while making its time derivative vanish. The difficulty is that
	convergence on punctured annuli does not account for distance
	concentrated at the vertex. We therefore pass to the limit in both
	the heights and their squared-distance measures.
	
	Proposition~\ref{prop:normalized} gives, after pullback to $\C$ and
	passage to a subsequence,
	\[
	\begin{gathered}
		q_i/b_i\longrightarrow\xi\in K,\qquad
		u_i/b_i\longrightarrow v,\qquad L_\C v=r^{-1}\xi,\qquad
		\mathcal O(v)=0,\\
		b_i^{-2}\dist(X,\C_{a_i})^2\rho\dd\mu_{N_1^i}
		\rightharpoonup v^2\rho\dd\mu_\C+\alpha_{\rm vtx}\delta_0.
	\end{gathered}
	\]
	The forcing $r^{-1}\xi$ comes from the reference mean curvature.
	The punctured field need not yet belong to the Friedrichs form
	domain, and $\alpha_{\rm vtx}$ may be positive. These are the two
	quantities that must be estimated near the vertex. At infinity,
	Lemma~\ref{lem:residual-tail} gives relative tightness, including
	the additional term $|q_i|^2$ in the tail bound.
	
	\paragraph{(C) Weighted rigidity and the vertex estimates.}
	To estimate these two quantities, we dilate at the largest radius
	where the small single-graph description fails. Let $\varrho_i$ be
	this radius on a selected slice $N_{s_i}^i$, where
	$G(N_{s_i}^i)\le2\delta_i$. The annular estimate in
	Lemma~\ref{lem:residual-neck} gives
	\[
	\begin{aligned}
		\mathcal H_{\C_{a_i}}
		(u_i;\{2\varrho_i<r<r_0\})
		&\le C\varrho_i^m
		+Cr_0^{-2}\int_{\{r_0<r<2r_0\}}u_i^2\dd\mu
		+Cr_0^2G(N_{s_i}^i)^2+Cr_0^m|q_i|^2\\
		&\le C(\varrho_i^m+b_i^2).
	\end{aligned}
	\]
	The last term in the first line is the contribution of the nonzero
	reference mean curvature. It is controlled at the normalization
	scale because $|q_i|\le b_i$.
	
	If $\varrho_i^m$ were not $o(b_i^2)$, dilation at a comparable bad
	radius would give a nonzero stationary integral varifold with finite
	weighted Hardy distance. The Gaussian scale comparison in
	Lemma~\ref{lem:scale-comparison} gives its full multiplicity-one
	blowdown $\C$. Condition (S3) then identifies the stationary limit
	with $|\C|$. Lemma~\ref{lem:moving-compactness} upgrades this
	convergence to complete smooth convergence at the selected time,
	contradicting the failed graphical radius. Lemma~\ref{lem:first-bad}
	therefore gives
	\[
	\varrho_i^m=o(b_i^2),\qquad
	\int_{N_{s_i}^i\cap B_{\epsilon_H}}
	\frac{\dist(X,\C_{a_i})^2}{|X|^2}\dd\mu
	\le C(\varrho_i^m+b_i^2).
	\]
	
	The resulting gradient and Hardy bounds put $v$ in the Friedrichs
	form domain. Proposition~\ref{prop:spectrum} identifies the kernel
	there with $rK$. The residual columns $r^{-1}K$ have a nonsingular
	pairing with this kernel, and the compact observations are nonsingular
	on the same space. Proposition~\ref{prop:gaussian-bordered} consequently
	gives
	\[
	\|v\|_\V+|\xi|
	\le C\left(\|L_\C v-r^{-1}\xi\|_{\V^*}
	+|\mathcal O(v)|\right)=0.
	\]
	The localized first variation estimate \eqref{eq:lw-localized-transport}
	transfers the Hardy bound to the terminal time and gives
	$\alpha_{\rm vtx}=0$. Proposition~\ref{prop:residual-contraction}
	converts this contradiction argument into the required uniform
	contraction estimate.
	
	\paragraph{(D) The local distance and Gaussian area estimates.}
	The contraction estimate controls the terminal distance by the initial
	distance and the defect. To obtain an estimate by the defect alone,
	we also need a comparison in the other direction. The truncated
	Gaussian distance permits this comparison because its bounded
	integrand controls the quadratic term in first variation.
	Lemmas~\ref{l:lem:transport} and \ref{l:lem:comparison} give,
	on $[s,s+2]$ with $a=a(N_{s+2})$ fixed,
	\[
	d\le C_0(e+\Delta),\qquad
	d:=D_a(N_{s+1}),\quad e:=D_a(N_{s+2}),\quad
	\Delta^2:=E(s)-E(s+2).
	\]
	Apply the contraction to the second unit interval and choose its
	constant $\eta$ after $C_0$ has been fixed. Then we have
	\[
	\begin{aligned}
		e^2+|q(a)|^2&\le\eta d^2+C_\eta\Delta^2,\\
		d&\le C_0\sqrt{\eta}\,d+C\Delta,
		\qquad d+e+|q(a)|\le C\Delta.
	\end{aligned}
	\]
	Here $C_0\sqrt{\eta}<1/2$, so the first term in the second line
	is absorbed. This is the estimate in
	\eqref{eq:residual-two-unit-distance}.
	
	The reference Gaussian area is not assumed to equal $F(\C)$.
	Instead, the spatial estimate in Lemma~\ref{lem:residual-spatial}
	retains its contribution:
	\[
	|F(N)-F(\C_a)|
	\le C\left(D_a(N)^2+G(N)^2+|q(a)|^2\right).
	\]
	Averaging on the last quarter of the second unit interval and using
	monotonicity, we obtain \eqref{eq:residual-terminal-area}. Thus we obtain
	\[
	\begin{aligned}
		E(s)&\le |F(\C_a)-F(\C)|+C\Delta^2,\\
		E(s)^{1-\theta}
		&\le C\left(|q(a)|+\Delta^{2(1-\theta)}\right)\le C\Delta.
	\end{aligned}
	\]
	The second line uses the reduced gradient inequality and
	$2(1-\theta)\ge1$. After squaring the bound in the second line, we obtain
	the local discrete inequality in Proposition~\ref{prop:local-loja}. After decreasing the energy
	neighborhood, we obtain the same inequality with any fixed
	$0<\beta<1/2$ no larger than $\theta$.
	
	\paragraph{(E) Finite length and uniqueness.}
	The preceding estimates hold only while the flow is close to the
	reference family. Summing these estimates over all time intervals of length two would
	require the first exit time from this neighborhood to be infinite. Instead, we work
	up to a possible first exit time and obtain bounds independent of
	that time.
	
	On the interval $[t_0+2j,t_0+2(j+1)]$, we write
	$e_j=E(t_0+2j)$ and $\delta_j^2=e_j-e_{j+1}$. For $e_j>0$, the local inequality
	and concavity give
	\[
	e_j^{1-\beta}\le C\delta_j,\qquad
	e_j^\beta-e_{j+1}^\beta
	\ge\beta e_j^{\beta-1}\delta_j^2\ge c\delta_j.
	\]
	If $e_j=0$, the subsequent drops vanish and no division is needed.
	We sum only over the complete time intervals preceding a hypothetical first
	exit $T_*$. On the final incomplete time interval, exact dissipation and
	Cauchy--Schwarz give a bound without any local distance estimate.
	Therefore we have
	\[
	\int_{t_0}^{T_*}G(N_s)\dd s
	\le CE(t_0)^\beta+\sqrt{2E(t_0)}.
	\]
	By starting sufficiently close to the prescribed tangent subsequence,
	we make this bound smaller than the chosen neighborhood size. Forward
	transport relative to the fixed base cone then prevents the first
	exit. The annular-entry estimate in Proposition~\ref{prop:anchor}
	maintains the fitting chart throughout this argument.
	Proposition~\ref{prop:continuation} therefore gives
	$\int_{t_0}^\infty G(N_s)\dd s<\infty$.
	
	For a compactly supported smooth function $\varphi$, first variation
	gives
	\[
	\begin{aligned}
		\left|\int_{N_{s_2}}\varphi\rho\dd\mu
		-\int_{N_{s_1}}\varphi\rho\dd\mu\right|
		&\le C_\varphi\int_{s_1}^{s_2}G(N_s)\dd s
		+\|\varphi\|_\infty\int_{s_1}^{s_2}G(N_s)^2\dd s\\
		&\longrightarrow0\qquad(s_2>s_1\to\infty).
	\end{aligned}
	\]
	The local measure limit is unique, and the prescribed tangent
	subsequence identifies it with $|\C|$ in the original coordinates.
	Conical homogeneity then identifies every parabolic tangent flow
	with the same static cone. This proves Theorem~\ref{thm:main}.
	
	\paragraph{(F) Strictly minimizing and quadratic cones.}
	To apply the uniqueness theorem, we verify the spectral assumptions and
	weighted stationary rigidity. For a strictly
	minimizing cone, finite weighted distance gives an asymptotic rate at
	most $-(m-2)/2$ for any nontrivial stationary limit. The stationary
	rate theorem instead gives a rate at least $\gamma_0^+$. Since
	$\gamma_0^+>-(m-2)/2$ by (S1), these bounds are incompatible.
	Proposition~\ref{e:strict-s3} therefore gives
	\[
	\text{strict minimization and (S1)}
	\quad\Longrightarrow\quad\textup{(S3)}.
	\]
	For quadratic cones, the product structure of the link makes the
	Jacobi spectrum explicit. Proposition~\ref{e:spectrum} verifies
	(S1) and (S2). The stationary arguments in
	Propositions~\ref{e:minimizing-s3} and \ref{e:exceptional-s3}
	prove (S3), including the nonminimizing pairs $(1,5)$ and $(5,1)$.
	These verifications prove Corollary~\ref{cor:examples}.
	
	\paragraph{(G) Quantitative convergence and rigidity.}
	Once continuation is established, the local discrete inequality
	applies on every sufficiently late interval $[s,s+2]$. The same
	concavity calculation, now summed from any such $s$, gives
	\[
	E(s)^{2(1-\beta)}\le C\left(E(s)-E(s+2)\right),\qquad
	\int_s^\infty G(N_\sigma)\dd\sigma\le CE(s)^\beta.
	\]
	The recurrence for the Gaussian area excess gives
	\[
	E(s)\le C(1+s-s_0)^{-1/(1-2\beta)},
	\qquad 0<\beta<1/2.
	\]
	This proves Theorem~\ref{thm:energy-loja}. The parameter-speed
	estimate in Lemma~\ref{lem:gauge} and the limit $a(s)\to0$ also give
	\[
	|a(s)|\le\int_s^\infty|a'(\sigma)|\dd\sigma
	\le C\int_s^\infty G(N_\sigma)\dd\sigma
	\le CE(s)^\beta.
	\]
	These quantitative estimates are proved together with continuation
	in Proposition~\ref{prop:continuation}.
	
	For a sequence of smooth self-shrinkers converging locally to the
	cone, the shrinker equation and Gaussian concentration give the
	global Gaussian and area-ratio bounds needed for the local distance
	estimate. A constant rescaled flow has $e=d$ and $\delta=0$.
	The residual contraction therefore gives
	\[
	d^2+|q(a)|^2\le\eta d^2,\qquad 0<\eta<1,
	\quad\Longrightarrow\quad d=0,\quad q(a)=0.
	\]
	The hypersurface must coincide with its fitted cone, contradicting
	smoothness at the vertex when the cone is nonflat. It proves Corollary~\ref{cor:shrinker-rigidity}.
	
	\paragraph{Organization.}
	Section~\ref{sec:preliminaries} develops the Gaussian first-variation,
	link, spectral, and annular estimates. Section~\ref{sec:compactness}
	treats normalized compactness, the vertex estimates, and contraction
	under (S3). Section~\ref{sec:main-proofs} proves the main conclusions:
	Subsection~\ref{sec:discrete} compares the truncated and unbounded
	Gaussian distances and proves the local discrete inequality;
	Subsection~\ref{sec:uniqueness} establishes continuation,
	fixed-coordinate uniqueness, and the quantitative convergence estimates.
	The final subsection proves rigidity among smooth self-shrinkers.
	Section~\ref{sec:examples} verifies the hypotheses for the examples.
	
	Theorems~\ref{thm:main} and \ref{thm:energy-loja} follow from
	Proposition~\ref{prop:continuation} in Subsection~\ref{sec:uniqueness}.
	Corollary~\ref{cor:examples} is proved in Section~\ref{sec:examples},
	after the spectral and stationary verifications. Corollary~\ref{cor:shrinker-rigidity} is proved in the final
	subsection of Section~\ref{sec:main-proofs}.
	
	\paragraph{Acknowledgement.}
	Hanbing Fang is supported by Hong Kong RGC grants No.~14300623 and No.~14304225. Yongheng Han is supported by the Project of Stable Support for Youth Team in Basic Research Field, Chinese Academy of Sciences (YSBR-001).
	
	\paragraph{AI disclosure.}
	This project began with the study of uniqueness of tangent flows modeled on
	generalized Simons cones. That work
	used parabolic differential inequalities for the square of the quadratic
	defining function $q|x|^2-p|y|^2$ and was developed with AI assistance.
	We subsequently used AI assistance to investigate general minimal cones
	and the corresponding Lojasiewicz inequalities and rigidity results in
	this paper. The principal AI tools were ChatGPT-6 Astra and the Danus
	system. The authors take full responsibility for the content of this
	article, including its mathematical claims, arguments, and conclusions.
	
	\section{Preliminary estimates}\label{sec:preliminaries}
	Throughout the paper, $B_R(x)$ denotes the Euclidean ball and $B_R:=B_R(0)$;
	$r:=|X|$ is the ambient radius. The induced area measure on a hypersurface
	$N$ is $\dd\mu_N$, and $\dd\|V\|$ is the weight measure of a varifold $V$. When the hypersurface has no ambiguity, we also use $\dd\mu$ for simplicity.
	Unless a different measure is displayed, on a hypersurface $N$ we write
	\[
	\langle v,w\rangle_{L^2_\rho(N)}:=\int_Nvw\rho\dd\mu_N,
	\qquad \|v\|_2=\|v\|_{L^2_\rho(N)}
	:=\left(\int_N|v|^2\rho\dd\mu_N\right)^{1/2}.
	\]
	For fields, $vw$ is replaced by their metric inner product. The underlying
	hypersurface is the one on which the functions are defined. In Euclidean
	reference hypersurface estimates we explicitly use the unweighted $L^2$ norm.
	Gradients and norms use the indicated hypersurface metric. The letters $c$ and $C$ denote positive constants that may change from line
	to line; we take $0<c\le1\le C$. Unless specified otherwise, these
	constants depend only on the fixed cone, the area-ratio
	bound $\Lambda$, and the fixed cutoffs and parameter neighborhood.
	The symbol $D$ denotes the Euclidean derivative for ambient functions and the
	Fr\'echet derivative for maps between function spaces; $D_a$ specifies
	parameter differentiation. Tangential derivatives are denoted by $\nabla$.
	We reserve $t$ for original flow time, $s$ (and $\sigma$ as an integration
	variable) for rescaled time, $r=|X|$ for spatial radius, and
	$\upsilon>0$ for Gaussian variance (defined in \eqref{eq:gaussian-variance}).
	Radial cutoffs are denoted by $\zeta$ with subscripts as needed, angular
	cutoffs by $\xi$, and compact observation densities by $\chi$.
	We write $f\asymp g$ if $C^{-1}g\le f\le Cg$ for such a uniform $C$.
	All limits and little-oh estimates are taken along the sequence under
	consideration, with any fixed auxiliary parameters held fixed. We always assume the volume upper bound \eqref{equ:areabd}.

	In this section, we establish the estimates used to fit the evolving
	hypersurface to the family $\C_a$ constructed in
	Proposition~\ref{prop:link-reduction}. Its links are denoted by
	$\Sigma_a$, and $a(s)$ is the fitted parameter along the flow.
	The reference cones need not be minimal; their reduced-link residual
	is denoted by $q(a)=\nabla f(a)$. The Gaussian first-variation estimates (Lemmas~\ref{lem:transport} and \ref{lem:exterior-tail}) control changes
	in positional distance; the link analysis (Proposition~\ref{prop:spectrum}) identifies the full kernel,
	and an annular observation (Lemma~\ref{lem:gauge}) fixes the family parameter. The final area
	estimate (Proposition~\ref{prop:anchorarea}) converts distance and shrinker quantity into Gaussian area excess.

	\subsection{Positional distance and Gaussian first-variation estimates}\label{sec:transport}
	The area-ratio bound controls all Gaussian polynomial moments uniformly
	in time. In the integrations below, we use compactly supported cutoffs
	and let their radii tend to infinity. The moment bounds justify this
	limit, together with Cauchy--Schwarz for the terms involving $\Phi$.

	\begin{lem}\label{lem:transport}
		Let $(N_s)_{s_1\le s\le s_2}$ be a smooth proper boundaryless rescaled
		mean curvature flow with the uniform area-ratio bound of Theorem~\ref{thm:main},
		and suppose $s_1\le s_2$. If $U$ is globally Lipschitz, $U(0)=0$, and $\operatorname{Lip}(U)\le L$, then we have
		\[
		\left(\int_{N_{s_2}}U^2\rho\dd\mu\right)^{1/2}
		\le\left(\int_{N_{s_1}}U^2\rho\dd\mu\right)^{1/2}+L\int_{s_1}^{s_2}G(N_{\sigma})\dd\sigma.
		\]
		In particular, $D_a(N_{s_2})\le D_a(N_{s_1} )+\sqrt{s_2-s_1}\,
		\left(\int_{s_1}^{s_2}G(N_s)^2\dd s\right)^{1/2}$ holds on $[s_1,s_2]$.
	\end{lem}
	\begin{proof}
		For smooth $U$, by the standard first variation formula, we have
		\[
		\frac{\dd}{\dd s}\int_{N_s}U^2\rho\dd\mu
		=\int_{N_s}(2U\langle DU,\Phi\rangle-U^2|\Phi|^2)\rho\dd\mu
		\le2L\left(\int_{N_s}U^2\rho\dd\mu\right)^{1/2}G(N_s).
		\]
		Writing $I(s)=\int_{N_s}U^2\rho\dd\mu$, for every $\epsilon>0$ we have
		\[
		\frac{\mathrm d}{\mathrm ds}\sqrt{I(s)+\epsilon}
		=\frac{I'(s)}{2\sqrt{I(s)+\epsilon}}
		\le L\frac{\sqrt{I(s)}}{\sqrt{I(s)+\epsilon}}G(N_s)
		\le LG(N_s).
		\]
		By integrating and then letting $\epsilon\downarrow0$, we obtain the assertion.
		For Lipschitz $U$, smooth approximation with uniformly bounded Lipschitz
		constants and $|U(X)|\le L|X|$ permits passage to the limit by the Gaussian
		moment bound and dominated convergence. This proves the lemma.
	\end{proof}
	
	\begin{lem}\label{lem:comparison}
		For a small smooth family of cones and hypersurfaces with the fixed area
		bound, we have 
		\[|D_a(N)-D_b(N)|\le C|a-b|.\]
	\end{lem}
	\begin{proof}
		By the smooth dependence of the compact links, it follows that
		\begin{align*}
			|\dist(X,\C_a)-\dist(X,\C_b)|\le C|a-b||X|.
		\end{align*}
		By the reverse triangle inequality and \eqref{equ:areabd}, we obtain
		\[
		\begin{aligned}
			|D_a(N)-D_b(N)|
			&\le\|\dist(\cdot,\C_a)-\dist(\cdot,\C_b)\|_{L^2_\rho(N)}\\
			&\le C|a-b|\left(\int_N|X|^2\rho\dd\mu\right)^{1/2}
			\le C_\Lambda|a-b|.
		\end{aligned}
		\]
		Here the last moment is bounded by summing
		$C\Lambda(k+1)^{m+2}e^{-k^2/4}$ over the annuli $k\le|X|<k+1$.
		This finishes the proof.
	\end{proof}

	To estimate distance by smooth ambient functions, we first construct a
	homogeneous defining function for the cone.
	\begin{lem}\label{lem:positional-function}
		Let \(\mathcal C\subset\mathbb R^{m+1}\), \(m\ge2\), be a two-sided embedded minimal hypersurface cone with smooth, compact, connected link. Then there is a function \(q:\mathbb R^{m+1}\setminus\{0\}\to\mathbb R\), smooth and homogeneous of degree two, whose zero set is \(\mathcal C\setminus\{0\}\), and constants \(c,C>0\), depending only on \(\mathcal C\), such that we have
		\[
		cr\operatorname{dist}(X,\mathcal C)\le |q(X)|\le Cr\operatorname{dist}(X,\mathcal C),\qquad\mathrm{where}\qquad r=|X|.
		\]
		For every \(m\)-dimensional hyperplane \(P\subset\mathbb R^{m+1}\) and every \(X\ne0\), we have
		\begin{equation}\label{lem:positional-function:1}
			-\operatorname{tr}_{P}D^2(q^2)(X)+2q(X)^2
			\le \left(2+Cr^{-2}\right)q(X)^2-|\operatorname{proj}_{P}Dq(X)|^2. 
		\end{equation}
		Consequently \(U(X):=q(X)/|X|\), extended by \(U(0)=0\), is globally Lipschitz, homogeneous of degree one, and satisfies \(c\operatorname{dist}(X,\mathcal C)\le|U(X)|\le C\operatorname{dist}(X,\mathcal C)\).
	\end{lem}

	\begin{proof}
		Since the two-sided closed embedded link \(\Sigma=\mathcal C\cap S^m\) separates the sphere, we can choose a smooth signed defining function \(f:S^m\to\mathbb R\) which agrees with signed spherical distance in a neighborhood of \(\Sigma\) and has no other zeros. We define \(q(r\theta):=r^2f(\theta)\). Compactness of the sphere, transversality along \(\Sigma\), and homogeneity give the asserted two-sided comparison with \(r\operatorname{dist}(X,\mathcal C)\).
		
		By homogeneity, it suffices to prove \eqref{lem:positional-function:1} on \(r=1\). For a hyperplane \(P\), direct differentiation gives
		\begin{equation}\label{lem:positional-function:2}
			-\operatorname{tr}_{P}D^2(q^2)
			=-2q\operatorname{tr}_{P}D^2q-2|\operatorname{proj}_{P}Dq|^2. 
		\end{equation}
		In a fixed neighborhood of \(\Sigma\), \(|Dq|\ge c_1>0\). Put \(n_q:=Dq/|Dq|\), choose a unit normal \(n_P\) to \(P\) with \(n_P\cdot n_q\ge0\), and set \(\tau^2:=1-(n_P\cdot n_q)^2\). Then we get \(|n_P-n_q|\le\sqrt2\tau\) and
		\begin{equation}\label{lem:positional-function:3}
			|\operatorname{proj}_{P}Dq|^2=|Dq|^2\tau^2. 
		\end{equation}
		Define \(T_0:=\Delta q-D^2q(n_q,n_q)\). On the level set \(q=0\), this is \(|Dq|\) times its scalar mean curvature, up to the fixed sign convention, and therefore vanishes because \(\mathcal C\setminus\{0\}\) is minimal. Smoothness and transversality give \(|T_0|\le C|q|\) near \(\Sigma\). Since \(D^2q\) is bounded on the unit sphere, it follows that
		\begin{equation}\label{lem:positional-function:4}
			|\operatorname{tr}_{P}D^2q|=|\Delta q-D^2q(n_P,n_P)|\le C|q|+C\tau. 
		\end{equation}
		Inserting \eqref{lem:positional-function:3}--\eqref{lem:positional-function:4} into \eqref{lem:positional-function:2}, we bound its right-hand side by
		\[
		Cq^2+C|q|\tau-2|Dq|^2\tau^2.
		\]
		Since \(|Dq|\ge c_1\), Young's inequality gives \(C|q|\tau\le |Dq|^2\tau^2+C'q^2\). After enlarging \(C\), we bound the preceding expression by \(Cq^2-|Dq|^2\tau^2=Cq^2-|\operatorname{proj}_{P}Dq|^2\). On the compact complement of this neighborhood, \(|q|\) is bounded below and all derivatives are bounded, so we obtain the same estimate follows after increasing the constant. This proves \eqref{lem:positional-function:1}.
		
		Finally, \(U=q/r\) is degree one and smooth away from zero. Its spherical value and spherical derivative are bounded, so \(|DU|\) is bounded above on the unit sphere and homogeneity makes the same bound valid at every radius. Since \(|U(X)|\le C|X|\), the extension \(U(0)=0\) is globally Lipschitz. Dividing the comparison for \(q\) by \(r\), we obtain the final distance comparison. Thus we conclude the proof.
	\end{proof}

	The next lemma bounds the normalized squared-distance integrals over
	$\{r>R\}$ uniformly as $R\to\infty$. Its moment estimate is related to the Gaussian
	Poincar\'e inequality of Colding--Minicozzi \cite[Section~3.1, equation~(3.3)]{CM}.
	Distance non-concentration along rescaled flow is established for cylinders
	by Sun--Wang--Xue \cite[Section~3, Corollary~3.3]{SWX} and for general convex
	cylindrical flows in \cite[Lemma~3.3]{SWXRelative}. We prove the version
	needed here for a singular minimal cone.
	
	\begin{lem}\label{lem:exterior-tail}
		Let \(\mathcal C\subset\mathbb R^{m+1}\), \(m\ge2\), be as in Lemma~\ref{lem:positional-function}, and let \(U=q/|X|\) be the globally Lipschitz degree-one positional function from Lemma~\ref{lem:positional-function}. Let \((N_s)_{0\le s\le1}\) be a smooth proper rescaled mean-curvature flow with a finite uniform area-ratio bound. Define
		\[
		I(s):=\int_{N_s}U^2\rho\dd\mu,
		\]
		\[
		\delta^2:=\int_0^1G(N_s)^2\dd s.
		\]
		Then there is \(C_{\mathcal C}<\infty\) such that we have
		\begin{equation}\label{lem:exterior-tail:1}
			\int_0^1\int_{N_s\cap\{|X|>1\}}|X|^2U^2\rho\dd\mu\dd s
			\le C_{\mathcal C}(I(0)+\delta^2), 
		\end{equation}
		and, for every \(R\ge2\),
		\begin{equation}\label{lem:exterior-tail:2}
			\int_{N_1\cap\{|X|>R\}}U^2\rho\dd\mu
			\le \frac{C_{\mathcal C}}{R^2}(I(0)+\delta^2). 
		\end{equation}
	\end{lem}
	
	\begin{proof}
		Let \(W:=U^2=q^2/r^2\), \(Q:=q^2\) and \(h:=r^2\), where \(r=|X|\), and set
		\[
		P:=\partial_s-\Delta_{N_s}+\tfrac12X^\top\cdot\nabla_{N_s}.
		\]
		By \eqref{lem:positional-function:1},
		$PQ\le(2+Cr^{-2})Q-|\nabla q|^2$ and $Ph=h-2m$.
		The quotient rule gives
		\[
		\begin{aligned}
			P(Q/h)
			&=h^{-1}PQ-Qh^{-2}Ph
			+2h^{-2}\nabla Q\cdot\nabla h-2Qh^{-3}|\nabla h|^2\\
			&=h^{-1}PQ-Qh^{-2}Ph+2h^{-1}\nabla h\cdot\nabla W\\
			&\le\left(1+\frac{C+2m}{r^2}\right)W
			+4\frac{X^\top}{r^2}\cdot\nabla W-r^{-2}|\nabla q|^2\\
			&\le\left(1+\frac{C+2m+1}{r^2}\right)W
			+4\frac{X^\top}{r^2}\cdot\nabla W-\tfrac12|\nabla U|^2.
		\end{aligned}
		\]
		Thus we have
		\begin{equation}\label{lem:exterior-tail:3}
			PW\le\left(1+\frac{C+2m+1}{r^2}\right)W+4\mathbf b\cdot\nabla W-\tfrac12|\nabla U|^2,\qquad\mathrm{where} \qquad \mathbf b:=X^\top/r^2. 
		\end{equation}
		Here one uses \(r^{-2}|\nabla q|^2\ge\tfrac12|\nabla U|^2-W/r^2\), which follows from \(\nabla q=r\nabla U+U\nabla r\) and \(|v+w|^2\ge\tfrac12|v|^2-|w|^2\).
		
		For a tangential field $Z$, we define
		$\operatorname{div}_{N_s,\rho}Z:=\rho^{-1}\operatorname{div}_{N_s}(\rho Z)
		=\operatorname{div}_{N_s}Z-\tfrac12\langle X^\top,Z\rangle$.
		Using $\operatorname{div}_{N_s}X^\top=m+\langle\mathbf H,X^\perp\rangle$
		and $\mathbf H=\Phi-X^\perp/2$, we obtain
		\[
		\begin{aligned}
			\operatorname{div}_{N_s,\rho}(r^{-2}X^\top)
			&=r^{-2}(m+\langle\mathbf H,X^\perp\rangle)
			-2r^{-4}|X^\top|^2-\tfrac12r^{-2}|X^\top|^2\\
			&=\frac m{r^2}-\frac{2(r^2-|X^\perp|^2)}{r^4}
			+\frac{\langle\Phi,X^\perp\rangle}{r^2}-\frac12.
		\end{aligned}
		\]
		Set $\mathbf Z:=\Phi+2X^\perp/r^2$. The preceding calculation gives
		\begin{equation}\label{lem:exterior-tail:4}
			\operatorname{div}_{N_s,\rho}\mathbf b
			=\frac{m-2}{r^2}+\frac{2|X^\perp|^2}{r^4}-\frac12+\frac{\langle\Phi,X^\perp\rangle}{r^2}. 
		\end{equation}
		Let $\chi=\zeta^2$ be a smooth ambient cutoff supported away from zero.
		The Gaussian first-variation identity and weighted integration by parts give
		\[
		\begin{aligned}
			\frac{\mathrm d}{\mathrm ds}\int_{N_s}\chi W\rho\dd\mu
			&=\int_{N_s}\left(\chi PW+W(\partial_s\chi+
			\Delta_{N_s}\chi-\tfrac12X^\top\cdot\nabla\chi)
			-\chi W|\Phi|^2\right)\rho\dd\mu,\\
			\int_{N_s}4\chi\mathbf b\cdot\nabla W\rho\dd\mu
			&=-4\int_{N_s}W(\mathbf b\cdot\nabla\chi+
			\chi\operatorname{div}_{N_s,\rho}\mathbf b)\rho\dd\mu,\\
			|\Phi|^2+4r^{-2}\langle\Phi,X^\perp\rangle+8r^{-4}|X^\perp|^2
			&=|\mathbf Z|^2+4r^{-4}|X^\perp|^2.
		\end{aligned}
		\]
		Substituting \eqref{lem:exterior-tail:3} and
		\eqref{lem:exterior-tail:4} into these identities, we obtain
		\begin{equation}\label{lem:exterior-tail:5}
			\begin{aligned}
				\frac{\mathrm d}{\mathrm ds}\int_{N_s}\chi W\rho\dd\mu
				&+\frac12\int_{N_s}\chi|\nabla U|^2\rho\dd\mu
				+\int_{N_s}\chi W\left(|\mathbf Z|^2+\frac{4|X^\perp|^2}{r^4}\right)\rho\dd\mu\\
				&\le\int_{N_s}\chi\left(3+\frac Br^2\right)W\rho\dd\mu
				+\int_{N_s}W(\mathcal R_0\chi+2\mathbf Z\cdot D\chi)\rho\dd\mu, 
			\end{aligned}
		\end{equation}
		where \(B:=C-2m+9\) and
		\[
		\mathcal R_0\chi=\partial_s^{\mathrm{amb}}\chi+\operatorname{tr}_{TN_s}D^2\chi-(X/2+4X/r^2)\cdot D\chi.
		\]
		Here $\partial_s^{\mathrm{amb}}\chi$ means differentiation at fixed
		$X$, so $\partial_s\chi=\partial_s^{\mathrm{amb}}\chi+D\chi\cdot\Phi$;
		it vanishes for all time-independent cutoffs used below. In rewriting the
		cutoff terms we used $\Delta_{N_s}\chi=\operatorname{tr}_{TN_s}D^2\chi+
		\mathbf H\cdot D\chi$.
		Every term in \eqref{lem:exterior-tail:5} follows by expanding the preceding displayed formulas; in particular
		\begin{equation}\label{lem:exterior-tail:6}
			2W \mathbf Z\cdot D\chi=4W\zeta \mathbf Z\cdot D\zeta\le\tfrac12\chi W|\mathbf Z|^2+8W|D\zeta|^2. 
		\end{equation}
		
		Choose a fixed radial function \(\zeta_0 \), zero on \(r\le1/2\) and one on \(r\ge1\), and apply \eqref{lem:exterior-tail:5} with \(\zeta=\zeta_0\), initially with an additional outer cutoff. The coefficients \(r^{-2}\), \(D\zeta_0\), and \(D^2\zeta_0 \) are bounded on the support of \(\zeta_0 \), and all fixed-transition errors are bounded by \(C\int_{N_s\cap\{1/2<r<1\}}W\rho\dd\mu\). In particular, we have
		\[
		|\mathcal R_0(\zeta_0^2)|+|D\zeta_0|^2
		\le C\mathbf1_{\{1/2<r<1\}},\qquad
		\zeta_0^2r^{-2}\le4\zeta_0^2.
		\]
		Lemma~\ref{lem:transport} applied to $U$ gives
		\[
		\sup_{0\le s\le1}\sqrt{I(s)}
		\le\sqrt{I(0)}+C\int_0^1G(N_s)\dd s
		\le\sqrt{I(0)}+C\delta.
		\]
		Thus we obtain
		\begin{equation}\label{lem:exterior-tail:7}
			\sup_{0\le s\le1}I(s)\le C(I(0)+\delta^2). 
		\end{equation}
		After integrating \eqref{lem:exterior-tail:5} in time,
		we drop the nonnegative terminal cutoff integral and
		apply \eqref{lem:exterior-tail:7}. Outer-cutoff errors vanish by the area-ratio bound, Gaussian tails, the linear growth of \(U\), and the nonnegative square term. We thus obtain
		from \eqref{lem:exterior-tail:5}--\eqref{lem:exterior-tail:7} that
		\begin{align}\label{lem:exterior-tail:8}
			\tfrac12\int_0^1\!\int_{N_s}\zeta_0^2
			(|\nabla U|^2+W|\mathbf Z|^2)\rho\dd\mu\dd s
			&\le \int_{N_0}\zeta_0^2W\rho\dd\mu
			+C\int_0^1 I(s)\dd s\nonumber\\
			&\le I(0)+C\sup_{[0,1]}I(s)
			\le C(I(0)+\delta^2).
		\end{align}

		Let	\(u:=\zeta_0 U\). Using $\nabla\rho=-X^\top\rho/2$ and
		$\operatorname{div}_{N_s}X^\top=m+\langle\Phi,X^\perp\rangle-|X^\perp|^2/2$,
		we have
		\[
		0=\int_{N_s}\operatorname{div}_{N_s}(X^\top u^2\rho)\dd\mu
		=\int_{N_s}\left((m-r^2/2+\langle\Phi,X^\perp\rangle)u^2
		+2uX^\top\cdot\nabla u\right)\rho\dd\mu.
		\]
		Rearranging, we obtain
		\begin{equation}
			\frac12\int_{N_s}r^2u^2\rho\dd\mu
			=m\int_{N_s}u^2\rho\dd\mu+\int_{N_s}u^2\langle\Phi,X^\perp\rangle\rho\dd\mu+2\int_{N_s}uX^\top\cdot\nabla u\rho\dd\mu. 
		\end{equation}
		Young's inequality absorbs fixed fractions of the left side in the last two terms and yields
		\begin{equation}\label{lem:exterior-tail:10}
			\int_{N_s}r^2u^2\rho\dd\mu\le C\int_{N_s}(u^2+|\nabla u|^2+|\Phi|^2u^2)\rho\dd\mu. 
		\end{equation}
		On the support of \(\zeta_0 \), we have
		\[
		|\Phi|^2\le2|\mathbf Z|^2+8|X^\perp|^2/r^4\le2|\mathbf Z|^2+32.
		\]
		Moreover, we know \(|\nabla(\zeta_0 U)|^2\le2\zeta_0^2|\nabla U|^2+2U^2|D\zeta_0|^2\), and the second term is supported in the fixed transition annulus. By combining \eqref{lem:exterior-tail:7} and \eqref{lem:exterior-tail:8}
		with the time integral of \eqref{lem:exterior-tail:10}, we obtain
		\[
		\begin{aligned}
			\int_0^1\!\int_{N_s\cap\{r>1\}}r^2U^2\rho\dd\mu\dd s
			&\le C\int_0^1\!\int_{N_s}
			\left(\zeta_0^2|\nabla U|^2+\zeta_0^2W|\mathbf Z|^2+(1+|D\zeta_0|^2)W\right)\rho\dd\mu\dd s\\
			&\le C(I(0)+\delta^2)+C\sup_{[0,1]}I(s)
			\le C(I(0)+\delta^2),
		\end{aligned}
		\]
		which proves \eqref{lem:exterior-tail:1}.
		
		For \eqref{lem:exterior-tail:2}, choose a radial cutoff \(\zeta_R\) which is zero on \(r\le R/2\), one on \(r\ge R\), and satisfies \(|D\zeta_R|\le C/R\), \(|D^2\zeta_R|\le C/R^2\). Let \(I_R(s):=\int_{N_s}\zeta_R^2W\rho\dd\mu\). Then we can apply \eqref{lem:exterior-tail:5}, absorb the \(\mathbf Z\)-cutoff term using \eqref{lem:exterior-tail:6}, and drop all remaining nonnegative terms. The cutoff coefficients on \(r>R/2\) are bounded uniformly for \(R\ge2\), so we get
		\begin{equation}\label{lem:exterior-tail:11}
			I_R'(s)\le CI_R(s)+C\int_{N_s\cap\{r>R/2\}}W\rho\dd\mu. 
		\end{equation}
		Now we integrate \eqref{lem:exterior-tail:11} from any \(s\in[1/2,3/4]\) to time one, apply Gronwall's inequality, and average in \(s\). Since \(r^2\ge R^2/4\) on \(r>R/2\), \eqref{lem:exterior-tail:1} thus gives
		\[
		\begin{aligned}
			I_R(1)&\le 4\int_{1/2}^{3/4}\left[e^{C(1-\sigma)}I_R(\sigma)
			+C\int_\sigma^1 e^{C(1-\tau)}\int_{N_\tau\cap\{r>R/2\}}W\rho\dd\mu\dd\tau\right]\dd\sigma\\
			&\le C\int_{1/2}^1\int_{N_s\cap\{r>R/2\}}W\rho\dd\mu\dd s
			\le \frac C{R^2}(I(0)+\delta^2).
		\end{aligned}
		\]
		Because \(\zeta_R=1\) on \(r>R\), this is \eqref{lem:exterior-tail:2}. The proof is complete.

	\end{proof}

	\subsection{Analytic reduction of the link and the Gaussian Jacobi operator}
	We first construct the analytic reduction of the spherical area. The
	reference links need not be minimal; their reduced gradient is retained
	in the subsequent estimates. We then identify the Gaussian kernel and
	establish an inverse which
	retains the reduced-link residual. The form-domain condition at the
	vertex will be justified in Section~\ref{sec:compactness}.
	We use spherical normal coordinates
	$\mathcal P_u:=(\omega+u n_\Sigma)/\sqrt{1+u^2}$ and let
	$\mathcal A(u):=\operatorname{Area}(\Sigma_u)$, where
	$u\in C^{2,\alpha}(\Sigma)$ is a small scalar graph height and $\omega\in\Sigma$
	is the base point. Here $\Sigma=\C\cap S^m$ is the fixed minimal link and
	$K=\ker L_\Sigma$. We also regard $A_\Sigma$ as the associated shape
	endomorphism, with one index raised. If
	$\nabla_Xn_\Sigma=-A_\Sigma X$, the gradient of $\mathcal A$ with respect to the fixed $L^2(\Sigma)$ pairing is
	\begin{equation}\label{eq:linkgradient}
		\mathcal G(u):=-H_{\Sigma_u}\det(I-uA_\Sigma)(1+u^2)^{-(m+1)/2},\qquad
		D\mathcal G(0)=-L_\Sigma.
	\end{equation}
	Here $H_{\Sigma_u}$ denotes scalar mean curvature. Denote $B_u=I-uA_\Sigma$. The area and normal-speed
	factors are
	\[
	\begin{aligned}
		J_u&=(1+u^2)^{-(m-1)/2}\det B_u
		\sqrt{1+\frac{|B_u^{-1}\nabla u|^2}{1+u^2}},\\
		\langle D_u\mathcal P_u[v],n_{\Sigma_u}\rangle
		&=\frac{v(1+u^2)^{-1}}
		{\sqrt{1+|B_u^{-1}\nabla u|^2/(1+u^2)}},\\
		D\mathcal A(u)[v]
		&=-\int_\Sigma H_{\Sigma_u}
		\langle D_u\mathcal P_u[v],n_{\Sigma_u}\rangle J_u\dd\mu_\Sigma
		=\int_\Sigma\mathcal G(u)v\dd\mu_\Sigma.
	\end{aligned}
	\]
	Thus we obtain \eqref{eq:linkgradient} directly from first variation; here
	$v\in C^{2,\alpha}(\Sigma)$ is an arbitrary scalar variation.

	\begin{prop}
		\label{prop:link-reduction}
		Let $\Sigma$ be a smooth closed embedded two-sided minimal hypersurface
		of $S^m$, and let $K=\ker L_\Sigma$. Fix an orthonormal basis of
		$K$ and identify $K\simeq\R^{d_K}$ for the family parameters; the
		element $a\in K$ below is the kernel component of the scalar graph
		height. Fix $0<\alpha<1$. There are a neighborhood $U\subset K$ of
		zero and an analytic map
		\[
		\eta:U\longrightarrow C^{2,\alpha}(\Sigma)\cap K^\perp,
		\qquad \eta(0)=D\eta(0)=0,
		\]
		such that, with $k(a):=a+\eta(a)$,
		\[
		P_{K^\perp}\mathcal G(k(a))=0.
		\]
		Define
		\[
		\mathcal P_a:=\mathcal P_{k(a)},\qquad
		\Sigma_a:=\mathcal P_a(\Sigma),\qquad
		\C_a:=\{r\mathcal P_a(\omega):r\ge0,\ \omega\in\Sigma\}.
		\]
		Each link $\Sigma_a$ is smooth and embedded. Denote
		$f(a):=\mathcal A(k(a))$, with gradient
		taken in the fixed $L^2(\Sigma)$ inner product on $K$. Then we have
		\[
		\nabla f(a)=\mathcal G(k(a))\in K.
		\]
		All sufficiently small minimal normal graphs are exactly the members
		with $\nabla f(a)=0$. There are $C<\infty$ and
		$0<\theta\le1/2$ such that, after shrinking $U$,
		\begin{equation}\label{eq:reduced-link-gradient}
			|f(a)-f(0)|^{1-\theta}
			\le C|\nabla f(a)|.
		\end{equation}
		More generally, every sufficiently small $u\in C^{2,\alpha}(\Sigma)$
		satisfies
		\[
		|\mathcal A(u)-\mathcal A(0)|^{1-\theta}
		\le C\|\mathcal G(u)\|_{L^2(\Sigma)}.
		\]
	\end{prop}
	\begin{proof}
		The complementary map
		\[
		(a,b)\longmapsto P_{K^\perp}\mathcal G(a+b),\qquad
		b\in C^{2,\alpha}(\Sigma)\cap K^\perp,
		\]
		is analytic. Its derivative in $b$ at zero is
		$-L_\Sigma|_{K^\perp}$, an isomorphism onto
		$C^{0,\alpha}(\Sigma)\cap K^\perp$. The analytic implicit-function
		theorem gives $\eta$. Its derivative vanishes at zero because
		$L_\Sigma a=0$ for $a\in K$. Elliptic regularity makes the links
		smooth, and smallness gives embeddedness. For $\zeta\in K$, we have
		\[
		Df(a)[\zeta]
		=\langle\mathcal G(k(a)),
		\zeta+D\eta(a)[\zeta]\rangle_{L^2}
		=\langle\mathcal G(k(a)),\zeta\rangle_{L^2}.
		\]
		This yields the gradient identity. Every sufficiently small critical graph
		satisfies the complementary equation and therefore has the stated form by
		 the implicit function theorem. The converse follows from the same
		gradient identity. The finite-dimensional analytic gradient inequality
		\cite{Lojasiewicz} gives the estimate for $f$. If
		$f$ is locally constant, \eqref{eq:reduced-link-gradient} is immediate.
		
		For the final assertion, write $a=P_Ku$ and
		$v=u-k(a)\in K^\perp$. The elliptic estimate for
		$-L_\Sigma$ on $K^\perp$ gives
		\[
		\|v\|_{H^2(\Sigma)}\le C\|L_\Sigma v\|_{L^2(\Sigma)}.
		\]
		Here $H^j$ denotes the ordinary Sobolev norm on the fixed compact link.
		The quasilinear expression for $\mathcal G$ and its linearization at
		zero imply, uniformly in a sufficiently small $C^{2,\alpha}$ ball,
		\[
		\|\mathcal G(u)-\mathcal G(k(a))+L_\Sigma v\|_{L^2}
		\le\varepsilon\|v\|_{H^2}.
		\]
		Projecting to $K^\perp$ and absorbing the last term, we therefore obtain
		\[
		\|v\|_{H^2}\le C\|\mathcal G(u)\|_{L^2},\qquad
		|\nabla f(a)|\le C\|\mathcal G(u)\|_{L^2}.
		\]
		The linear term in the area difference vanishes, since
		$\mathcal G(k(a))\in K$ and $v\perp K$.
		The second variation of the area functional has uniformly bounded
		coefficients on this neighborhood, so Taylor's formula gives
		\[
		|\mathcal A(u)-f(a)|
		\le C\|v\|_{H^1}^2
		\le C\|\mathcal G(u)\|_{L^2}^2.
		\]
		Combining this estimate with the reduced gradient inequality, we obtain
		the assertion. Indeed, $2(1-\theta)\ge1$, and
		$\|\mathcal G(u)\|_{L^2}\le1$ after shrinking the neighborhood.
		This concludes the proof.
	\end{proof}
	
	\begin{cor}\label{cor:reduced-cone-loja}
		Assume $m\ge3$, and let $\C_a$ be the conical family in
		Proposition~\ref{prop:link-reduction}. With $\theta$ as in \eqref{eq:reduced-link-gradient}, we have
		\[
		|F(\C_a)-F(\C)|^{1-\theta}
		\le C G(\C_a).
		\]
		Here the mean curvature and the Gaussian integral defining $G$ are
		taken on the regular cone.
	\end{cor}
	\begin{proof}
		Define $A_m:=\int_0^\infty r^{m-1}e^{-r^2/4}\dd r$.
		Polar integration gives
		\[
		F(\C_a)=c_m A_m f(a),\qquad
		G(\C_a)^2
		=c_m\left(\int_0^\infty r^{m-3}e^{-r^2/4}\dd r\right)
		\int_{\Sigma_a}|H_{\Sigma_a}|^2\dd\mu.
		\]
		The second identity uses that the position vector is tangent to a cone
		and its mean curvature has size $r^{-1}|H_{\Sigma_a}|$.
		The radial integral is finite for $m\ge3$. The smooth positive factor
		in \eqref{eq:linkgradient} and the link area Jacobian are uniformly
		bounded above and below. Using \eqref{eq:linkgradient}, we obtain
		\begin{equation}\label{eq:reduced-cone-gradient-comparison}
			\begin{aligned}
				|F(\C_a)-F(\C)|^{1-\theta}
				&=(c_mA_m)^{1-\theta}
				|f(a)-f(0)|^{1-\theta}
				\le C|\nabla f(a)|,\\
				C^{-1}|\nabla f(a)|^2
				&\le G(\C_a)^2
				\le C|\nabla f(a)|^2.
			\end{aligned}
		\end{equation}
		Here the first inequality uses \eqref{eq:reduced-link-gradient}, and
		the second line follows from the displayed polar formula and the
		fixed-measure gradient identity. This proves the corollary.
		This completes the proof.
	\end{proof}
	
	We henceforth write $q(a)=\nabla f(a)$. With the chosen
	normal and the fixed link parametrization, \eqref{eq:linkgradient} gives
	\begin{equation}\label{eq:reduced-cone-curvature}
		H_{\C_a}(r,\omega)
		=-\frac{q(a)(\omega)}{r\mathfrak w_a(\omega)},\qquad
		\mathfrak w_a=\det(I-k(a)A_\Sigma)
		(1+k(a)^2)^{-(m+1)/2}.
	\end{equation}
	This factor is smooth, positive, and uniformly bounded above and below
	for small $a$.
	
	We now pass from the compact link to the cone. For the ordinary Jacobi operator $J_\C:=\Delta_\C+|A_\C|^2$, we have
	\[
	J_\C(r^\gamma\phi_j)
	=r^{\gamma-2}[\gamma(\gamma+m-2)-\mu_j]\phi_j.
	\]
	For scalar functions supported away from the vertex, define the ordinary
	Jacobi quadratic form and its polarization by
	\begin{equation}
		\begin{aligned}
			Q^0(v,w)&:=\int_\C(\langle\nabla v,\nabla w\rangle-|A_\C|^2vw)\dd\mu_\C,\\
			Q^0(v)&:=Q^0(v,v)=-\int_\C vJ_\C v\dd\mu_\C.
		\end{aligned}
	\end{equation}
	For $v\in C_c^\infty(\C\setminus\{0\})$, the \textbf{radial Hardy inequality}
	and the \textbf{angular Rayleigh inequality} are, respectively,
	\begin{align}
		\int_0^\infty|\partial_r v|^2r^{m-1}\dd r
		&\ge\frac{(m-2)^2}{4}\int_0^\infty v^2r^{m-3}\dd r,
		\label{eq:radial-hardy}\\
		\int_\Sigma(|\nabla_\Sigma v|^2-|A_\Sigma|^2v^2)\dd\mu_\Sigma
		&\ge\mu_0\int_\Sigma v^2\dd\mu_\Sigma.
		\label{eq:angular-rayleigh}
	\end{align}
	Indeed the difference in \eqref{eq:radial-hardy} is
	$\int_0^\infty|v_r+(m-2)v/(2r)|^2r^{m-1}\dd r$.
	Polar integration and \eqref{eq:angular-rayleigh} now imply
	\[
	\int_\C(|\nabla v|^2-|A_\C|^2v^2)\dd\mu_\C
	\ge\nu_0^2\int_\C r^{-2}v^2\dd\mu_\C.
	\]
	Since $|A_\C|^2\le Cr^{-2}$, this also controls the full norm
	$\int_\C(|\nabla v|^2+r^{-2}v^2)\dd\mu_\C$ with the explicit constant
	\[
	Q^0(v)\ge c_H\int_\C(|\nabla v|^2+r^{-2}v^2)\dd\mu_\C,
	\qquad\mathrm{where}\qquad c_H:=\frac{\nu_0^2}{\nu_0^2+\max_\Sigma|A_\Sigma|^2+1}>0.
	\]
	The constants
	remain uniform for the small family.
	
	The Gaussian Jacobi operator is
	\[
	L_\C=J_\C-\tfrac12X\cdot\nabla_\C+\tfrac12
	=J_\C-\tfrac r2\partial_r+\tfrac12.
	\]
	In polar coordinates this reads
	\[
	L_\C=\partial_{rr}+\left(\frac{m-1}{r}-\frac r2\right)\partial_r
	+r^{-2}(L_\Sigma-(m-1))+\frac12.
	\]
	The associated form domain is the completion
	$\V:=\overline{C_c^\infty(\C\setminus\{0\})}^{\|\cdot\|_\V}$, where
	\[
	\norm v_\V^2:=\int_\C
	(|\nabla v|^2+(1+r^{-2}+r^2)v^2)\rho\dd\mu.
	\]
	The Gaussian quadratic form, initially on $C_c^\infty(\C\setminus\{0\})$, is
	\begin{equation}\label{eq:gaussian-jacobi-form}
		Q_\rho(v):=-\langle v,L_\C v\rangle_{L^2_\rho}
		=\int_\C\left(|\nabla v|^2-|A_\C|^2v^2-\tfrac12v^2\right)\rho\dd\mu_\C.
	\end{equation}
	The \textbf{Friedrichs realization} of $-L_\C$ is associated to the closure
	of this lower-semibounded form; $L_\C^F$ denotes its negative. We now prove
	that its form norm is equivalent to $\|\cdot\|_\V$.
	
	Denote $a_H=(m-2)/2$. Completing the radial square and integrating by parts, we obtain
	\[
	\begin{aligned}
		0&\le\int_\C|v_r+a_Hv/r|^2\rho\dd\mu_\C\\
		&=\int_\C|v_r|^2\rho\dd\mu_\C
		-a_H^2\int_\C r^{-2}v^2\rho\dd\mu_\C
		+\tfrac{a_H}{2}\|v\|_2^2,
	\end{aligned}
	\]
	where we used
	$\partial_r(r^{m-2}e^{-r^2/4})
	=((m-2)r^{m-3}-r^{m-1}/2)e^{-r^2/4}$ in the cross term.
	Thus we get the \textbf{Gaussian radial Hardy inequality}
	\begin{equation}\label{eq:gaussian-radial-hardy}
		\int_\C|v_r|^2\rho\dd\mu_\C
		\ge\frac{(m-2)^2}{4}\int_\C r^{-2}v^2\rho\dd\mu_\C
		-\frac{m-2}{4}\|v\|_2^2.
	\end{equation}
	Combining \eqref{eq:gaussian-jacobi-form}, the Gaussian radial Hardy
	inequality \eqref{eq:gaussian-radial-hardy}, and the angular Rayleigh
	inequality \eqref{eq:angular-rayleigh}, we obtain
	\begin{equation}\label{eq:gaussian-hardy-coercivity}
		\begin{aligned}
			Q_\rho(v)
			&=\int_\C|v_r|^2\rho\dd\mu_\C
			+c_m\int_0^\infty r^{m-3}e^{-r^2/4}
			\int_\Sigma(|\nabla_\Sigma v|^2-|A_\Sigma|^2v^2)\dd\mu_\Sigma\dd r
			-\tfrac12\|v\|_2^2\\
			&\ge\left(\frac{(m-2)^2}{4}+\mu_0\right)
			\int_\C r^{-2}v^2\rho\dd\mu_\C
			-\left(\frac{m-2}{4}+\frac12\right)\|v\|_2^2\\
			&=\nu_0^2\int_\C r^{-2}v^2\rho\dd\mu_\C-\frac m4\|v\|_2^2,
			\qquad\mathrm{where}\qquad c_m=(4\pi)^{-m/2}.
		\end{aligned}
	\end{equation}
	Denote $A_0=\max_\Sigma|A_\Sigma|^2$. Equations
	\eqref{eq:gaussian-jacobi-form} and \eqref{eq:gaussian-hardy-coercivity} give
	\begin{equation}\label{eq:gaussian-gradient-control}
		\begin{aligned}
			\int_\C|\nabla v|^2\rho\dd\mu_\C
			&=Q_\rho(v)+\int_\C|A_\C|^2v^2\rho\dd\mu_\C+\tfrac12\|v\|_2^2\\
			&\le Q_\rho(v)+A_0\int_\C r^{-2}v^2\rho\dd\mu_\C+\tfrac12\|v\|_2^2\\
			&\le\left(1+\frac{A_0}{\nu_0^2}\right)Q_\rho(v)
			+\left(\frac{mA_0}{4\nu_0^2}+\frac12\right)\|v\|_2^2.
		\end{aligned}
	\end{equation}
	For the quadratic radial term, $X$ is tangent to $\C$, so
	\[
	\begin{aligned}
		\tfrac12\int_\C r^2v^2\rho\dd\mu_\C
		&=m\|v\|_2^2+2\int_\C vX\cdot\nabla v\rho\dd\mu_\C\\
		&\le m\|v\|_2^2+\tfrac14\int_\C r^2v^2\rho\dd\mu_\C
		+4\int_\C|\nabla v|^2\rho\dd\mu_\C.
	\end{aligned}
	\]
	Rearranging and applying \eqref{eq:gaussian-gradient-control}, we obtain
	\begin{equation}\label{eq:gaussian-moment-control}
		\int_\C r^2v^2\rho\dd\mu_\C
		\le4m\|v\|_2^2+16\int_\C|\nabla v|^2\rho\dd\mu_\C
		\le C\left(Q_\rho(v)+C\|v\|_2^2\right).
	\end{equation}
	Combining \eqref{eq:gaussian-hardy-coercivity}--\eqref{eq:gaussian-moment-control},
	and using $|A_\C|^2\le A_0r^{-2}$ for the reverse bound, we obtain, for a fixed
	sufficiently large $C_0$,
	\begin{equation}\label{eq:gaussian-form-equivalence}
		c\|v\|_\V^2\le Q_\rho(v)+C_0\|v\|_2^2\le C\|v\|_\V^2.
	\end{equation}
	Thus it follows that $\V$ is the closed form domain. The superscript $F$ below denotes
	the self-adjoint Friedrichs realization on $L^2_\rho$;
	$\ker L_\C^F=\{v\in\operatorname{Dom}(L_\C^F):L_\C^Fv=0\}$.
	\begin{prop}\label{prop:spectrum}
		Under (S1) and (S2), we have
		\[
		\ker L_\C^F=\{r\psi:\psi\in K\}.
		\]
		For every $v$ in the operator domain perpendicular to this kernel, we have
		\[
		\|v\|_{L^2_\rho}+\|v\|_\V\le C\|L_\C^Fv\|_{L^2_\rho}.
		\]
	\end{prop}
	\begin{proof}
		Separating variables and setting $z=r^2/4$, we reduce the radial equation
		to the generalized Laguerre equation. More precisely, if
		$f(r)=r^\gamma y(z)$, $\gamma=\gamma_j^+$, then
		\[
		L_\C(f\phi_j)=r^\gamma\phi_j
		\left(zy''+(1+\nu_j-z)y'+\frac{1-\gamma}{2}y\right).
		\]
		The Laguerre polynomials are complete for $z^{\nu_j}e^{-z}\dd z$;
		see \cite[Chapter~V, Sections~5.1 and 5.7]{Szego}.
		This gives the complete
		Friedrichs eigenfunctions and eigenvalues
		\[
		r^{\gamma_j^+}L_k^{\nu_j}(r^2/4)\phi_j,\qquad
		\lambda_{j,k}=\frac{1-\gamma_j^+}2-k.
		\]
		Here $L_k^\nu$ is the generalized Laguerre polynomial. The plus root is
		selected by requiring a finite Hardy integral at zero. A zero eigenvalue requires
		$\gamma_j^+=1-2k$, with $k<m/4$. Condition (S2) removes $k\ge1$;
		$k=0$ gives exactly $rK$ at the base cone.
		Only finitely many pairs $(j,k)$ have $|\lambda_{j,k}|\le1$,
		since this inequality bounds $k$ and then bounds $\mu_j$ above. Thus
		there is a positive gap from zero on the kernel complement.
		Spectral expansion gives the
		$L^2$ estimate. Form coercivity up to a bounded $L^2$ term gives
		$\norm v_\V^2\le C(-\langle v,L^Fv\rangle+C\norm v_2^2)$.
		By \eqref{eq:gaussian-form-equivalence} and the spectral gap $g>0$, we get
		\[
		\begin{aligned}
			\|v\|_2^2&=\sum_{\lambda\ne0}|v_\lambda|^2
			\le g^{-2}\sum_{\lambda\ne0}\lambda^2|v_\lambda|^2
			=g^{-2}\|L^Fv\|_2^2,\\
			\|v\|_\V^2&\le C\left(|\langle v,L^Fv\rangle|+C_0\|v\|_2^2\right)
			\le C\left(\|v\|_2\|L^Fv\|_2+\|v\|_2^2\right)
			\le C\|L^Fv\|_2^2.
		\end{aligned}
		\]
		This proves the second estimate and finishes the proof.
	\end{proof}

	\begin{prop}
		\label{prop:gaussian-bordered}
		Assume (S1) and (S2), and choose an $L^2(\Sigma)$ orthonormal basis
		$\psi_1,\ldots,\psi_{d_K}$ of $K$. Define $g_l=r\psi_l$ and choose a
		nonnegative smooth radial cutoff $\chi$ supported in
		$\{1<r<2\}$ and positive on $\{5/4<r<7/4\}$. We define
		\[
		\mathcal O_l(v):=\int_\C\chi g_lv\rho\dd\mu_\C,\qquad
		\mathcal B\xi:=r^{-1}\sum_{l=1}^{d_K}\xi_l\psi_l.
		\]
		Let $\V^*$ denote the continuous dual of the form space $\V$, with norm
		\[
		\|\Lambda\|_{\V^*}:=
		\sup_{0\ne w\in\V}\frac{|\Lambda(w)|}{\|w\|_\V}.
		\]
		We use the $L^2_\rho$ pairing to regard $L^2_\rho$ as a subspace of
		$\V^*$. The operator $L_\C$ is interpreted weakly as the bounded map
		$\V\to\V^*$ given by
		\begin{equation}\label{eq:gaussian-weak-duality}
			\langle L_\C v,w\rangle_{\V^*,\V}
			=-\int_\C\left(\langle\nabla v,\nabla w\rangle
			-( |A_\C|^2+\tfrac12)vw\right)\rho\dd\mu_\C.
		\end{equation}
		Its boundedness follows from \eqref{eq:gaussian-form-equivalence}
		and Cauchy--Schwarz in the form norm. Then the map
		\[
		(v,\xi)\longmapsto
		\bigl(L_\C v+\mathcal B\xi,\mathcal O(v)\bigr)
		:\V\times\R^{d_K}\longrightarrow\V^*\times\R^{d_K}
		\]
		is an isomorphism. In particular, we have
		\[
		\|v\|_\V+|\xi|
		\le C\bigl(\|L_\C v+\mathcal B\xi\|_{\V^*}
		+|\mathcal O(v)|\bigr).
		\]
		The estimate persists uniformly under sufficiently small smooth
		perturbations of the link, the residual columns, and the observations,
		after radial pullback to the base cone.
	\end{prop}
	\begin{proof}
		The map $\mathcal B$ takes values in $L^2_\rho$, since $m\ge3$ and
		$\int_0^1r^{m-3}\dd r<\infty$. It is therefore bounded into $\V^*$.
		By the base-cone assertion of Proposition~\ref{prop:spectrum},
		$\mathcal K:=\ker L_\C^F=rK$. The same separated eigenbasis shows that
		the weak operator induces an isomorphism
		\[
		L_\C:\V\cap\mathcal K^\perp\longrightarrow
		\{F\in\V^*:F(g)=0\text{ for all }g\in\mathcal K\}.
		\]
		Indeed, choose $B>\sup\operatorname{Spec}(L_\C^F)$ so that
		$\sum(B-\lambda_j)|v_j|^2$ is an equivalent squared form norm.
		For every nonzero eigenvalue, the ratio
		$(B-\lambda_j)/|\lambda_j|$ is uniformly bounded: the spectrum is
		discrete, has a gap from zero away from $\mathcal K$, and tends to
		$-\infty$. Dividing the spectral coefficients by $\lambda_j$, we obtain
		the claimed bounded inverse between the form space and its dual.
		
		For given $(F,b)\in\V^*\times\R^{d_K}$, we pair the equation
		$L_\C v+\mathcal B\xi=F$ with each $g_l$. Self-adjointness gives
		\[
		\langle\mathcal B\xi,g_l\rangle
		=c_m\left(\int_0^\infty r^{m-1}e^{-r^2/4}\dd r\right)\xi_l
		=F(g_l).
		\]
		Thus it follows that $\xi$ is uniquely determined and bounded by $C\|F\|_{\V^*}$.
		The functional $F-\mathcal B\xi$ annihilates $\mathcal K$, so the
		preceding inverse gives a unique $v_\perp\perp\mathcal K$.
		Write $v=v_\perp+\sum c_lg_l$. The observation matrix is
		\[
		\mathcal O_l(g_j)
		=c_m\left(\int_0^\infty\chi(r)r^{m+1}e^{-r^2/4}\dd r\right)
		\delta_{lj},
		\]
		which is positive definite. Therefore it follows that $\mathcal O(v)=b$ fixes
		$c$ uniquely, with the asserted bound. If $K=\{0\}$, both finite
		blocks are absent and the same proof is simply the inverse for $L_\C$.
		
		For the perturbation statement, let $\Sigma_a$ be a small
		smooth family of links, and identify its cone with $\C$ by the
		radial parametrization. We pull back the weak form
		\[
		\langle\mathscr L_av,w\rangle
		=-\int_{\C_a}
		\bigl(\langle\nabla v_a,\nabla w_a\rangle
		-( |A_{\C_a}|^2+\tfrac12)v_aw_a\bigr)\rho\dd\mu,
		\]
		where $v_a,w_a$ are the transported functions. Smoothness of the link
		metrics, angular Jacobians, and second fundamental forms gives
		\[
		\|\mathscr L_a-L_\C\|_{\V\to\V^*}\longrightarrow0.
		\]
		For example, the curvature coefficient changes by $o(1)r^{-2}$,
		which is controlled by the Hardy term in the form norm. The radial
		coordinate and Gaussian are unchanged. Smooth perturbations of the
		columns $r^{-1}\psi_l$ converge in $L^2_\rho$, and compact observations
		converge in the operator norm $\V\to\R^{d_K}$.
		The full block operator therefore converges in operator norm to the
		isomorphism at zero. A Neumann series gives a uniformly bounded inverse.
		Thus we conclude the proof.
	\end{proof}
	
	\begin{lem}
		\label{lem:residual-tail}
		Let $\C_a$ be the reduced conical family, and set
		$q(a)=\nabla f(a)$. Let $(N_s)_{0\le s\le1}$ be a smooth
		proper rescaled flow with a fixed all-center area bound. Define
		\[
		d_a(s)^2:=\int_{N_s}\operatorname{dist}(X,\C_a)^2\rho\dd\mu,
		\qquad \delta^2:=\int_0^1G(N_s)^2\dd s.
		\]
		There is $C$, independent of small $a$, such that for every $R\ge2$,
		\[
		\int_{N_1\cap\{r>R\}}\operatorname{dist}(X,\C_a)^2\rho\dd\mu
		\le CR^{-2}\left(d_a(0)^2+\delta^2+|q(a)|^2\right).
		\]
	\end{lem}
	\begin{proof}
		We construct the degree-two defining function $Q_a$ by signed spherical
		distance near the link and a smooth extension on the two complementary
		components, as in Lemma~\ref{lem:positional-function}. Its comparison
		with $r\operatorname{dist}(X,\C_a)$ and all derivative bounds
		are uniform in small $a$. The function $U_a=Q_a/r$ is uniformly
		Lipschitz. Let $W_a=U_a^2$ and
		$I_a(s)=\int_{N_s}W_a\rho\dd\mu$.
		Lemma~\ref{lem:transport}, applied to $U_a$, gives
		\[
		\sup_{[0,1]}I_a(s)\le C\left(I_a(0)+\delta^2\right).
		\]
		
		We then identify the additional term in the positional differential
		inequality. On $r=1$ at the link, the trace of $D^2Q_a$ over the cone
		tangent plane is bounded by $C|H_{\Sigma_a}|\le C|q(a)|$.
		For a nearby point and a plane $P$, its absolute value is therefore
		at most $C(|q(a)|+|Q_a|+\tau)$, where
		$\tau=|\operatorname{proj}_P DQ_a|/|DQ_a|$.
		The proof of \eqref{lem:positional-function:1}, using Young's inequality
		on $|Q_a|\tau$ and $|Q_a||q(a)|$, consequently gives
		\[
		-2Q_a\operatorname{tr}_P D^2Q_a
		-2|\operatorname{proj}_P DQ_a|^2
		\le Cr^{-2}Q_a^2-|\operatorname{proj}_P DQ_a|^2
		+Cr^2|q(a)|^2.
		\]
		Away from the relative normal tube, $|Q_a|/r^2$ is bounded below,
		so the same inequality follows from homogeneity and the derivative
		bounds. Division by $r^2$ in the calculation leading to
		\eqref{lem:exterior-tail:3} adds only $C|q(a)|^2$ to its right side.
		All transport, divergence, and completed-square identities thereafter
		concern the actual flow and are unchanged.
		
		The calculations in \eqref{lem:exterior-tail:5}--\eqref{lem:exterior-tail:10}
		therefore apply with the additional source $C|q(a)|^2$.
		Its integral is bounded by $C|q(a)|^2\int_{N_s}\rho\dd\mu$,
		which is uniform by the area bound. In particular, the argument for
		\eqref{lem:exterior-tail:1} gives
		\begin{equation}\label{eq:residual-tail-moment}
			\int_0^1\!\int_{N_s\cap\{r>1\}}r^2W_a\rho\dd\mu\dd s
			\le C\left(I_a(0)+\delta^2+|q(a)|^2\right).
		\end{equation}
		Choose $\zeta_R$ as in the proof of \eqref{lem:exterior-tail:11}, and put
		$I_{a,R}(s):=\int_{N_s}\zeta_R^2W_a\rho\dd\mu$.
		The cutoff calculation in \eqref{lem:exterior-tail:11} now gives
		\begin{equation}\label{eq:residual-tail-differential}
			I_{a,R}'(s)\le CI_{a,R}(s)
			+C\int_{N_s\cap\{r>R/2\}}W_a\rho\dd\mu
			+C|q(a)|^2\int_{N_s\cap\{r>R/2\}}\rho\dd\mu.
		\end{equation}
		The area bound and Gaussian decay give
		$\int_{N_s\cap\{r>R/2\}}\rho\dd\mu\le CR^{-2}$ uniformly in $s$.
		Integrating \eqref{eq:residual-tail-differential} from
		$\sigma\in[1/2,3/4]$ to time one, applying Gronwall's inequality,
		and averaging in $\sigma$, we obtain
		\[
		\begin{aligned}
			\int_{N_1\cap\{r>R\}}W_a\rho\dd\mu
			&\le I_{a,R}(1)\\
			&\le C\int_{1/2}^1\!\int_{N_s\cap\{r>R/2\}}W_a\rho\dd\mu\dd s
			+CR^{-2}|q(a)|^2\\
			&\le CR^{-2}\int_0^1\!\int_{N_s\cap\{r>1\}}r^2W_a\rho\dd\mu\dd s
			+CR^{-2}|q(a)|^2\\
			&\le CR^{-2}\left(I_a(0)+\delta^2+|q(a)|^2\right),
		\end{aligned}
		\]
		where the last inequality is \eqref{eq:residual-tail-moment}.
		The two-sided comparison of $W_a$ with squared distance proves the claim.
		The proof is complete.
	\end{proof}
	
	\subsection{Fitting the parameter and estimating Gaussian area}
	In the static estimates of this subsection, $N$ is a smooth properly
	embedded boundaryless hypersurface with area ratios bounded by $\Lambda$
	and $G(N)<\infty$. Along a flow, the area-ratio bound is uniform,
	whereas exact dissipation gives only
	\[
	\int_{s_1}^{s_2}G(N_s)^2\dd s=E(s_1)-E(s_2)<\infty
	\]
	on compact time intervals. Thus $G(N_s)$ is finite for almost every
	$s$; no uniform pointwise bound in time is assumed. We use the static
	estimates at such times and on slices selected by averaging the
	dissipation. We first determine the family parameter by requiring
	orthogonality on a fixed regular annulus.

	\begin{lem}\label{lem:gauge}
		On the fixed regular annulus $\{1/2<|X|<3\}$, a sufficiently small
		complete single graph $N$ has a unique fitted parameter $a(N)$.
		Along a rescaled flow in this graph neighborhood, it satisfies
		\[
		|a'(s)|\le CG(N_{s}),\qquad
		|a(s_2)-a(s_1 )|\le C\left(\int_{s_1}^{s_2}G(N_\sigma)^2\dd\sigma\right)^{1/2}\quad(s_2-s_1 \le1).
		\]
	\end{lem}

	\begin{proof}
		Let $\mathcal P_a:\Sigma\to\Sigma_a$ denote the parametrization in
		Proposition~\ref{prop:link-reduction}, with the fixed orthonormal
		basis of $K$ used for the coordinates $a_l$. We define
		\[
		\psi_{a,l}(\mathcal P_a(\omega))
		:=\langle\partial_{a_l}\mathcal P_a(\omega),n_{\Sigma_a}(\mathcal P_a(\omega))\rangle,
		\qquad 1\le l\le d_K.
		\]
		These are the normal parameter velocities on $\Sigma_a$. Let
		$\Psi_{a,l}:=r\psi_{a,l}$. We use the radial cutoff $\chi$ fixed in
		Proposition~\ref{prop:gaussian-bordered}, supported in $\{1<r<2\}$
		and positive on $\{5/4<r<7/4\}$.
		We first write $N$ as the graph of a small function
		$u\in C^{2,\alpha}(\C\cap\{1/2<r<3\})$ over $\C$.
		For every sufficiently small parameter $a$, normal projection onto $\C_a$
		expresses the same hypersurface as a graph with height $h(u,a)$ on a
		neighborhood of $\operatorname{supp}\chi$. The regraphed height and its first derivatives depend continuously on
		$(u,a)$; the observation map below is $C^1$, and $h(0,0)=0$.
		We define the observation map
		\begin{equation}\label{eq:gauge-observation-map}
			\mathcal O_l(u,a):=\int_{\C_a}\chi\Psi_{a,l}h(u,a)\rho\dd\mu_{\C_a},
			\qquad \mathcal O=(\mathcal O_1,\ldots,\mathcal O_{d_K}).
		\end{equation}
		Thus we determine $a(N)$ by solving the finite system
		$\mathcal O(u,a)=0$ for $a$, with $u$ fixed.
		
		To compute the parameter derivative, fix a point $Y=r\omega$ in the
		fitting annulus and put $P_a(Y)=r\mathcal P_a(\omega)$. Write
		$\widehat h(a,Y)=h(0,a)(P_a(Y))$ for the height pulled back to the base
		cone. Let $\ell$ be the signed distance to $\C$ in a tubular
		neighborhood of the fitting annulus, chosen so that
		$D\ell(Y)=n_\C(Y)$ on $\C$. Since $u=0$ represents the fixed
		hypersurface $\C$, the regraphing equation is
		\[
		\ell\bigl(P_a(Y)+\widehat h(a,Y)n_{\C_a}(P_a(Y))\bigr)=0.
		\]
		At $a=0$ and height zero, the derivative of the left-hand side
		with respect to the height is
		$\langle D\ell(Y),n_\C(Y)\rangle=1$.
		The implicit function theorem thus gives differentiability in $a$.
		Differentiating at $a=0$ with $Y$ fixed, we obtain
		\[
		\partial_{a_k}\widehat h(0,Y)
		=-\left\langle
		\left.\partial_{a_k}P_a(Y)\right|_{a=0},n_\C(Y)
		\right\rangle
		=-\Psi_{0,k}(Y).
		\]
		Thus, under the radial identification,
		$\partial_{a_k}h(0,0)=-\Psi_{0,k}$. Since $h(0,0)=0$, only the derivative of the height contributes to
		$D_a\mathcal O(0,0)$ in \eqref{eq:gauge-observation-map}.
		Thus we get
		\begin{equation}\label{eq:fitting-matrix}
			\mathcal O(0,0)=0,\qquad
			D_a\mathcal O(0,0)=-\mathsf M,\qquad
			\mathsf M_{lk}:=\int_\C\chi\Psi_{0,l}\Psi_{0,k}\rho\dd\mu_\C.
		\end{equation}
		For every $b\in\R^{d_K}$, the full rank of the family gives
		\begin{equation}\label{eq:fitting-matrix-positive}
			b^T\mathsf M b
			=\int_\C\chi r^2\left(\sum_l b_l\psi_{0,l}\right)^2\rho\dd\mu_\C
			\ge c|b|^2.
		\end{equation}
		We apply the implicit-function theorem to
		$\mathcal O:C^{2,\alpha}(\C\cap\{1/2<r<3\})\times\R^{d_K}\to\R^{d_K}$.
		There are neighborhoods of $u=0$ and $a=0$ and a unique $C^1$ map
		$\mathfrak a(u)$ such that
		\[
		\mathfrak a(0)=0,\qquad \mathcal O(u,\mathfrak a(u))=0.
		\]
		We define $a(N):=\mathfrak a(u)$. In particular, by differentiating at $u=0$, we obtain the first-order formula
		\[
		a(N)=\mathsf M^{-1}
		\left(\int_\C\chi\Psi_{0,l}u\rho\dd\mu_\C\right)_{l=1}^{d_K}
		+o(\|u\|_{C^{2,\alpha}}).
		\]
		Thus it follows that the parameter removes the component of the height measured by the
		family velocities. After shrinking the neighborhoods, we have
		$\|(D_a\mathcal O(u,a))^{-1}\|\le C$ uniformly.
		
		We now write $a(s)=a(N_s)$ and use $\mathcal O(N,a)$ for
		\eqref{eq:gauge-observation-map} in hypersurface notation.
		The derivative $D_N\mathcal O[\Phi]$ is a compactly supported linear
		integral of the normal velocity, with uniformly bounded coefficients.
		By differentiating $\mathcal O(N_s,a(s))=0$ and applying Cauchy--Schwarz,
		we obtain
		\[
		\begin{aligned}
			0&=\frac{\mathrm d}{\mathrm ds}\mathcal O(N_s,a(s))
			=D_a\mathcal O(N_s,a(s))a'(s)+D_N\mathcal O(N_s,a(s))[\Phi],\\
			|a'(s)|&\le C|D_N\mathcal O[\Phi]|
			\le C\left(\int_{N_s\cap\{1/2<r<3\}}|\Phi|^2\rho\dd\mu\right)^{1/2}
			\le CG(N_s),\\
			|a(s_2)-a(s_1)|&\le C\int_{s_1}^{s_2}G(N_s)\dd s
			\le C\sqrt{s_2-s_1}\left(\int_{s_1}^{s_2}G(N_s)^2\dd s\right)^{1/2}.
		\end{aligned}
		\]
		At zero these velocities are precisely the chosen orthonormal basis
		of $K$. The family velocities need not be Jacobi fields away from
		critical parameters; only their full rank at zero and smooth dependence
		were used. This concludes the proof.
	\end{proof}

	\begin{lem}\label{lem:radial-tilt}
		Let $N^m\subset\R^{m+1}$ be a smooth properly embedded boundaryless
		hypersurface with area ratios bounded by $\Lambda$ and $G(N)<\infty$.
		Denote $D(N):=D_0(N)$. Then we have
		\begin{equation}\label{eq:radial-tilt}
			J(N):=\int_N\frac{|X^\perp|^2}{(1+|X|^2)^2}\rho\dd\mu
			\le C_{\mathcal C}\bigl(D(N)^2+D(N)G(N)\bigr).
		\end{equation}
		
	\end{lem}
	
	\begin{proof}
		Choose $\eta>0$ sufficiently small that the normal map
		$(Y,l)\mapsto Y+ln_\C(Y)$ is one-to-one for $|l|<2\eta|Y|$.
		For $0<\tau\le2\eta$, we set
		\[
		\mathcal T_\tau:=\{Y+ln_\C(Y):Y\in\C\setminus\{0\},\ |l|<\tau|Y|\}.
		\]
		On $\mathcal T_{2\eta}$, let $l$ be signed distance to $\C$ and put
		$n_l=Dl$. With $r=|X|$, normal
		coordinates and minimality give
		\begin{equation}\label{eq:distance-tube-identities}
			|D^2l|\le C/r,\quad |\Delta l|\le C|l|/r^2,\quad
			D^2l(n_l,\cdot)=0,\quad X\cdot n_l=l.
		\end{equation}
		To verify \eqref{eq:distance-tube-identities}, write
		$X=Y+ln_\C(Y)$ and let $\kappa_j$ be the eigenvalues of $A_\C(Y)$.
		Since $\sum_j\kappa_j=0$, $|\kappa_j|\le C/|Y|$, and
		$|Y|^2+l^2=r^2$, we have $|Y|\asymp r$ and
		\[
		\begin{aligned}
			Dl&=n_\C(Y),\qquad |Dl|=1,\qquad
			D^2l(n_l,\cdot)=0,\qquad X\cdot n_l=l,\\
			|D^2l|&\le\max_j\frac{|\kappa_j|}{|1-l\kappa_j|}\le Cr^{-1},\\
			|\Delta l|&=\left|\sum_j\frac{\kappa_j}{1-l\kappa_j}\right|
			=\left|l\sum_j\frac{\kappa_j^2}{1-l\kappa_j}\right|
			\le C|l|r^{-2}.
		\end{aligned}
		\]
		The common sign in the Hessian eigenvalues depends on orientation and
		has no effect on these absolute-value estimates.
		
		Set $T:=|\nabla_Nl|^2=1-(n_N\cdot n_l)^2$, using a local normal $n_N$
		on $N$.
		For a smooth cutoff $\zeta$ compactly supported in the tube away from
		zero, we apply Gaussian first variation to $\zeta^2ln_l$. Thus we obtain
		\begin{align*}
			\int_N\bigl(&\zeta^2T+\zeta^2l\Delta l-\zeta^2lD^2l(n_N,n_N)
			+2\zeta l\nabla_N\zeta\cdot\nabla_Nl
			-\tfrac12\zeta^2l^2\bigr)\rho\dd\mu\\
			&=-\int_N\zeta^2l\,\Phi\cdot n_l\,\rho\dd\mu.
		\end{align*}
		The third term has absolute value at most $C\eta\zeta^2T$ if the tube
		has width $|l|<\eta r$, because the Hessian annihilates $n_l$.
		We choose $\eta$ so that this term is absorbed. Young's inequality absorbs
		the cutoff cross term, and \eqref{eq:distance-tube-identities} yields
		\begin{equation}\label{eq:local-tilt-test}
			\int_N\zeta^2T\rho\dd\mu
			\le C\int_N\left(\zeta^2|l||\Phi|
			+\bigl[\zeta^2(1+r^{-2})+|D\zeta|^2\bigr]l^2\right)\rho\dd\mu.
		\end{equation}
		
		Choose a smooth degree-zero angular cutoff $\xi$ with
		\[
		0\le\xi\le1,\qquad \operatorname{supp}\xi\subset\mathcal T_\eta,
		\qquad \xi=1\text{ on }\mathcal T_{\eta/2},\qquad |D\xi|\le C_\eta/r.
		\]
		Outside $\mathcal T_{\eta/2}$, compactness of the sphere gives
		$\operatorname{dist}(X,\mathcal C)\ge c_{\mathcal C}r$.
		Inside it, we write $X=Y+ln_l$ with $Y$ tangent to the cone and
		orthogonal to $n_l$. The projection of $n_N$ on the distance level
		has length $\sqrt T$, so
		$|X\cdot n_N|^2\le2r^2T+2l^2$. These observations imply
		\begin{equation}\label{eq:radial-to-normal-tilt}
			J(N)\le CD(N)^2+
			C\int_N\xi^2\frac{r^2}{(1+r^2)^2}T\rho\dd\mu.
		\end{equation}
		Now we use \eqref{eq:local-tilt-test} with
		$\zeta=\xi r/(1+r^2)$, initially with additional inner and outer
		radial cutoffs. Away from those transitions, we have
		\[
		|D\zeta|\le\frac C{1+r^2},\qquad
		\zeta^2(1+r^{-2})\le\frac C{1+r^2},\qquad \zeta^2\le C.
		\]
		Since $|l|=\operatorname{dist}(X,\mathcal C)$ in the tube,
		Cauchy--Schwarz bounds the resulting right side by $C(D^2+DG)$.
		For completeness, an inner cutoff changing on $\varepsilon<r<2\varepsilon$
		has derivative bounded by $C/\varepsilon$, which is compensated by the factor $r/(1+r^2)$.
		Its extra error is at most
		\[
		C\int_{N\cap B_{2\varepsilon}}l^2\rho\dd\mu
		\le C\Lambda\varepsilon^{m+2}\longrightarrow0.
		\]
		For an outer cutoff on $R<r<2R$, all extra errors are bounded by
		Gaussian tails with a fixed polynomial weight. These tend to zero by
		the area-ratio bound. The velocity term is integrable by
		Cauchy--Schwarz and $G<\infty$. Fatou's lemma removes the two cutoffs
		on the left, while dominated convergence applies to the remaining
		right side. Combining \eqref{eq:local-tilt-test} and
		\eqref{eq:radial-to-normal-tilt}, we obtain
		\[
		\begin{aligned}
			J(N)&\le CD(N)^2+C\int_N\zeta^2T\rho\dd\mu\\
			&\le CD(N)^2+C\int_N\left(\zeta^2|l||\Phi|
			+[\zeta^2(1+r^{-2})+|D\zeta|^2]l^2\right)\rho\dd\mu\\
			&\le CD(N)^2+C\left(\int_Nl^2\rho\dd\mu\right)^{1/2}
			\left(\int_N|\Phi|^2\rho\dd\mu\right)^{1/2}
			\le C(D(N)^2+D(N)G(N)).
		\end{aligned}
		\]
		Here terms containing $l$ are integrated only in its tube, where
		$|l|=\dist(X,\C)$. This proves \eqref{eq:radial-tilt} and completes the proof.
	\end{proof}

	To recover Gaussian area from an annular observation, we use the following
	radial first-variation identity. Combined with Lemma~\ref{lem:radial-tilt},
	it yields the area estimate in Proposition~\ref{prop:anchorarea}.
	
	\begin{lem}
		Let $N^m\subset\R^{m+1}$ be a smooth proper boundaryless hypersurface
		with a finite uniform area-ratio bound and $G(N)<\infty$.
		Then there are smooth radial functions $f$ and $w\ge0$, depending only on
		$m$, with $\operatorname{supp}w\Subset\{1<r<2\}$, such that we have
		\begin{equation}\label{eq:radial-coefficient-bounds}
			|f(r)|\le C_m(1+r^2)^{-1},\qquad
			|f'(r)/r|\le C_m(1+r^2)^{-2},
		\end{equation}
		and
		\begin{equation}\label{eq:radial-area-identity}
			F(N)=\int_Nw\rho\dd\mu
			+\int_N\frac{f'(r)}r|X^\perp|^2\rho\dd\mu
			-\int_N f(r)\langle\Phi,X^\perp\rangle\rho\dd\mu.
		\end{equation}
		In particular, for every $m$-dimensional minimal cone, we have
		\begin{equation}\label{eq:radial-cone-normalization}
			\int_{\mathcal C}w\rho\dd\mu_{\mathcal C}=F(\mathcal C).
		\end{equation}
	\end{lem}
	\begin{proof}
		
		We choose a smooth nondecreasing function $\zeta_0$ which is zero near
		$[0,1]$ and one near $[2,\infty)$, with derivative compactly supported
		in $(1,2)$. We define
		\[
		A_m:=\int_0^\infty s^{m-1}e^{-s^2/4}\dd s,
		\]
		\[
		f(r):=r^{-m}e^{r^2/4}
		\left(\int_0^r s^{m-1}e^{-s^2/4}\dd s-A_m\zeta_0(r)\right),
		\qquad w(r):=A_m\zeta_0'(r)r^{1-m}e^{r^2/4}.
		\]
		At zero, $f(r)=m^{-1}+r^2/(2m(m+2))+O(r^4)$, with a convergent
		smooth even expansion, so the apparent singularity is removable.
		By differentiating $f$, we obtain
		\begin{equation}\label{eq:radial-ode}
			mf+rf'-\tfrac12r^2f=1-w.
		\end{equation}
		For $r\ge2$, we have
		\[
		f(r)=-r^{-m}e^{r^2/4}\int_r^\infty s^{m-1}e^{-s^2/4}\dd s
		=-2r^{-2}+O(r^{-4}).
		\]
		By integration by parts and the Gaussian tail bound, we obtain the
		last expansion with its derivative. Thus we obtain $f'(r)=O(r^{-3})$.
		Together with smoothness on compact intervals these estimates prove
		\eqref{eq:radial-coefficient-bounds}.
		
		We apply Gaussian first variation to $f(r)X$. We have
		\[
		\operatorname{div}_N(fX)-\tfrac12X\cdot(fX)
		=mf+rf'-\tfrac12r^2f-\frac{f'}r|X^\perp|^2.
		\]
		The integral of this expression equals
		$-\int_N f\langle\Phi,X^\perp\rangle\rho\dd\mu$.
		To justify this calculation, we use the compactly supported field
		$\zeta_R f(r)X$, where $\zeta_R=1$ on $B_R$, $\zeta_R=0$ outside
		$B_{2R}$, and $|D\zeta_R|\le C/R$. By \eqref{eq:radial-coefficient-bounds},
		the additional divergence term satisfies
		\[
		\left|\int_N fX^\top\cdot\nabla\zeta_R\rho\dd\mu\right|
		\le C\int_{N\cap\{R<r<2R\}}|f|\rho\dd\mu
		\le C\int_{N\cap\{r>R\}}\rho\dd\mu\longrightarrow0.
		\]
		The remaining terms converge by the area-ratio bound,
		\eqref{eq:radial-coefficient-bounds}, and Cauchy--Schwarz. In particular,
		$\|fX^\perp\|_{L^2_\rho}\le C$ and $G(N)<\infty$ justify the term
		containing $\Phi$. We let $R\to\infty$ and substitute \eqref{eq:radial-ode}
		into the first-variation identity. We obtain
		\[
		\int_N(1-w)\rho\dd\mu
		-\int_N\frac{f'}r|X^\perp|^2\rho\dd\mu
		=-\int_Nf\langle\Phi,X^\perp\rangle\rho\dd\mu,
		\]
		which proves \eqref{eq:radial-area-identity}.

		For the cone, by polar integration, we have
		$F(\mathcal C)=(4\pi)^{-m/2}\mathcal H^{m-1}(\Sigma)A_m$.
		The definition of $w$ gives the same expression multiplied by
		$\int_0^\infty\zeta_0'(r)\dd r=1$ for its annular integral.
		This proves \eqref{eq:radial-cone-normalization} directly, without a
		smooth first variation at the vertex.
		This finishes the proof.
	\end{proof}
	
	\begin{prop}\label{prop:anchorarea}
		Let $\varepsilon_{\mathcal C}>0$ be sufficiently small. Suppose $N$ is a smooth proper boundaryless hypersurface with the fixed
		area-ratio bound, has $G(N)<\infty$, and is a
		complete single normal graph $Y\mapsto Y+u(Y)n_{\mathcal C}(Y)$ over
		the cone annulus $\C\cap\{1/2<r<3\}$. Assume that
		$\|u\|_{C^1}\le\varepsilon_\C$ and, with $F_u(Y)=Y+u(Y)n_\C(Y)$,
		\begin{equation}\label{eq:annular-single-graph}
			N\cap\{3/4<r<5/2\}
			=F_u(\C\cap\{1/2<r<3\})\cap\{3/4<r<5/2\}.
		\end{equation}
		Thus every point of $N$ in $\{3/4<r<5/2\}$ belongs to this graph. Then there
		is $C_{\mathcal C}<\infty$ such that we have
		\begin{equation}\label{eq:anchor-area}
			|F(N)-F(\mathcal C)|\le C_{\mathcal C}\bigl(D(N)^2+G(N)^2\bigr).
		\end{equation}
		
	\end{prop}
	\begin{proof}
		We put $g:=w\rho$. Since $g$ is radial and $Y\cdot n_{\mathcal C}=0$, we know
		\[
		g(Y+un_{\mathcal C})=g(Y)+O(u^2)
		\]
		uniformly on $\C\cap\{1/2<r<3\}$. The graph area factor is
		\[
		\det(I-uA_{\mathcal C})
		\sqrt{1+|(I-uA_{\mathcal C})^{-1}\nabla_{\mathcal C}u|^2}.
		\]
		Minimality makes its linear term vanish. The small graph assumption,
		compact support of $w$, and absence of additional portions of the hypersurface consequently
		imply, for the compact annulus
		$\mathcal U_{\mathrm{ann}}:=\C\cap\{7/8\le r\le17/8\}$,
		\begin{equation}\label{eq:annular-area-expansion}
			\left|\int_Nw\rho\dd\mu-\int_{\mathcal C}w\rho\dd\mu_{\mathcal C}\right|
			\le C\int_{\mathcal U_{\mathrm{ann}}}(u^2+|\nabla_{\mathcal C}u|^2)\dd\mu_{\mathcal C}.
		\end{equation}
		Here $\mathcal U_{\mathrm{ann}}$ is chosen to contain every base point whose graph image meets
		the support, as well as the cone support itself.
		
		Let $l$ be signed distance to $\C$ on $\mathcal T_\eta$, and put
		$T:=|\nabla_Nl|^2$, as in the proof of Lemma~\ref{lem:radial-tilt}.
		We choose an ambient cutoff $\zeta$ supported in
		$\mathcal T_\eta\cap\{3/4<r<5/2\}$ and equal to one on the graph
		over $\mathcal U_{\mathrm{ann}}$. By \eqref{eq:local-tilt-test}, we have
		\[
		\int_{\mathrm{graph}(\mathcal U_{\mathrm{ann}})}T\rho\dd\mu\le C(D^2+DG),\qquad
		\int_{\mathrm{graph}(\mathcal U_{\mathrm{ann}})}l^2\rho\dd\mu\le D^2.
		\]
		All coefficients here are bounded because the cutoff stays in the fixed
		annulus. In graph coordinates, $l\circ F_u=u$, and the Gaussian weights and
		area Jacobians are uniformly comparable on these annuli. If
		$b=(I-uA_\C)^{-1}\nabla_\C u$, then
		\[
		T\circ F_u=\frac{|b|^2}{1+|b|^2},\qquad
		c|\nabla_\C u|^2\le T\circ F_u\le C|\nabla_\C u|^2.
		\]
		Therefore we have
		\[
		\int_{\mathcal U_{\mathrm{ann}}}(u^2+|\nabla_\C u|^2)\dd\mu_\C
		\le C\int_{\operatorname{graph}(\mathcal U_{\mathrm{ann}})}(l^2+T)\rho\dd\mu
		\le C(D^2+DG).
		\]
		Together with \eqref{eq:annular-area-expansion}, this estimates the annular
		integral of $w\rho$.

		We subtract \eqref{eq:radial-cone-normalization} from
		\eqref{eq:radial-area-identity}. By \eqref{eq:radial-coefficient-bounds},
		Cauchy--Schwarz, Lemma~\ref{lem:radial-tilt}, and Young's inequality,
		we obtain
		\[
		\begin{aligned}
			|F(N)-F(\C)|
			&\le C(D^2+DG)+CJ+CG\sqrt J\\
			&\le C(D^2+DG)+CJ+\tfrac C2(G^2+J)\\
			&\le C(D^2+DG+G^2)
			\le C(D^2+G^2).
		\end{aligned}
		\]
		Thus we obtain \eqref{eq:anchor-area}.
		Therefore we conclude the proof.
	\end{proof}

	\begin{prop}\label{prop:anchor}
		Fix $0<\alpha<1$. Along the rescaled flow of Theorem~\ref{thm:main},
		for every sufficiently small $\epsilon_{\mathrm{graph}}>0$, there is $\epsilon_{\mathrm{entry}}>0$
		with the following property. If
		$E(s)+D_0(N_s)^2<\epsilon_{\mathrm{entry}}$, then for every
		$\sigma\in[s+1,s+2]$ there is a normal graph height $u_\sigma$ on
		$\C\cap\{1/2<r<3\}$ satisfying
		\begin{equation}\label{eq:annular-entry-conclusion}
			\sup_{\sigma\in[s+1,s+2]}\|u_\sigma\|_{C^{2,\alpha}(\C\cap\{1/2<r<3\})}
			<\epsilon_{\mathrm{graph}},
		\end{equation}
		and \eqref{eq:annular-single-graph} with $N=N_\sigma$ and $u=u_\sigma$.
		In particular, the observations of Lemma~\ref{lem:gauge} are defined on
		this entire time interval.
	\end{prop}
	\begin{proof}
		We argue by contradiction. Suppose that
		$E(s_i)+D_0(N_{s_i})^2\to0$, but \eqref{eq:annular-entry-conclusion}
		or \eqref{eq:annular-single-graph} fails at some
		$\sigma_i\in[s_i+1,s_i+2]$. By monotonicity and Lemma~\ref{lem:transport}, we obtain
		\[
		\begin{aligned}
			\int_{s_i}^{s_i+3}G(N_\sigma)^2\dd\sigma
			&=E(s_i)-E(s_i+3)\le E(s_i)\longrightarrow0,\\
			\sup_{\sigma\in[s_i,s_i+3]}D_0(N_\sigma)
			&\le D_0(N_{s_i})+\sqrt3
			\left(\int_{s_i}^{s_i+3}G(N_\sigma)^2\dd\sigma\right)^{1/2}
			\longrightarrow0.
		\end{aligned}
		\]
		The compactness argument in Proposition~\ref{prop:normalized}, applied
		without division by the distance, identifies the limiting flow on
		$\{1/4<r<4\}\times(0,3)$ as the static cone. Its multiplicity is one:
		$F(N_{s_i+\sigma})=F(\C)+E(s_i+\sigma)\to F(\C)$, and the area-ratio
		bound gives uniform Gaussian integral estimates on $\{r>R\}$.
		
		We apply the local regularity theorem for mean curvature flow
		\cite[Lemma~2.1]{Schulze} on neighborhoods with closure contained in
		$\{1/4<r<4\}\times(0,3)$. A finite covering of
		$\{1/2\le r\le3\}\times[1,2]$ gives smooth normal graphs on that set and
		\[
		\sup_{1\le\sigma\le2}
		\|u_{i,\sigma}\|_{C^{2,\alpha}(\C\cap\{1/2<r<3\})}\longrightarrow0.
		\]
		The argument excluding additional points in
		Proposition~\ref{prop:normalized} shows that every point of
		$N_{s_i+\sigma}\cap\{3/4<r<5/2\}$ belongs to this graph. Thus
		\eqref{eq:annular-single-graph} holds as well. Both conclusions contradict
		the choice of $\sigma_i$. This proves the entry assertion and thus completes the proof.
	\end{proof}
	
	\section{Compactness and weighted stationary rigidity}\label{sec:compactness}
	In this section, we assume (S1), (S2), and (S3). The reference cones in this section
	are the reduced cones $\C_a$, which need not be minimal.
	We first normalize the graph and its link residual, then use the
	Hardy estimate and (S3) to rule out concentration at the vertex.
	The Gaussian inverse with residual columns gives the contraction in
	Proposition~\ref{prop:residual-contraction}.
	
	On each unit interval the reference parameter $a_i\to0$ is fixed.
	All cone functions are pulled back to $\C$ by the homogeneous link
	parametrizations. The constants in the tubular, Hardy, and parabolic
	estimates are uniform. Denote $\Theta:=F(\C)$ and, for Gaussian variance $\upsilon>0$, define
	\begin{equation}\label{eq:gaussian-variance}
		\rho_\upsilon:=(4\pi\upsilon)^{-m/2}e^{-|X|^2/(4\upsilon)},\qquad
		F_\upsilon(N):=\int_N\rho_\upsilon\dd\mu.
	\end{equation}
	Thus we have $F_1=F$.
	The fixed range of Gaussian areas is
	\begin{equation}\label{cpt:eq:window}
		0<f_*\le F(N_s^i)\le F_*<2\Theta.
	\end{equation}

	We first record the compactness statement that isolates these two possible
	obstructions to contraction.
	\begin{prop}\label{prop:normalized}
		Assume (S1) and (S2). Let $(N_s^i)_{0\le s\le1}$ be smooth proper
		rescaled flows with the fixed area-ratio bound and
		\eqref{cpt:eq:window}. Fix $a_i\to0$, let
		\[
		q_i=\nabla f(a_i),\qquad d_i=D_{a_i}(N_0^i),\qquad
		b_i^2=d_i^2+|q_i|^2>0,\qquad
		\delta_i^2=\int_0^1G(N_s^i)^2\dd s,
		\]
		and suppose $b_i\to0$ and $\delta_i/b_i\to0$.
		After passing to a subsequence, we have $q_i/b_i\to\xi\in K$. At positive times the
		hypersurfaces are complete single graphs of heights $u_i$ over
		$\C_{a_i}$ on each fixed compact regular annulus, and
		\[
		u_i/b_i\longrightarrow v\quad\hbox{smoothly},\qquad
		L_\C v=r^{-1}\xi,\qquad v\in L^2_\rho(\C).
		\]
		The limit is independent of time. For some $\alpha_{\rm vtx}\ge0$, we have
		\begin{equation}\label{cpt:eq:defect}
			b_i^{-2}\dist(X,\C_{a_i})^2\rho\dd\mu_{N_1^i}
			\longrightarrow v^2\rho\dd\mu_\C+\alpha_{\rm vtx}\delta_0,
		\end{equation}
		including convergence of total measures. If the terminal compact
		fitting observations vanish, then $\mathcal O(v)=0$, with
		$\mathcal O$ as in Proposition~\ref{prop:gaussian-bordered}.
	\end{prop}
	\begin{proof}
		By Lemma~\ref{lem:transport}, we have
		\[
		\sup_{0\le s\le1}D_{a_i}(N_s^i)\le b_i+\delta_i\le2b_i
		\]
		for large $i$. By averaging, we choose $s_i\in[1/2,3/4]$ such that
		$G(N_{s_i}^i)\le2\delta_i\to0$. Since
		$\mathbf H_i=\Phi_i-X^\perp/2$, on each fixed ball we have
		\[
		\int_{N_{s_i}^i\cap B_R}|\mathbf H_i|^2\dd\mu
		\le C_R G(N_{s_i}^i)^2+\tfrac12R^2\mu_{N_{s_i}^i}(B_R)
		\le C_R.
		\]
		By integral-varifold compactness \cite{Allard}, we obtain an integral
		limit $V$. The distance estimate and $a_i\to0$ give
		$\operatorname{spt}\|V\|\subset\C$. On $\C\setminus\{0\}$, its tangent
		plane agrees almost everywhere with $T\C$. Since $G(N_{s_i}^i)\to0$,
		$V$ is Gaussian stationary. The position vector is tangent to $\C$, so
		$V$ is also stationary for ordinary area. The constancy theorem \cite[Section~41]{SimonGMT} on the connected regular cone gives
		$V=k|\C|$ away from the vertex, for an integer $k\ge0$. The area bound
		excludes an atom at zero. By Gaussian integral estimates and
		\eqref{cpt:eq:window}, we have
		\[
		\|V\|(B_\epsilon)\le C\epsilon^m\longrightarrow0,\qquad
		0<f_*\le F(V)=k\Theta\le F_*<2\Theta.
		\]
		Thus we obtain $k=1$ and $V=|\C|$.
		
		For every compactly supported smooth ambient function $\psi$, the first variation formula implies
		\[
		\frac{\mathrm d}{\mathrm ds}\int_{N_s^i}\psi\rho\dd\mu
		=\int_{N_s^i}(D\psi\cdot\Phi_i-\psi|\Phi_i|^2)\rho\dd\mu.
		\]
		By Cauchy--Schwarz, we therefore have
		\[
		\sup_{s,\sigma\in[0,1]}\left|\int_{N_s^i}\psi\rho\dd\mu-\int_{N_\sigma^i}\psi\rho\dd\mu\right|
		\le C_\psi(\delta_i+\delta_i^2)\longrightarrow0.
		\] Thus it follows that the weight measures converge to the
		unit cone at every time, uniformly when tested against a fixed $\psi$.
		We undo the rescaling on compact spacetime cylinders away from zero.
		By Brakke compactness \cite{Brakke}, the ordinary flows converge to the
		static multiplicity-one cone. We apply local regularity
		\cite[Lemma~2.1]{Schulze} on its regular part and obtain smooth convergence.
		
		For every $0<R_1<R_2<\infty$ and $b>0$, all points of
		$N_s^i\cap\{R_1\le|X|\le R_2\}$, $b\le s\le1$, lie in the normal
		coordinate neighborhoods where the preceding smooth graphical convergence holds,
		for sufficiently large $i$. To prove this assertion, we suppose in the corresponding ordinary
		flow coordinates that terminal points $x_i\to x$ stay a positive distance
		from the limiting cone. If $\mathcal P_\tau^i$ denotes that ordinary flow
		and $\tau_0$ its terminal time, Huisken monotonicity gives, for fixed $h>0$,
		\[
		\begin{aligned}
			1&\le(4\pi h)^{-m/2}\int_{\mathcal P_{\tau_0-h}^i}
			e^{-|Y-x_i|^2/(4h)}\dd\mu\\
			&\longrightarrow(4\pi h)^{-m/2}\int_\C e^{-|Y-x|^2/(4h)}\dd\mu_\C
			\le C e^{-\dist(x,\C)^2/(8h)}\longrightarrow0\quad(h\downarrow0).
		\end{aligned}
		\]
		The area-ratio bound controls Gaussian tails in the first limit and the
		last Gaussian integral. This contradiction excludes additional points
		outside the regular graph neighborhoods.
		On the regular part away from the vertex, the graphical function $u_i$ satisfies
		\[
		\partial_su_i=\mathfrak a_i^{kl}\nabla_{kl}u_i+\mathfrak b_i^k\nabla_ku_i+c_i u_i
		-r^{-1}q_i/\mathfrak w_{a_i}.
		\]
		Its coefficients are uniformly smooth and uniformly parabolic and converge
		to those of $L_{\mathcal C}$. The residual divided by $b_i$ is
		uniformly smooth on each such annulus and converges to $-r^{-1}\xi$.
		Normal coordinates, Jacobians and Gaussian
		weights are uniformly comparable on this cylinder; hence the forward
		bound gives a uniform spacetime $L^2$ bound for $u_i/b_i$.
		Interior parabolic estimates, using a backward time interval of fixed positive length, including at
		the final time, give bounds for every derivative on smaller cylinders.
		A diagonal subsequence converges smoothly to $v(Y,s)$ satisfying
		$\partial_sv=L_{\mathcal C}v-r^{-1}\xi$. For the normal velocity, we have
		\[
		\partial_su_i\langle n_{\C_{a_i}},n_{N_s^i}\rangle
		=\langle\Phi_i,n_{N_s^i}\rangle.
		\]
		On a compact cylinder $\mathcal Q$, the normal angle is bounded away
		from zero and the coordinate measures are uniformly comparable. Hence we get
		\[
		\int_{\mathcal Q}|\partial_s(u_i/b_i)|^2\dd\mu_\C\dd s
		\le Cb_i^{-2}\int_0^1\int_{N_s^i}|\Phi_i|^2\rho\dd\mu\dd s
		=C\delta_i^2/b_i^2\longrightarrow0.
		\] Thus it follows that $v$ is
		time independent, so $L_{\mathcal C}v=r^{-1}\xi$. By exhausting the punctured cone and using the forward
		bound, we obtain $\int_\C v^2\rho\dd\mu_\C\le4$.
		
		We put
		\[
		\mu_i^{\mathrm{pos}}:=b_i^{-2}\dist(X,\C_{a_i})^2\rho\dd\mu_{N_1^i}.
		\]
		By Lemmas~\ref{lem:transport} and \ref{lem:residual-tail}, we have
		\[
		\mu_i^{\mathrm{pos}}(\R^{m+1})\le(1+\delta_i/b_i)^2\le4,\qquad
		\mu_i^{\mathrm{pos}}(\{r>R\})
		\le CR^{-2}(1+\delta_i^2/b_i^2+|q_i|^2/b_i^2).
		\]
		Thus it follows that the measures are tight. On each compact
		set away from zero their limit is exactly $v^2\rho\dd\mu_{\mathcal C}$,
		by the complete graphical convergence and the area formula. The remaining
		positive measure is supported at zero, proving \eqref{cpt:eq:defect} and
		convergence of total measures. Finally, we pass to the limit in the annular orthogonality conditions
		of Lemma~\ref{lem:gauge} for the normalized graph functions $u_i(1)/b_i$
		on $\operatorname{supp}\chi$. The family normal velocities converge to $r\psi_l$,
		so $\mathcal O(v)=0$. This concludes the proof.
	\end{proof}
	
	Proposition \ref{prop:normalized} does not determine the behavior at the vertex.
	We next obtain a Hardy estimate on a small graphical annulus and compare
	the original scale with a shrinking scale. We use this comparison to obtain smooth graphical convergence with
	multiplicity one on regular annuli at the selected terminal time.
	We call
	\begin{equation}
		\mathcal H_\C(u;\Omega):=\int_{\C\cap\Omega}
		(|\nabla u|^2+r^{-2}u^2)\dd\mu_\C
	\end{equation}
	the \textbf{graph Hardy integral}. For an ambient hypersurface or varifold,
	its \textbf{weighted Hardy distance} is
	\begin{equation}\label{eq:weighted-hardy-distance}
		\mathcal E_H(V):=\int\frac{\dist(X,\C)^2}{1+|X|^2}\dd\|V\|(X).
	\end{equation}
	\begin{lem}\label{lem:hardy-neck}
		Fix sufficiently small $r_*>0$ and $\eta_*>0$, depending only on the cone.
		Suppose $0<\varrho <r_0/4$, $r_0\le r_*$, and
		$N$ is the normal graph of $u$ over $\varrho <r<2r_0$, with
		\[
		|u|/r+|\nabla u|+r|\nabla^2u|\le\eta_*.
		\]
		We define
		\[
		G_{\rm ann}^2:=\int_{\operatorname{graph}_{\C\cap\{\varrho<r<2r_0\}}u}
		|\Phi|^2\rho\dd\mu_N.
		\]
		Then there is $C<\infty$, depending only on the cone, such that 
		\begin{equation}\label{cpt:eq:neck}
			\int_{\C\cap\{2\varrho<r<r_0\}}\left(|\nabla u|^2+r^{-2}u^2\right)\dd\mu_\C
			\le C\varrho^m+Cr_0^{-2}\int_{\C\cap\{r_0<r<2r_0\}}u^2\dd\mu_\C+Cr_0^2G_{\rm ann}^2.
		\end{equation}
		
	\end{lem}
	\begin{proof}
		Let $H_m:=(m-2)^2/4$ and $\lambda:=H_m+\mu_0>0$. For compactly
		supported $z$ on the punctured cone, \eqref{eq:radial-hardy} and
		\eqref{eq:angular-rayleigh} give
		\[
		Q^0(z)=\int_\C(|\nabla z|^2-|A_{\mathcal C}|^2z^2)\dd\mu_\C
		\ge\lambda\int_\C r^{-2}z^2\dd\mu_\C.
		\]
		Also we have $Q^0(z)\ge\int_\C|\nabla z|^2\dd\mu_\C-A_0\int_\C r^{-2}z^2\dd\mu_\C$, where
		$A_0:=\max_\Sigma|A_\Sigma|^2$. Set $\tau=\lambda/(\lambda+A_0+1)\in(0,1)$. Then we get
		\[
		\begin{aligned}
			Q^0(z)&=(1-\tau)Q^0(z)+\tau Q^0(z)\\
			&\ge\tau\int_\C|\nabla z|^2\dd\mu_\C
			+[(1-\tau)\lambda-\tau A_0]\int_\C r^{-2}z^2\dd\mu_\C\\
			&=\tau\int_\C(|\nabla z|^2+r^{-2}z^2)\dd\mu_\C.
		\end{aligned}
		\]
		Thus we have, with $c_H=\tau$,
		\begin{equation}\label{cpt:eq:coercive}
			Q^0(z)\ge c_H\int_\C(|\nabla z|^2+r^{-2}z^2)\dd\mu_\C.
		\end{equation}

		We write the area of a graph on the annulus as
		\[
		\mathcal A_{\mathrm{ann}}(u):=
		\int_{\C\cap\{\varrho<r<2r_0\}}A(Y,u,\nabla u)\dd\mu_\C,
		\]
		where the area density is given by
		\[
		A(Y,z,p):=\det(I-zA_\C(Y))
		\sqrt{1+|(I-zA_\C(Y))^{-1}p|^2}.
		\]
		It is smooth for $|z|/r+|p|\le\eta_*$. Its Hessian differs from the
		Hessian at $(0,0)$ by at most $C\eta_*$ in the scaled variables $(z/r,p)$.
		Minimality gives $D\mathcal A_{\mathrm{ann}}(0)=0$ and $D^2\mathcal A_{\mathrm{ann}}(0)=Q^0$.
		We apply the fundamental theorem of calculus to $D\mathcal A_{\mathrm{ann}}(\tau u)$.
		For every compactly supported test function $\psi$, we obtain
		\begin{align}
			|D\mathcal A_{\mathrm{ann}}(u)[\psi]-Q^0(u,\psi)|
			\le C\eta_*\int_{\C\cap\{\varrho<r<2r_0\}}\bigl(|\nabla u||\nabla\psi|
			+r^{-1}|u||\nabla\psi| +r^{-1}|\nabla u||\psi|+r^{-2}|u\psi|\bigr)\dd\mu_\C.
			\label{cpt:eq:secant}
		\end{align}

		In fact, for $0\le\tau\le1$, we know
		\[
		\begin{aligned}
			D\mathcal A_{\mathrm{ann}}(u)[\psi]-Q^0(u,\psi)
			&=\int_0^1\left(D^2\mathcal A_{\mathrm{ann}}(\tau u)-D^2\mathcal A_{\mathrm{ann}}(0)\right)[u,\psi]\dd\tau,\\
			|D\mathcal A_{\mathrm{ann}}(u)[\psi]-Q^0(u,\psi)|
			&\le C\eta_*\int (|\nabla u|+r^{-1}|u|)
			(|\nabla\psi|+r^{-1}|\psi|)\dd\mu_\C,
		\end{aligned}
		\]
		which expands into \eqref{cpt:eq:secant}.

		The exact area first variation is
		$D\mathcal A_{\mathrm{ann}}(u)[\psi]=-\int_{\C\cap\{\varrho<r<2r_0\}}(H_N\circ F_u)\psi\Theta_u\dd\mu_\C$,
		where $F_u(Y):=Y+un_{\mathcal C}$ and
		$J_u=\det(I-uA_\C)\sqrt{1+|(I-uA_\C)^{-1}\nabla u|^2}$
		is the graph area Jacobian and
		$\Theta_u=\langle n_{\mathcal C},n_N\rangle J_u$ lies between
		$1/2$ and $2$ after decreasing $\eta_*$. The graph-normal formula and
		$A_{\mathcal C}(\partial_r)=0$ (the radial direction is a zero-curvature
		direction of a cone) give
		\[
		n_N\circ F_u=\frac{n_\C-(I-uA_\C)^{-1}\nabla u}
		{\sqrt{1+|(I-uA_\C)^{-1}\nabla u|^2}},\qquad
		\langle F_u,n_N\rangle
		=\frac{u-r\partial_ru}{\sqrt{1+|(I-uA_\C)^{-1}\nabla u|^2}}.
		\]
		Thus we obtain
		\[
		|\langle F_u,n_N\rangle|\le C(|u|+r|\nabla u|).
		\]
		Writing $\phi=\langle\Phi_N,n_N\rangle$, we thus obtain
		\begin{equation}\label{cpt:eq:forcing}
			|D\mathcal A_{\mathrm{ann}}(u)[\psi]|
			\le C\int_{\C\cap\{\varrho<r<2r_0\}}(|\phi|+|u|+r|\nabla u|)|\psi|\dd\mu_\C.
		\end{equation}
		
		We choose a radial cutoff $\zeta$ supported in $(\varrho,2r_0)$ and equal to one
		on $(2\varrho,r_0)$, with derivative bounds $C/\varrho $ and $C/r_0$ on its inner
		and outer annuli. Set $\psi:=\zeta^2u$. Polarization gives the exact identity
		\[
		Q^0(u,\zeta^2u)=Q^0(\zeta u)-\int_\C u^2|\nabla\zeta|^2\dd\mu_\C.
		\]
		By expanding $\nabla(\zeta^2u)$ and applying Young's inequality, we bound
		the error in \eqref{cpt:eq:secant} by
		\[
		C\eta_*\int_\C(|\nabla(\zeta u)|^2+r^{-2}\zeta^2u^2)\dd\mu_\C
		+C\eta_*\int_\C u^2|\nabla\zeta|^2\dd\mu_\C.
		\]
		For an arbitrarily small $\sigma>0$, the first term in
		\eqref{cpt:eq:forcing} is bounded by
		\[
		C\int_{\C\cap\{\varrho<r<2r_0\}}\zeta^2|\phi u|\dd\mu_\C
		\le\sigma\int_\C r^{-2}\zeta^2u^2\dd\mu_\C+C_\sigma r_0^2\int_{\C\cap\{\varrho<r<2r_0\}}\phi^2\dd\mu_\C.
		\]
		The two terms $C\int\zeta^2u^2$ and
		$C\int\zeta^2r|\nabla u||u|$ in \eqref{cpt:eq:forcing} are bounded by
		\[
		Cr_0^2\int_\C(|\nabla(\zeta u)|^2+r^{-2}\zeta^2u^2)\dd\mu_\C
		+Cr_0^2\int_\C u^2|\nabla\zeta|^2\dd\mu_\C.
		\]
		We choose $\eta_*$ and $\sigma$ small compared with $c_H$, and then
		choose $r_*$ small. Equation~\eqref{cpt:eq:coercive} absorbs all terms involving the
		Hardy integral of $\zeta u$. Hence we have
		\[
		\int_\C(|\nabla(\zeta u)|^2+r^{-2}\zeta^2u^2)\dd\mu_\C
		\le C\int_\C u^2|\nabla\zeta|^2\dd\mu_\C+Cr_0^2\int_{\C\cap\{\varrho<r<2r_0\}}\phi^2\dd\mu_\C.
		\]
		The integral over the inner annulus is bounded by $C\varrho^m$, since $|u|\le\eta_*r$ there.
		The outer annulus gives the middle term of \eqref{cpt:eq:neck}.
		Finally the graph Jacobian is comparable to one and the Gaussian has a
		fixed positive lower bound on the graph, so $\int_{\C\cap\{\varrho<r<2r_0\}}\phi^2\dd\mu_\C\le CG_{\rm ann}^2$.
		In particular, \eqref{cpt:eq:coercive}--\eqref{cpt:eq:forcing} give
		\[
		\begin{aligned}
			\mathcal H_\C(u;\{2\varrho<r<r_0\})
			&\le\mathcal H_\C(\zeta u;\C)
			\le C\int_\C u^2|\nabla\zeta|^2\dd\mu_\C+Cr_0^2G_{\rm ann}^2\\
			&\le C\varrho^{-2}\int_{\{\varrho<r<2\varrho\}}u^2\dd\mu_\C
			+Cr_0^{-2}\int_{\{r_0<r<2r_0\}}u^2\dd\mu_\C+Cr_0^2G_{\rm ann}^2\\
			&\le C\eta_*^2\varrho^m
			+Cr_0^{-2}\int_{\{r_0<r<2r_0\}}u^2\dd\mu_\C+Cr_0^2G_{\rm ann}^2.
		\end{aligned}
		\]
		This is \eqref{cpt:eq:neck}, which completes the proof.
	\end{proof}
	
	\begin{lem}
		\label{lem:residual-neck}
		Assume (S1), and let $\C_a$ be the reduced conical family of
		Proposition~\ref{prop:link-reduction}. Define
		$q(a):=\nabla f(a)$. There are $r_*,\eta_*>0$ and $C<\infty$,
		uniform for small $a$, with the following property. If
		$0<\varrho<r_0/4$, $r_0\le r_*$, and $N$ is a normal graph of height
		$u$ over $\C_a\cap\{\varrho<r<2r_0\}$ satisfying
		\[
		|u|/r+|\nabla u|+r|\nabla^2u|\le\eta_*,
		\]
		then we have
		\[
		\begin{aligned}
			\int_{\C_a\cap\{2\varrho<r<r_0\}}
			(|\nabla u|^2+r^{-2}u^2)\dd\mu \le C\varrho^m
			+Cr_0^{-2}\int_{\C_a\cap\{r_0<r<2r_0\}}u^2\dd\mu
			+Cr_0^2G_{\rm ann}^2+Cr_0^m|q(a)|^2.
		\end{aligned}
		\]
		Here $G_{\rm ann}$ is the Gaussian norm of the actual shrinker quantity
		on the graph. No minimality of $\Sigma_a$ is assumed.
	\end{lem}
	\begin{proof}
		For compactly supported $v$ on the reference cone, set
		\[
		Q_a^0(v):=\int_{\C_a}
		(|\nabla v|^2-|A_{\C_a}|^2v^2)\dd\mu.
		\]
		Strict Hardy stability persists under small smooth link perturbations.
		The radial Hardy inequality and angular Rayleigh inequality therefore
		give, uniformly in $a$,
		\[
		Q_a^0(v)\ge c_H\int_{\C_a}
		(|\nabla v|^2+r^{-2}v^2)\dd\mu.
		\]
		Let $H_a$ denote the scalar mean curvature of the reference cone.
		The gradient identity for the reduced link gives
		$|H_a|\le C|q(a)|/r$.
		Let $\mathcal A_a(u)$ be the ordinary graph area on the annulus.
		Its first two derivatives at zero are
		\[
		D\mathcal A_a(0)[\psi]=-\int H_a\psi\dd\mu,
		\qquad
		D^2\mathcal A_a(0)[u,\psi]
		=Q_a^0(u,\psi)+\int H_a^2u\psi\dd\mu.
		\]
		These formulas follow by expanding
		$\det(I-uA_{\C_a})
		\sqrt{1+|(I-uA_{\C_a})^{-1}\nabla u|^2}$;
		the quadratic coefficient in $u$ is $\frac12(H_a^2-|A_{\C_a}|^2)$.
		The same secant estimate as in \eqref{cpt:eq:secant} holds after these
		two derivatives are subtracted, with a uniform error bounded by
		\[
		C\eta_*\int (|\nabla u|+r^{-1}|u|)
		(|\nabla\psi|+r^{-1}|\psi|)\dd\mu.
		\]
		
		Choose the radial cutoff $\zeta$ from Lemma~\ref{lem:hardy-neck} and
		put $\psi=\zeta^2u$. The term $\int H_a^2\zeta^2u^2$ is nonnegative.
		The only new first-order term is bounded, for any $\varepsilon>0$, by
		\[
		\left|\int H_a\zeta^2u\dd\mu\right|
		\le\varepsilon\int r^{-2}\zeta^2u^2\dd\mu
		+C_\varepsilon|q(a)|^2
		\mu_{\C_a}(\{r<2r_0\})
		\le\varepsilon\int r^{-2}\zeta^2u^2\dd\mu
		+C_\varepsilon r_0^m|q(a)|^2.
		\]
		The exact first variation at $u$ is controlled by the actual mean
		curvature of $N$. Because the radial direction of every cone has zero
		second fundamental form, the graph-normal formula still gives
		$|X\cdot n_N|\le C(|u|+r|\nabla u|)$.
		Thus we obtain the same forcing estimate \eqref{cpt:eq:forcing}.
		Coercivity, the secant bound, and Young's inequality now give
		\[
		\begin{aligned}
			\int (|\nabla(\zeta u)|^2+r^{-2}\zeta^2u^2)\dd\mu
			&\le C\int u^2|\nabla\zeta|^2\dd\mu
			+Cr_0^2G_{\rm ann}^2+Cr_0^m|q(a)|^2,
		\end{aligned}
		\]
		provided we first choose $\eta_*,\varepsilon$ small and then choose $r_*$ small.
		The inner cutoff term is $O(\varrho^m)$ since $|u|\le\eta_*r$;
		the outer term is the displayed annular positional term. The Gaussian
		and graph Jacobian are uniformly comparable to the unweighted cone
		measure on this fixed small ball. This proves the estimate and thus finishes the proof.
	\end{proof}
	
	We now apply the annular Hardy estimate \eqref{cpt:eq:neck} at moving spatial scales. The next
	two lemmas supply the Gaussian comparison and terminal regularity
	needed for the blowup argument.
	\begin{lem}\label{lem:scale-comparison}
		Suppose smooth hypersurfaces $N_i$ have uniform area ratios, $\lambda_i\to0$ and
		\[
		F_{b^2}(N_i)\longrightarrow\Theta\quad\hbox{for some }0<b<1,
		\qquad \lambda_i^{2-m}G(N_i)^2\longrightarrow0.
		\]
		Suppose $\lambda_i^{-1}N_i$ converges as integral varifolds to a nonzero
		stationary $V$ with
		$\int_{\R^{m+1}}\operatorname{dist}(Y,\mathcal C)^2(1+|Y|^2)^{-1}\dd \|V\|<\infty$.
		Then every blowdown of $V$ equals the multiplicity-one cone $\mathcal C$.
		The same conclusion holds if the finite weighted-distance bound is
		obtained as a limit of bounds relative to $\C_{a_i}$, with $a_i\to0$.
	\end{lem}
	\begin{proof}
		For a fixed $N_i$, we apply first variation to the position field which gives
		\[
		\frac{\mathrm d}{\mathrm d\upsilon}F_\upsilon(N_i)=\frac{1-\upsilon}{4\upsilon^2}\int_{N_i}|X^\perp|^2\rho_\upsilon\dd\mu
		+\frac1{2\upsilon}\int_{N_i}\Phi_i\cdot X^\perp\rho_\upsilon\dd\mu.
		\]
		For $0<\upsilon<1$, we complete the square to obtain
		\[
		\frac{1-\upsilon}{4\upsilon^2}|X^\perp|^2+\frac1{2\upsilon}\Phi_i\cdot X^\perp
		=\frac{1-\upsilon}{4\upsilon^2}\left|X^\perp+\frac \upsilon{1-\upsilon}\Phi_i\right|^2
		-\frac{|\Phi_i|^2}{4(1-\upsilon)}
		\ge-\frac{|\Phi_i|^2}{4(1-\upsilon)}.
		\]
		Since $\rho_\upsilon/\rho=\upsilon^{-m/2}e^{-(1-\upsilon)|X|^2/(4\upsilon)}\le \upsilon^{-m/2}$,
		we obtain
		\[
		\frac{\mathrm d}{\mathrm d\upsilon}F_\upsilon(N_i)\ge-\frac{\upsilon^{-m/2}}{4(1-\upsilon)}G(N_i)^2.
		\]
		Integrating in $\upsilon$ from $r^2$ to $b^2$, we get
		\[
		F_{r^2}(N_i)\le F_{b^2}(N_i)
		+\frac{G(N_i)^2}{2(m-2)(1-b^2)}(r^{2-m}-b^{2-m}).
		\]
		We take $r=L\lambda_i$ and let $i\to\infty$. By compact monotone
		truncations of the Gaussian, for every fixed $L>0$, we obtain
		\begin{equation}\label{eq:blowup-gaussian-upper-bound}
			F_{L^2}(V)\le\Theta\qquad(L>0).
		\end{equation}
		
		We write $\mathcal E_H=\mathcal E_H(V)<\infty$ as defined in
		\eqref{eq:weighted-hardy-distance}. For the blowdown $V_R=R^{-1}V$ and every fixed $A$, by homogeneity,
		we have
		\[
		\begin{aligned}
			\int_{B_A}\dist(Z,\C)^2\dd\|V_R\|(Z)
			&=R^{-m-2}\int_{B_{AR}}\dist(X,\C)^2\dd\|V\|(X)\\
			&\le R^{-m-2}(1+A^2R^2)
			\int_{B_{AR}}\frac{\dist(X,\C)^2}{1+|X|^2}\dd\|V\|(X)\\
			&\le(R^{-m-2}+A^2R^{-m})\mathcal E_H(V)\longrightarrow0.
		\end{aligned}
		\]
		Stationarity, the area bound and integral-varifold compactness give
		stationary integral subsequential limits supported on $\mathcal C$.
		The constancy theorem \cite[Section~41]{SimonGMT} makes each such limit $k|\mathcal C|$. There is no vertex
		atom by the area bound. Stationary monotonicity gives a density at infinity
		independent of the subsequence; it is positive since $V\ne0$. Thus the
		same integer $k\ge1$ occurs for every subsequence. Uniform Gaussian tails
		imply $F_{L^2}(V)\to k\Theta$ as $L\to\infty$. The preceding upper
		bound \eqref{eq:blowup-gaussian-upper-bound} gives
		\[
		\Theta\le k\Theta=\lim_{L\to\infty}F_{L^2}(V)\le\Theta,
		\]
		so $k=1$.
		
		For the moving-cone version, define $\widetilde V_i=(D_{1/\lambda_i})_\#|N_i|$.
		If the weighted-distance hypothesis is obtained from uniform bounds
		relative to $\C_{a_i}$, then for every compactly supported cutoff
		$0\le\zeta\le1$,
		\[
		\int\zeta(Y)\frac{\dist(Y,\C)^2}{1+|Y|^2}\dd\|V\|(Y)
		=\lim_i\int\zeta(Y)\frac{\dist(Y,\C_{a_i})^2}{1+|Y|^2}
		\dd\|\widetilde V_i\|(Y)\le C.
		\]
		Here the distance functions converge uniformly on compact sets because
		$a_i\to0$, and the local masses are uniformly bounded. Exhaustion by such
		cutoffs gives $\mathcal E_H(V)<\infty$. Thus the preceding blowdown and
		multiplicity argument applies without change. Thus we conclude the proof.
	\end{proof}
	
	\begin{lem}\label{lem:moving-compactness}
		Let the flows and frozen reduced reference cones be as in
		Proposition~\ref{prop:normalized}.
		Let $s_i\in[b_0,1]$, $b_0>0$, and $r_i\to0$, with
		\[
		r_i^{-m}\delta_i^2\longrightarrow0.
		\]
		Assume that $r_i^{-1}N_{s_i}^i\to|\mathcal C|$ as varifolds.
		Then the convergence is smooth, complete, and multiplicity one on every
		compact punctured ambient annulus. In particular, all points there belong to the corresponding
		single cone graph for sufficiently large $i$.
	\end{lem}
	\begin{proof}
		We write the rescaling as
		$N_s^i=e^{(\tau_i+s)/2}(M^i_{T_i-e^{-(\tau_i+s)}}-p_i)$.
		Put $t_i:=-e^{-(\tau_i+s_i)}$ and $\ell_i:=\sqrt{-t_i}\,r_i$.
		The exact ordinary flows
		\[
		P_\zeta^i:=\ell_i^{-1}(M^i_{T_i+t_i+\ell_i^2\zeta}-p_i)
		\]
		are defined on $[-1,0]$ for large $i$. Their relation to rescaled time is
		\[
		s=s_i-\log(1-r_i^2\zeta),\qquad
		P_\zeta^i=c_i(\zeta)N_s^i,\quad
		c_i(\zeta):=\frac{\sqrt{1-r_i^2\zeta}}{r_i}.
		\]
		The factors $\sqrt{1-r_i^2\zeta}$ converge uniformly to one.
		At corresponding points and times, we have
		\[
		|\mathbf H_P|^2=c_i^{-2}|\mathbf H_N|^2,\quad
		\mathrm d\mu_P=c_i^m \mathrm d\mu_N,\quad \mathrm d\zeta=c_i^2\mathrm ds.
		\]
		On a fixed $B_L$ in ordinary coordinates, the corresponding points have
		$|X|\le C_Lr_i$. Using the identity $\mathbf H_N=\Phi_i-X^\perp/2$, we obtain
		\begin{equation}\label{cpt:eq:movingenergy}
			\int_{-1}^0\int_{P_\zeta^i\cap B_L}|\mathbf H_P|^2\dd\mu_{P_\zeta^i}\dd\zeta
			\le C_L(r_i^{-m}\delta_i^2+r_i^4)\longrightarrow0.
		\end{equation}
		For the position term, let
		$\mathcal I_i=[s_i-\log(1+r_i^2),s_i]$, so $|\mathcal I_i|\le r_i^2$.
		By the area-ratio bound, we have
		\[
		r_i^{-m}\int_{\mathcal I_i}\int_{N_s^i\cap B_{C_Lr_i}}|X|^2\dd\mu\dd s
		\le C_Lr_i^{-m}\,r_i^2\,r_i^{m+2}=C_Lr_i^4.
		\]
		This gives the last term in \eqref{cpt:eq:movingenergy}.
		
		Brakke compactness applies by the uniform local area bounds. At almost
		every time, a bounded on the integral of squared mean curvature and integral-varifold compactness
		identify integral limiting slices; \eqref{cpt:eq:movingenergy} makes them
		stationary. More directly, for every compact smooth $\psi$, combining the identity
		\[
		\frac{\mathrm d}{\mathrm d\zeta}\int_{P_\zeta^i}\psi\dd\mu_{P_\zeta^i}
		=\int_{P_\zeta^i}(D\psi\cdot\mathbf H_P-\psi|\mathbf H_P|^2)\dd\mu_{P_\zeta^i}
		\]
		with \eqref{cpt:eq:movingenergy}, we get
		\[
		\sup_{-1\le\zeta_1<\zeta_2\le0}
		\left|\int_{P_{\zeta_2}^i}\psi\dd\mu-\int_{P_{\zeta_1}^i}\psi\dd\mu\right|
		\le C_\psi\left(\varepsilon_i(L)^{1/2}+\varepsilon_i(L)\right)\longrightarrow0,
		\]
		where $\operatorname{supp}\psi\subset B_L$ and $\varepsilon_i(L)$ is
		the left side of \eqref{cpt:eq:movingenergy}. The assumed terminal limit thus identifies
		all limiting weight measures as $|\mathcal C|$. In particular the Brakke
		limit is the static unit cone. On any spatial neighborhood compactly
		contained in its regular part, the localized unit-density regularity
		theorem \cite[Lemma~2.1]{Schulze}, also used in
		Proposition~\ref{prop:normalized}, gives
		smooth convergence through time zero on a smaller neighborhood.
		
		To exclude additional terminal points, we apply the backward Gaussian
		argument in Proposition~\ref{prop:normalized} to $P^i$ on a fixed
		backward interval. The cone Gaussian integral based a fixed positive
		distance from the cone tends uniformly to zero as the time gap tends to
		zero, while a smooth terminal point has density one. Uniform area ratios
		control the tails. Thus every terminal point in the smaller annulus belongs to the single
		graph. This proves the claimed complete convergence.
	\end{proof}
	
	We now define the inner radius used in the dilation argument. Fix $0<r_0\le\min\{r_*/8,1/8\}$ and
	$0<\eta\le\eta_*/4$, where $r_*,\eta_*$ are from
	Lemma~\ref{lem:residual-neck}. A radius $b\in(0,r_0]$ is \emph{good}
	for $N$ relative to $\C_a$ if there is a normal graph
	$F_u(Y)=Y+u(Y)n_{\C_a}(Y)$ on $\C_a\cap\{b/2<r<4b\}$ such that
	\[
	\begin{gathered}
		N\cap\{b<r<2b\}
		=F_u(\C_a\cap\{b/2<r<4b\})\cap\{b<r<2b\},\\
		\sup_{\C_a\cap\{b/2<r<4b\}}
		(|u|/r+|\nabla u|+r|\nabla^2u|)<\eta.
	\end{gathered}
	\]
	A radius not satisfying these
	conditions is \emph{bad}. Define the \textbf{inner radius} by
	\begin{equation}\label{eq:inner-radius-definition}
		\varrho(N;\C_a):=\sup(\{b\in(0,r_0]:b\text{ is bad}\}\cup\{0\}).
	\end{equation}
	Thus every $b>\varrho$ up to $r_0$ is good; if $\varrho>0$, some bad
	radius lies in $(\varrho/2,\varrho]$.
	\begin{lem}\label{lem:first-bad}
		Assume (S1), (S2), and (S3). Let the sequence be as in
		Proposition~\ref{prop:normalized}, and choose
		$s_i\in[1/2,3/4]$ with $G(N_{s_i}^i)\le2\delta_i$.
		Use the good-radius definition above relative to $\C_{a_i}$,
		and put $\varrho_i=\varrho(N_{s_i}^i;\C_{a_i})$.
		Then $\varrho_i^m=o(b_i^2)$, and, for every fixed sufficiently small
		$\epsilon_H>0$,
		\begin{equation}\label{cpt:eq:wholehardy}
			\int_{N_{s_i}^i\cap B_{\epsilon_H}}
			\frac{\dist(X,\C_{a_i})^2}{|X|^2}\dd\mu
			\le C(\varrho_i^m+b_i^2).
		\end{equation}
	\end{lem}
	
	\begin{proof}
		By Proposition~\ref{prop:normalized}, every fixed compact range of radii
		is eventually good. Thus we obtain $\varrho_i\to0$. The graphs on overlapping
		annuli agree by uniqueness of normal projection. We choose the graphical
		threshold below that of Lemma~\ref{lem:residual-neck}. On the fixed outer
		annulus, normalized smooth convergence gives
		\[
		r_0^{-2}\int_{\C_{a_i}\cap\{r_0<r<2r_0\}}u_i^2\dd\mu_{\C_{a_i}}
		\le Cb_i^2,\qquad G(N_{s_i}^i)^2\le4\delta_i^2=o(b_i^2).
		\]
		By Lemma~\ref{lem:residual-neck}, we have
		\[
		\begin{aligned}
			\mathcal H_{\C_{a_i}}(u_i;\{2\varrho_i<r<r_0\})
			&\le C\varrho_i^m+Cr_0^{-2}\int_{\{r_0<r<2r_0\}}u_i^2\dd\mu_{\C_{a_i}}
			+Cr_0^2G(N_{s_i}^i)^2+Cr_0^m|q_i|^2\\
			&\le C\varrho_i^m+Cb_i^2+C\delta_i^2
			\le C(\varrho_i^m+b_i^2).
		\end{aligned}
		\]
		Choose $\epsilon_H>0$ sufficiently small relative to $r_0$ that the normal
		projection of the graphical part in $B_{\epsilon_H}$ lies in $\{r<r_0\}$.
		Since $|F_{u_i}(Y)|\asymp|Y|$, the graph Jacobians are bounded and
		$\dist(F_{u_i}(Y),\C_{a_i})\le|u_i(Y)|$, we have
		\begin{equation}
			\begin{aligned}
				\int_{N_{s_i}^i\cap B_{\epsilon_H}}
				\frac{\dist(X,\C_{a_i})^2}{|X|^2}\dd\mu
				&\le\mu_{N_{s_i}^i}(B_{4\varrho_i})
				+C\int_{\C_{a_i}\cap\{2\varrho_i<r<r_0\}}r^{-2}u_i^2\dd\mu_{\C_{a_i}}\\
				&\le C\varrho_i^m+C\mathcal H_{\C_{a_i}}(u_i;\{2\varrho_i<r<r_0\})\\
				&\le C(\varrho_i^m+b_i^2).
			\end{aligned}
		\end{equation}
		Here the graph smallness ensures that $|X|>4\varrho_i$ implies that its
		normal projection has radius greater than $2\varrho_i$.
		This proves \eqref{cpt:eq:wholehardy}. If $\varrho_i=0$, we apply Lemma~\ref{lem:residual-neck}
		with inner cutoff radius $\sigma$ and let $\sigma\downarrow0$.
		The inner-annulus bound $C\sigma^m$ disappears by monotone convergence. Thus it follows that the
		estimate \eqref{cpt:eq:wholehardy} holds also in this case, since
		\[
		\int_{N_{s_i}^i\cap B_{\epsilon_H}}\frac{\dist(X,\C_{a_i})^2}{|X|^2}\dd\mu_{N_{s_i}^i}
		\le\liminf_{\sigma\downarrow0}C(\sigma^m+b_i^2)=Cb_i^2.
		\]
		
		Next, we argue by contradiction. Suppose that $\varrho_i\ge c b_i^{2/m}$ along a
		subsequence. Choose a bad radius $r_i\in(\varrho_i/2,\varrho_i]$.
		Since $4r_i>\varrho_i$, the corresponding annulus is a complete small
		graph. After dilation by the factor $r_i^{-1}$, the graph has area bounded below
		by a positive constant independent of $i$. On each fixed ball, we have
		\[
		\int_{r_i^{-1}N_{s_i}^i\cap B_L}|\mathbf H|^2\dd\mu_{r_i^{-1}N_{s_i}^i}
		\le C_L(r_i^{2-m}G(N_{s_i}^i)^2+r_i^4)\longrightarrow0.
		\]
		Indeed, $r_i^m\ge c'b_i^2$ and $\delta_i/b_i\to0$ give
		\[
		r_i^{-m}\delta_i^2\le C\delta_i^2/b_i^2\longrightarrow0,\qquad
		r_i^{2-m}G(N_{s_i}^i)^2\le4r_i^2\,r_i^{-m}\delta_i^2\longrightarrow0.
		\]
		By integral-varifold compactness, we obtain a stationary integral limit
		$V$. The annular area lower bound gives $V\ne0$.
		On every fixed ball, the rescaled form of \eqref{cpt:eq:wholehardy} gives
		\[
		\int_{r_i^{-1}N_{s_i}^i}
		\frac{\dist(Y,\C_{a_i})^2}{1+|Y|^2}\dd\mu_{r_i^{-1}N_{s_i}^i}
		\le C r_i^{-m}(\varrho_i^m+b_i^2)\le C,
		\]
		where the integral is restricted to the region corresponding to
		$B_{\epsilon_H}$ before scaling. Every fixed ball eventually lies there.
		Lower semicontinuity and exhaustion show $\mathcal E_H(V)<\infty$.
		
		We check the multiplicity hypothesis needed for rigidity. For a fixed
		$0<b<1$, the test $\rho_{b^2}/\rho$ and its gradient are bounded.
		Gaussian first-variation estimates and $\delta_i\to0$ imply that $F_{b^2}(N_s^i)$
		varies by $o(1)$ over the whole time interval. Its value at one fixed interior
		time tends to $\Theta$ by Proposition~\ref{prop:normalized} and uniform
		Gaussian tails. Thus we obtain $F_{b^2}(N_{s_i}^i)\to\Theta$.
		We can now apply Lemma~\ref{lem:scale-comparison} with
		$\lambda_i=r_i$, so the blowdown of $V$ is the unit cone.
		By (S3), we obtain $V=|\mathcal C|$.
		
		Finally $r_i^{-m}\delta_i^2\le C\delta_i^2/b_i^2\to0$.
		Lemma~\ref{lem:moving-compactness} upgrades the selected-slice convergence
		to complete smooth convergence through its terminal time on the larger
		annulus defining badness. The radius $r_i$ must therefore be good for
		large $i$, a contradiction. This proves the asserted estimate and therefore completes the proof.
	\end{proof}

	\begin{prop}
		\label{prop:residual-contraction}
		Assume (S1), (S2), and (S3).
		Use the reduced links and compact fitting of
		Proposition~\ref{prop:link-reduction} and
		Lemma~\ref{lem:gauge}. Fix the area-ratio bound and
		$0<f_*<F(\C)<F_*<2F(\C)$.
		For every $\eta>0$ there are $\varepsilon_\eta,C_\eta>0$ with the
		following property. Let $(N_s)_{0\le s\le1}$ be a unit rescaled-flow
		interval with the fixed area-ratio bound and
		$f_*\le F(N_s)\le F_*$. Assume that the terminal slice lies in the
		fitting neighborhood of Lemma~\ref{lem:gauge}, and put
		\[
		a:=a(N_1),\qquad q(a):=\nabla f(a),\qquad
		\delta^2:=\int_0^1G(N_s)^2\dd s.
		\]
		If $|a|+D_a(N_0)+\delta<\varepsilon_\eta$, then we have
		\[
		D_a(N_1)^2+|q(a)|^2
		\le\eta D_a(N_0)^2+C_\eta\int_0^1G(N_s)^2\dd s.
		\]
	\end{prop}
	\begin{proof}
		Suppose that the assertion fails. There are $a_i\to0$ and unit intervals
		with $d_i=D_{a_i}(N_0^i)\to0$ such that
		\[
		e_i^2+|q_i|^2>\eta d_i^2+i\delta_i^2,\qquad
		e_i=D_{a_i}(N_1^i),\quad q_i=\nabla f(a_i).
		\]
		Set $b_i^2=d_i^2+|q_i|^2>0$. Lemma~\ref{lem:transport}, applied to the distance from the frozen cone $\C_{a_i}$, gives
		$e_i\le d_i+\delta_i$, so we obtain
		\[
		(i-2)\delta_i^2<2d_i^2+|q_i|^2\le2b_i^2,
		\qquad \delta_i/b_i\longrightarrow0.
		\]
		Proposition~\ref{prop:normalized} gives a time-independent punctured
		limit $v$ with $L_\C v=\mathcal B\xi=r^{-1}\xi$ and
		$\mathcal O(v)=0$, where $q_i/b_i\to\xi$.
		
		By Lemma~\ref{lem:first-bad}, the selected-time bad radius satisfies
		$\varrho_i^m=o(b_i^2)$. Fix a sufficiently small $r_0$. For
		$0<\epsilon<r_0$, Lemma~\ref{lem:residual-neck} and normalized smooth
		convergence at the selected times give
		\begin{equation}\label{eq:residual-limit-hardy}
			\begin{aligned}
				&\int_{\C\cap\{\epsilon<r<r_0\}}
				(|\nabla v|^2+r^{-2}v^2)\dd\mu_\C\\
				&\quad\le\liminf_i b_i^{-2}
				\mathcal H_{\C_{a_i}}
				(u_i;\{2\varrho_i<r<r_0\})\\
				&\quad\le C\limsup_i\left(
				\frac{\varrho_i^m}{b_i^2}+1
				+\frac{\delta_i^2}{b_i^2}+\frac{|q_i|^2}{b_i^2}\right)
				\le C.
			\end{aligned}
		\end{equation}
		The constant is independent of $\epsilon$. Monotone convergence as
		$\epsilon\downarrow0$ gives finite gradient and Hardy integrals at
		the vertex.
		
		We next check infinity. Choose a radial cutoff $\chi_R$ supported in
		$\{r_0/2<r<2R\}$, equal to one on $\{r_0<r<R\}$, with
		$|D\chi_R|\le C_{r_0}$ on its inner transition and $C/R$ on its outer
		transition. The equation $L_\C v=\mathcal B\xi$ and Gaussian
		integration by parts give
		\[
		\begin{aligned}
			\int\chi_R^2|\nabla v|^2\rho\dd\mu
			&=-\int\chi_R^2v\mathcal B\xi\rho\dd\mu
			+\int(|A_\C|^2+\tfrac12)\chi_R^2v^2\rho\dd\mu
			-2\int\chi_Rv\nabla\chi_R\cdot\nabla v\rho\dd\mu\\
			&\le\tfrac12\int\chi_R^2|\nabla v|^2\rho\dd\mu
			+C_{r_0}\int_{\{r>r_0/2\}}
			(v^2+|\mathcal B\xi|^2)\rho\dd\mu.
		\end{aligned}
		\]
		The last integral is finite because $v,\mathcal B\xi\in L^2_\rho$.
		Absorption and Fatou's lemma give a finite Gaussian gradient integral
		at infinity. By \eqref{eq:gaussian-moment-control} applied to compactly supported
		cutoffs of $v$ and the bound in \eqref{eq:residual-limit-hardy}, we obtain
		\[
		\int_\C r^2v^2\rho\dd\mu
		\le C\int_\C(v^2+|\nabla v|^2)\rho\dd\mu<\infty.
		\]
		For an inner cutoff $\eta_\epsilon$, equal to zero on $r\le\epsilon$
		and one on $r\ge2\epsilon$, we have
		\[
		\begin{aligned}
			\|(1-\eta_\epsilon)v\|_\V^2
			&\le C\int_{\{r<2\epsilon\}}
			\left( |\nabla v|^2+(1+r^{-2}+r^2)v^2\right)\rho\dd\mu
			\longrightarrow0.
		\end{aligned}
		\]
		By the integrability just established, the outer cutoff approximations
		converge to $v$ in the same norm. Therefore we have $v\in\V$. Both sides of the equation are
		continuous in $\V^*$, so density extends its punctured weak form to
		all tests in $\V$. Proposition~\ref{prop:gaussian-bordered} now gives
		\begin{equation}\label{eq:residual-limit-zero}
			\|v\|_\V+|\xi|
			\le C\left(\|L_\C v-\mathcal B\xi\|_{\V^*}+|\mathcal O(v)|\right)=0.
		\end{equation}
		
		It remains to exclude the terminal vertex atom. Choose a smooth radial
		cutoff $\zeta$ equal to one on $[0,1]$ and zero on $[2,\infty)$, and put
		\[
		J_{i,\epsilon}(s)=\int_{N_s^i}\zeta(r/\epsilon)^2
		\dist(X,\C_{a_i})^2\rho\dd\mu.
		\]
		For $2\epsilon<\epsilon_H$, \eqref{cpt:eq:wholehardy} gives
		\[
		\begin{aligned}
			b_i^{-2}J_{i,\epsilon}(s_i)
			\le4\epsilon^2 b_i^{-2}
			\int_{N_{s_i}^i\cap B_{2\epsilon}}
			r^{-2}\dist(X,\C_{a_i})^2\rho\dd\mu \le C\epsilon^2\left(1+\frac{\varrho_i^m}{b_i^2}\right).
		\end{aligned}
		\]
		Thus we obtain $\limsup_i b_i^{-2}J_{i,\epsilon}(s_i)\le C\epsilon^2$.
		Since the distance to a cone is at most $r$, the derivative of the
		cutoff squared-distance function has size at most a constant times
		the distance. Gaussian first variation and Lipschitz approximation
		therefore give
		\begin{equation}\label{eq:lw-localized-transport}
			|J_{i,\epsilon}(1)-J_{i,\epsilon}(s_i)|
			\le C(d_i+\delta_i)\delta_i+4\epsilon^2\delta_i^2=o(b_i^2).
		\end{equation}
		No minimality of the fixed reference cone is used in this transport
		calculation. Passing to the limit as $\epsilon\downarrow0$, we conclude that the limiting
		measure at $s=1$ has no atom at the vertex.
		By the convergence of total measures in \eqref{cpt:eq:defect},
		$e_i^2/b_i^2\to0$. Since $q_i/b_i\to0$, also
		$d_i^2/b_i^2\to1$. After dividing the inequality by $b_i^2$ and passing to the limit,
		we obtain $0\ge\eta$, which contradicts $\eta>0$.
		This completes the proof.
	\end{proof}
	
	\section{Uniqueness, Lojasiewicz inequality and rigidity}\label{sec:main-proofs}
	In this section, we first prove a local two-unit inequality with the reduced-link
	residual. We then establish continuation on finite intervals before
	summing the inequality to infinity. This gives uniqueness and the
	quantitative statements. The final subsection treats constant
	rescaled flows.
	
	\subsection{Spatial comparisons and the local discrete inequality}\label{sec:discrete}
	Throughout this subsection, set $\rho:=(4\pi)^{-m/2}e^{-|X|^2/4}$ and $\Phi:=\mathbf H+\frac{X^\perp}{2}$.
	We hold the reference cone $\C_a$ fixed in each estimate and omit $a$
	when no confusion can arise. The constants are uniform for small $a$.
	The next two lemmas for the truncated Gaussian distance require only
	a closed reference cone with vertex at the origin. The subsequent discrete inequality and rigidity arguments
	use the full cone hypotheses of Theorem~\ref{thm:main}.
	We use $l$ for signed distance in a normal neighborhood, as in
	Section~\ref{sec:preliminaries}. For the family $\C_a$, we define the
	nonnegative pointwise distance and the truncated Gaussian distance by
	\begin{equation}\label{eq:pointwise-truncated-distance}
		j_a(X):=\dist(X,\C_a),\qquad
		B_a(N)^2:=\int_N\frac{j_a(X)^2}{1+j_a(X)^2}\rho\dd\mu.
	\end{equation}
	When the reference cone is fixed, we suppress the parameter and write
	$j(X)$, $B(N)$, and $D(N)$. Thus we have
	\[
	j(X)=\dist(X,\C),\qquad
	B(N)^2=\int_N\frac{j(X)^2}{1+j(X)^2}\rho\dd\mu,\qquad
	D(N)^2=\int_N j(X)^2\rho\dd\mu.
	\]
	In integrands we also write $j=j(X)$. We recall that
	$G(N)^2=\int_N|\Phi|^2\rho\dd\mu$.
	The truncated Gaussian distance $B$ permits a backward time estimate in addition to
	the forward estimate of Lemma~\ref{lem:transport}. To use it for Gaussian
	area, we must first recover the unbounded distance $D$ from it. The next
	lemma does this with an error controlled by the Gaussian norm $G$ of the
	shrinker quantity.
	
	\begin{lem}\label{l:lem:comparison}
		Let $\C_a$ be a closed cone with vertex at the origin, and let $N$ be a
		smooth properly embedded boundaryless hypersurface with the fixed
		area-ratio bound. Then we have
		\[
		B_a(N)\le D_a(N)\le C_m\bigl(B_a(N)+G(N)\bigr).
		\]
		The constant depends only on $m$. In particular, the estimate applies
		to every cone in the reduced family, without a minimality assumption.
	\end{lem}
	\begin{proof}
		The first inequality follows directly from
		$j_a^2/(1+j_a^2)\le j_a^2$. For the second inequality, we may assume
		$G(N)<\infty$. Define
		\[
		j=j_a(X)=\dist(X,\C_a),\qquad
		u:=\frac{j^2}{1+j^2},\qquad B=B_a(N),\quad D=D_a(N),\quad G=G(N).
		\]
		Since $0\in\C_a$, we have $j\le r$. The distance function is
		one-Lipschitz, so almost everywhere on $N$ we obtain
		\begin{equation}\label{eq:truncated-moment-pointwise}
			0\le u\le1,\qquad u^2\le u,\qquad
			|\nabla u|^2\le\frac{4j^2}{(1+j^2)^4}\le4u,\qquad
			j^2=u+(1+j^2)u^2\le2u+r^2u^2.
		\end{equation}
		The Gaussian moment estimate \eqref{lem:exterior-tail:10} holds for
		any compactly supported Lipschitz function on $N$: its proof uses
		only the divergence identity for $X^\top u^2\rho$ and Young's inequality.
		Apply it first to smooth approximations of $\zeta_Ru$, where
		$\zeta_R=1$ on $B_R$, $\zeta_R=0$ outside $B_{2R}$, and
		$|D\zeta_R|\le C/R$. Smooth approximation in the local Sobolev norm
		justifies the Lipschitz test. Then let $R\to\infty$.
		The area-ratio bound controls the Gaussian moments and the cutoff
		errors; boundedness of $u$ and $G<\infty$ control the velocity term.
		Thus \eqref{lem:exterior-tail:10} and
		\eqref{eq:truncated-moment-pointwise} give
		\begin{equation}\label{eq:truncated-distance-moment}
			\begin{aligned}
				\int_Nr^2u^2\rho\dd\mu
				&\le C_m\int_N(u^2+|\nabla u|^2+|\Phi|^2u^2)\rho\dd\mu\\
				&\le C_m\left(\int_Nu\rho\dd\mu+
				\int_N|\Phi|^2\rho\dd\mu\right)
				=C_m(B^2+G^2).
			\end{aligned}
		\end{equation}
		Using the last inequality in \eqref{eq:truncated-moment-pointwise},
		we consequently obtain
		\[
		\begin{aligned}
			D^2=\int_Nj^2\rho\dd\mu
			&\le2\int_Nu\rho\dd\mu+\int_Nr^2u^2\rho\dd\mu\\
			&\le2B^2+C_m(B^2+G^2)\le C_m(B^2+G^2).
		\end{aligned}
		\]
		After taking square roots, we obtain the desired estimate. This finishes the proof.
	\end{proof}
	
	Lemma~\ref{l:lem:comparison} provides a spatial comparison between the two
	distance quantities. We next estimate the change of the truncated Gaussian distance
	along the flow. Its bounded integrand controls the quadratic term in first
	variation and yields estimates in both time directions.
	
	\begin{lem}\label{l:lem:transport}
		Let $(N_s)_{s_1\le s\le s_2}$ satisfy the hypotheses of
		Lemma~\ref{lem:transport}, with $s_2-s_1\le1$. Fix the cone parameter
		$a$ on this interval and put $\delta^2:=\int_{s_1}^{s_2}G(N_\sigma)^2\dd\sigma$. Then we have
		\[
		B_a(N_{s_2})\le B_a(N_{s_1})+\delta,\qquad
		B_a(N_{s_1})\le B_a(N_{s_2})+3\delta.
		\]
	\end{lem}
	\begin{proof}
		We use the pointwise distance in \eqref{eq:pointwise-truncated-distance} to define the function
		$f=j_a/(1+j_a^2)^{1/2}$. This function satisfies $0\le f\le1$ and
		has Lipschitz constant at most one. Convolution with nonnegative smooth
		kernels gives smooth functions $f_\epsilon$ with the same bounds, converging
		uniformly to $f$. For $I_\epsilon(\sigma):=\int_{N_\sigma}f_\epsilon^2\rho\dd\mu$,
		the Gaussian first-variation identity with outer cutoffs gives
		\[
		I_\epsilon'=\int_{N_\sigma}
		(2f_\epsilon Df_\epsilon\cdot\Phi-f_\epsilon^2|\Phi|^2)\rho\dd\mu,
		\qquad
		-2\sqrt{I_\epsilon}\,G-G^2\le I_\epsilon'
		\le2\sqrt{I_\epsilon}\,G.
		\]
		The Gaussian area bound and $\int_{s_1}^{s_2}G(N_\sigma)^2\dd\sigma<\infty$ justify removal of the
		outer cutoff. Uniform convergence of $f_\epsilon$ and the uniform Gaussian
		area bound permit passage to the limit in the integrated inequalities.
		Thus, for $I(\sigma):=B_a(N_{\sigma})^2$, the same differential inequalities hold in the
		absolutely continuous sense:
		\[
		-2\sqrt I\,G-G^2\le I'\le2\sqrt I\,G.
		\]
		The upper differential inequality gives a bound for the derivative of
		$(I+\epsilon)^{1/2}$. Integrating this bound and letting
		$\epsilon\downarrow0$, we obtain the forward assertion and
		$\sup_{[s_1,s_2]}\sqrt I\le B_a(N_{s_1})+\delta$.
		Integration of the lower inequality yields
		\[
		B_a(N_{s_1})^2-B_a(N_{s_2})^2\le2B_a(N_{s_1})\delta+3\delta^2.
		\]
		Therefore we have $(B_a(N_{s_1})-\delta)^2\le B_a(N_{s_2})^2+4\delta^2$, which implies
		$B_a(N_{s_1})\le B_a(N_{s_2})+3\delta$.
		This concludes the proof.
	\end{proof}

	\begin{lem}
		\label{lem:residual-spatial}
		Let $\C_a$ be the reduced family, set
		$q(a)=\nabla f(a)$, and use the fixed all-center area bound.
		For a smooth proper boundaryless hypersurface $N$,  write $D_a$ and
		$B_a$ for its unbounded and truncated distances to
		$\C_a$. Uniformly for small $a$, we have
		\[
		D_a(N)\le C_m\bigl(B_a(N)+G(N)\bigr).
		\]
		If $N$ has the complete small graph on the fixed annulus
		required by Proposition~\ref{prop:anchorarea}, then
		\[
		|F(N)-F(\C_a)|
		\le C\left(D_a(N)^2+G(N)^2+|q(a)|^2\right).
		\]
	\end{lem}
	\begin{proof}
		The first assertion is Lemma~\ref{l:lem:comparison}, which applies
		also to the nonminimal cone $\C_a$.
		
		All relative tubes and their coefficients are uniform for the smooth
		compact family of links. Let $l_a$ be signed distance in such a tube,
		and let $r=|X|$. The cone identities used in the spatial estimates are
		\[
		X\cdot Dl_a=l_a,\qquad D^2l_a(Dl_a,\cdot)=0,\qquad
		|D^2l_a|\le C/r.
		\]
		The trace bound which uses minimality is replaced by
		\[
		|\Delta l_a|\le C\left(|l_a|/r^2+|q(a)|/r\right).
		\]
		Indeed the principal curvatures of a parallel hypersurface are the
		reference curvatures divided by $1-l_a\kappa_j$, up to the common
		sign convention. Their reference sum is bounded by $C|q(a)|/r$.
		
		For the area estimate, put $D=D_a(N)$, $G=G(N)$, and
		$Q=|q(a)|$. We first establish the counterpart of
		\eqref{eq:local-tilt-test}. Set $n_a=Dl_a$ and
		$T_a:=|\nabla_Nl_a|^2=1-(n_N\cdot n_a)^2$.
		For $\zeta$ compactly supported in the relative tube away from the origin,
		we apply the Gaussian first variation formula to the test field
		$\zeta^2 l_a n_a$ to obtain
		\begin{align*}
			\int_N\bigl(&\zeta^2T_a+\zeta^2l_a\Delta l_a
			-\zeta^2l_aD^2l_a(n_N,n_N)
			+2\zeta l_a\nabla_N\zeta\cdot\nabla_Nl_a
			-\tfrac12\zeta^2l_a^2\bigr)\rho\dd\mu
			&=-\int_N\zeta^2l_a\Phi\cdot n_a\rho\dd\mu.
		\end{align*}
		Since $D^2l_a$ annihilates $n_a$, its contribution is bounded by
		$C\eta\zeta^2T_a$ in a tube of relative width $\eta$.
		Choose $\eta$ small, and use Young's inequality on the cutoff term.
		The bound for $\Delta l_a$ then gives
		\begin{equation}\label{eq:residual-local-tilt-test}
			\int_N\zeta^2T_a\rho\dd\mu
			\le C\int_N\left(
			\zeta^2|l_a||\Phi|
			+[\zeta^2(1+r^{-2})+|D\zeta|^2]l_a^2
			+Q\zeta^2|l_a|/r\right)\rho\dd\mu.
		\end{equation}
		Thus \eqref{eq:residual-local-tilt-test} differs from
		\eqref{eq:local-tilt-test} only in its last term.
		The area bound implies
		\begin{equation}\label{eq:residual-inverse-radius-integral}
			\int_N r^{-2}\rho\dd\mu
			\le C\Lambda\sum_{k=0}^\infty2^{-k(m-2)}+C\Lambda
			\le C\Lambda.
		\end{equation}
		Here the sum estimates the annuli
		$2^{-k-1}<r<2^{-k}$, and the exterior integral is bounded by
		Gaussian decay. We use $m\ge3$ in the convergence of the sum.
		Consequently, for bounded $\zeta$, Cauchy--Schwarz gives
		\begin{equation}\label{eq:residual-tilt-source}
			Q\int_N\zeta^2|l_a|/r\,\rho\dd\mu
			\le CQ\left(\int_{N\cap\operatorname{supp}\zeta}l_a^2\rho\dd\mu\right)^{1/2}
			\left(\int_Nr^{-2}\rho\dd\mu\right)^{1/2}
			\le CQD.
		\end{equation}
		All terms involving $l_a$ are restricted to its tube, where
		$|l_a|=j_a(X)$.
		
		Choose the angular cutoff $\xi$ as in the proof of
		\eqref{eq:radial-to-normal-tilt}, now relative to $\C_a$.
		The geometric argument for \eqref{eq:radial-to-normal-tilt} uses only
		$X\cdot n_a=l_a$ and gives
		\[
		J(N)\le CD^2+C\int_N\xi^2\frac{r^2}{(1+r^2)^2}T_a\rho\dd\mu.
		\]
		Apply \eqref{eq:residual-local-tilt-test} with
		$\zeta=\xi r/(1+r^2)$. The bounds for $D\zeta$ used in the proof
		of \eqref{eq:radial-tilt} remain uniform in $a$. Inner and outer
		cutoffs are removed by the area bound and Gaussian decay, as in the proof of \eqref{eq:radial-tilt}; the additional source is integrable by
		\eqref{eq:residual-tilt-source}. Therefore we obtain
		\begin{equation}\label{eq:residual-radial-tilt}
			\begin{aligned}
				J(N)&\le CD^2+C\int_N\left(
				\zeta^2|l_a||\Phi|
				+[\zeta^2(1+r^{-2})+|D\zeta|^2]l_a^2
				+Q\zeta^2|l_a|/r\right)\rho\dd\mu\\
				&\le C(D^2+DG+QD).
			\end{aligned}
		\end{equation}
		
		We next estimate the annular observation. Write
		$F_u(Y)=Y+u(Y)n_{\C_a}(Y)$ for the given complete small graph,
		and let $\mathcal U_{\mathrm{ann},a}$ denote the fixed compact
		annulus used in \eqref{eq:annular-area-expansion}, with $\C$
		replaced by $\C_a$. Since $Y\cdot n_{\C_a}=0$, the radial function
		$g=w\rho$ satisfies $g(F_u(Y))=g(Y)+O(u^2)$.
		The graph area factor satisfies
		\[
		\det(I-uA_{\C_a})
		\sqrt{1+|(I-uA_{\C_a})^{-1}\nabla_{\C_a}u|^2}
		=1-u\operatorname{tr}A_{\C_a}
		+O(u^2+|\nabla_{\C_a}u|^2).
		\]
		On this annulus $|\operatorname{tr}A_{\C_a}|\le CQ$.
		The support condition and completeness of the graph thus give the
		counterpart of \eqref{eq:annular-area-expansion}:
		\begin{equation}\label{eq:residual-annular-area-expansion}
			\left|\int_Nw\rho\dd\mu-\int_{\C_a}w\rho\dd\mu_{\C_a}\right|
			\le C\int_{\mathcal U_{\mathrm{ann},a}}
			(u^2+|\nabla_{\C_a}u|^2+Q|u|)\dd\mu_{\C_a}.
		\end{equation}
		For a cutoff equal to one on this graph annulus,
		\eqref{eq:residual-local-tilt-test} and
		\eqref{eq:residual-tilt-source} give
		\[
		\int_{\operatorname{graph}(\mathcal U_{\mathrm{ann},a})}
		(l_a^2+T_a)\rho\dd\mu\le C(D^2+DG+QD).
		\]
		Indeed all cutoff coefficients are bounded on the fixed annulus.
		Moreover, $l_a\circ F_u=u$ and, with
		$b=(I-uA_{\C_a})^{-1}\nabla_{\C_a}u$, we have
		$T_a\circ F_u=|b|^2/(1+|b|^2)$.
		The small graph assumption makes $T_a\circ F_u$ comparable to
		$|\nabla_{\C_a}u|^2$, and makes the Gaussian weights and graph
		Jacobians uniformly comparable. Cauchy--Schwarz on the fixed
		annulus therefore yields
		\begin{equation}\label{eq:residual-annular-observation}
			\begin{aligned}
				\int_{\mathcal U_{\mathrm{ann},a}}
				(u^2+|\nabla_{\C_a}u|^2+Q|u|)\dd\mu_{\C_a}
				&\le C(D^2+DG+QD)
				+CQ\left(\int_{\mathcal U_{\mathrm{ann},a}}u^2\dd\mu_{\C_a}\right)^{1/2}\\
				&\le C(D^2+DG+QD).
			\end{aligned}
		\end{equation}
		
		Finally, polar integration in the definition of $w$ gives
		$\int_{\C_a}w\rho\dd\mu_{\C_a}=F(\C_a)$.
		This identity requires only that $\C_a$ is a cone.
		Subtract it from \eqref{eq:radial-area-identity} and apply
		\eqref{eq:radial-coefficient-bounds},
		\eqref{eq:residual-annular-area-expansion}, and
		\eqref{eq:residual-annular-observation}. Thus we obtain
		\[
		\begin{aligned}
			|F(N)-F(\C_a)|
			&\le C(D^2+DG+QD)
			+C\int_N\frac{|X^\perp|^2}{(1+r^2)^2}\rho\dd\mu
			+C\int_N\frac{|\Phi||X^\perp|}{1+r^2}\rho\dd\mu\\
			&\le C(D^2+DG+QD)+CJ+CG\sqrt J\\
			&\le C(D^2+DG+QD)+CJ+C(G^2+J)\\
			&\le C(D^2+DG+QD+G^2)
			\le C(D^2+G^2+Q^2).
		\end{aligned}
		\]
		Here we used Cauchy--Schwarz in the second line,
		Young's inequality in the third and last lines, and
		\eqref{eq:residual-radial-tilt} in the fourth line.
		This proves the second assertion and completes the proof.
	\end{proof}
	
	\begin{prop}
		\label{prop:local-loja}
		Assume (S1), (S2), and (S3). Let $\theta$ be a reduced-link gradient
		exponent in Proposition~\ref{prop:link-reduction}.
		For a rescaled flow on $[s,s+2]$, assume the fixed area bounds,
		$E\ge0$, exact dissipation, and complete sufficiently small regular
		annular graphs throughout the interval. Suppose that the reduced
		fitting parameters, fitted distances, and $E$ remain in a sufficiently
		small fixed neighborhood. Then we have
		\[
		E(s)^{2(1-\theta)}\le C\left(E(s)-E(s+2)\right).
		\]
	\end{prop}
	\begin{proof}
		Fix $a=a(N_{s+2})$ throughout the interval and put
		\[
		\delta^2=\int_s^{s+2}G(N_t)^2\dd t,
		\qquad d=D_a(N_{s+1}),\quad e=D_a(N_{s+2}),\quad q=q(a).
		\]
		Choose $t_0\in[s,s+1]$ with $G(N_{t_0})\le\delta$.
		The truncated-distance transport estimate in
		Lemma~\ref{l:lem:transport} does not require a minimal reference cone.
		By applying this estimate on at most two unit intervals, we obtain
		\[
		B_a(N_{t_0})
		\le B_a(N_{s+2})+3\sqrt2\delta.
		\]
		By Lemmas~\ref{lem:transport} and \ref{l:lem:comparison}, we have
		\begin{equation}\label{eq:residual-midpoint-comparison}
			\begin{aligned}
				d&\le D_a(N_{t_0})+\int_{t_0}^{s+1}G(N_t)\dd t\\
				&\le C\left(B_a(N_{t_0})+G(N_{t_0})\right)+\delta\\
				&\le C\left(B_a(N_{s+2})+\delta\right)
				\le C_0(e+\delta).
			\end{aligned}
		\end{equation}
		The constant $C_0$ is fixed before choosing a contraction parameter.
		Proposition~\ref{prop:residual-contraction} on the second unit interval gives
		\[
		e^2+|q|^2\le\eta d^2+C_\eta\delta^2.
		\]
		Choose $\eta>0$ so that $C_0\sqrt{\eta}<1/2$.
		Then, by \eqref{eq:residual-midpoint-comparison}, we have
		\begin{equation}\label{eq:residual-two-unit-distance}
			\begin{aligned}
				d&\le C_0(e^2+|q|^2)^{1/2}+C_0\delta\\
				&\le C_0\sqrt{\eta}\,d
				+C_0(\sqrt{C_\eta}+1)\delta,\\
				d+e+|q|&\le C\delta.
			\end{aligned}
		\end{equation}
		All thresholds may be decreased after this fixed choice of $\eta$.
		
		For $t\in[s+7/4,s+2]$, forward transport gives
		$D_a(N_t)\le d+\delta$. Monotonicity, the spatial area estimate in
		Lemma~\ref{lem:residual-spatial}, and
		\eqref{eq:residual-two-unit-distance} therefore imply
		\begin{equation}\label{eq:residual-terminal-area}
			\begin{aligned}
				F(N_{s+2})-F(\C_a)
				&\le4\int_{s+7/4}^{s+2}|F(N_t)-F(\C_a)|\dd t\\
				&\le C\int_{s+7/4}^{s+2}
				\left(D_a(N_t)^2+G(N_t)^2+|q|^2\right)\dd t\\
				&\le C(d^2+\delta^2+|q|^2)\le C\delta^2.
			\end{aligned}
		\end{equation}
		Thus, by exact dissipation, we obtain
		\begin{equation}\label{eq:residual-initial-area}
			\begin{aligned}
				E(s)=E(s+2)+\delta^2 \le F(\C_a)-F(\C)+C\delta^2
				\le |F(\C_a)-F(\C)|+C\delta^2.
			\end{aligned}
		\end{equation}
		Equation~\eqref{eq:reduced-cone-gradient-comparison} gives
		$|F(\C_a)-F(\C)|^{1-\theta}\le C|q|\le C\delta$.
		Since $2(1-\theta)\ge1$ and $\delta\le1$, \eqref{eq:residual-initial-area}
		and concavity of $x^{1-\theta}$ give
		\[
		E(s)^{1-\theta}\le C\left(|q|+\delta^{2(1-\theta)}\right)\le C\delta.
		\]
		We then square the estimate to prove the assertion.
	\end{proof}
	
	\subsection{Continuation and uniqueness}\label{sec:uniqueness}
	\begin{prop}\label{prop:continuation}
		Under the hypotheses of Theorem~\ref{thm:main}, the reduced fitting
		chart and the hypotheses of Proposition~\ref{prop:local-loja}
		remain valid for all sufficiently late times. The rescaled flow has
		finite Gaussian length and satisfies \eqref{eq:energy-loja} and
		\eqref{eq:energy-length}. Its fitted parameter also satisfies
		\begin{equation}\label{eq:fitted-parameter-rate}
			|a(s)|\le CE(s)^\beta
			\le C(1+s-s_0)^{-\beta/(1-2\beta)}\qquad(s\ge s_0).
		\end{equation}
		The rescaled flow converges locally as Radon measures to
		$|\C|$ in the original coordinates, smoothly on compact subsets of
		$\operatorname{Reg}(\C)$. Every tangent flow is the same static cone.
	\end{prop}
	\begin{proof}
		The energy identity and the prescribed tangent subsequence give
		$E\ge0$, $E\to0$, and $E'=-G^2$.
		Choose an exponent $0<\beta<1/2$ no larger than the reduced-link
		exponent $\theta$ and decrease the neighborhood so that $E\le1$.
		The estimate in Proposition~\ref{prop:local-loja} then holds with $\theta$
		replaced by $\beta$.
		
		We first justify continuation before summing over an infinite interval.
		The unnormalized annular-entry argument in Proposition~\ref{prop:anchor}
		uses only the smooth minimal base cone, the area bound, vanishing
		positional distance and dissipation, and interior local regularity.
		The reduced fitting map from
		Lemma~\ref{lem:gauge} is defined on those annular graphs.
		Choose a graph threshold so small that its reduced fit and fitted
		annular graph meet the hypotheses of
		Proposition~\ref{prop:local-loja} whenever the fixed-cone
		Gaussian distance $D_0$ and the energy are sufficiently small.
		This is possible by continuity of the fitting map and
		$D_a(N)\le D_0(N)+C|a|$.
		Fix $r_D>0$ and $\varepsilon_E>0$ below the resulting entry thresholds,
		so that $r_D^2+\varepsilon_E<\epsilon_{\mathrm{entry}}$. We also choose
		$\varepsilon_E$ small enough to keep the Gaussian areas in the interval
		required by Proposition~\ref{prop:residual-contraction}.
		On any interval where $D_0<r_D$ and $E<\varepsilon_E$, complete
		annular graphs with the chosen threshold are available after one
		unit of time, by Proposition~\ref{prop:anchor}.
		
		We choose a sufficiently large time $t$ from the prescribed sequence
		along which the rescaled flow converges to $\C$ and set $t_0=t+1$. Unnormalized convergence on a fixed larger time interval gives
		the required complete graphs on $[t_0,t_0+1]$ and small reduced fits.
		The uniform Gaussian tails and forward positional transport permit
		$D_0(N_{t_0})$ and $E(t_0)$ to be as small as needed. We require
		\[
		D_0(N_{t_0})<r_D/4,\qquad E(t_0)<\varepsilon_E,
		\qquad C E(t_0)^\beta+\sqrt{2E(t_0)}<r_D/4,
		\]
		where the constant $C$ below is fixed by the local discrete estimate.
		The entry argument and monotonicity maintain annular graphs and the
		reduced fitting chart up to any possible first exit of $D_0$ from
		$[0,r_D)$. For the first unit this follows from the chosen return;
		subsequently it follows by applying the entry argument one unit earlier.
		
		Assume that the first exit time $T_*$ is finite. Set $t_j=t_0+2j$ and take
		$J$ so that $t_J<T_*\le t_{J+1}$, using the preceding time interval when
		$T_*$ is a grid time. On every complete time interval with $j<J$, define
		\[
		e_j=E(t_j),\qquad \delta_j^2=e_j-e_{j+1}.
		\]
		By Proposition~\ref{prop:local-loja}, we have
		$e_j^{1-\beta}\le C\delta_j$. Concavity gives, at each positive energy,
		\begin{equation}\label{eq:continuation-concavity}
			\begin{aligned}
				e_j^\beta-e_{j+1}^\beta
				&=\beta\int_{e_{j+1}}^{e_j}x^{\beta-1}\dd x
				\ge\beta e_j^{\beta-1}(e_j-e_{j+1})\\
				&=\beta e_j^{\beta-1}\delta_j^2
				\ge (\beta/C)\delta_j.
			\end{aligned}
		\end{equation}
		Zero energies cause no division and have zero subsequent dissipation.
		Cauchy--Schwarz on each interval of length two and
		\eqref{eq:continuation-concavity} give
		\begin{equation}\label{eq:continuation-finite-length}
			\begin{aligned}
				\int_{t_0}^{t_J}G(N_t)\dd t
				&=\sum_{j=0}^{J-1}\int_{t_j}^{t_{j+1}}G(N_t)\dd t\\
				&\le\sqrt2\sum_{j=0}^{J-1}
				\left(\int_{t_j}^{t_{j+1}}G(N_t)^2\dd t\right)^{1/2}
				=\sqrt2\sum_{j=0}^{J-1}\delta_j\\
				&\le C\sum_{j=0}^{J-1}(e_j^\beta-e_{j+1}^\beta)
				=C(e_0^\beta-e_J^\beta)\le CE(t_0)^\beta.
			\end{aligned}
		\end{equation}
		The same bound holds with both sums empty when $J=0$.
		On the final incomplete time interval, exact dissipation alone gives
		$\int_{t_J}^{T_*}G\le\sqrt{2E(t_0)}$.
		Lemma~\ref{lem:transport}, applied to the fixed base cone and combined with \eqref{eq:continuation-finite-length}, now implies
		\[
		D_0(N_{T_*})\le D_0(N_{t_0})
		+C E(t_0)^\beta+\sqrt{2E(t_0)}<r_D/2,
		\]
		contradicting the first exit. Graphical entry and the reduced fitting
		therefore continue for all later times.
		
		Proposition~\ref{prop:local-loja} consequently applies at every
		sufficiently late real starting time, and proves
		\eqref{eq:energy-loja} after increasing $s_0$.
		Fix any $s\ge s_0$ and set $e_k=E(s+2k)$ and
		$\Delta_k=(e_k-e_{k+1})^{1/2}$. The calculation in
		\eqref{eq:continuation-concavity} applies at these starting times.
		For every integer $n\ge1$, we therefore obtain
		\begin{equation}\label{eq:continuation-tail-length}
			\begin{aligned}
				\int_s^{s+2n}G(N_\sigma)\dd\sigma
				&\le\sqrt2\sum_{k=0}^{n-1}\Delta_k
				\le C\sum_{k=0}^{n-1}(e_k^\beta-e_{k+1}^\beta)\\
				&=C(E(s)^\beta-E(s+2n)^\beta)\le CE(s)^\beta.
			\end{aligned}
		\end{equation}
		Passing to the limit as $n\to\infty$ in
		\eqref{eq:continuation-tail-length}, we obtain the length bound in
		\eqref{eq:energy-length}.
		
		To obtain the decay rate, we choose $\vartheta=1-2\beta>0$.
		By \eqref{eq:energy-loja}, we have
		$e_k-e_{k+1}\ge c e_k^{1+\vartheta}$.
		Whenever $e_{k+1}>0$, we obtain
		\begin{equation}\label{eq:continuation-reciprocal-energy}
			\begin{aligned}
				e_{k+1}^{-\vartheta}-e_k^{-\vartheta}
				&=\vartheta\int_{e_{k+1}}^{e_k}x^{-\vartheta-1}\dd x\\
				&\ge\vartheta e_k^{-\vartheta-1}(e_k-e_{k+1})
				\ge c\vartheta.
			\end{aligned}
		\end{equation}
		Applying \eqref{eq:continuation-reciprocal-energy} with $s=s_0$
		and sum over $k$, we obtain
		\[
		E(s_0+2k)\le
		\left(E(s_0)^{-\vartheta}+c\vartheta k\right)^{-1/\vartheta}
		\le C(1+k)^{-1/\vartheta}.
		\]
		If an energy vanishes, exact dissipation makes all later energies and
		lengths zero. Monotonicity between the times $s_0+2k$ gives the
		continuous-time decay bound in \eqref{eq:energy-length}.
		
		For $\psi\in C_c^1(\mathbb R^{m+1})$, Gaussian first variation and
		Cauchy--Schwarz give, for $t>s\ge t_0$,
		\[
		\begin{aligned}
			\left|\int_{N_t}\psi\rho\dd\mu-\int_{N_s}\psi\rho\dd\mu\right|
			&\le C_\psi\int_s^tG(N_\sigma)\dd\sigma
			+\|\psi\|_\infty\int_s^tG(N_\sigma)^2\dd\sigma\\
			&\le C_\psi\int_s^\infty G(N_\sigma)\dd\sigma
			+\|\psi\|_\infty E(s)\longrightarrow0.
		\end{aligned}
		\]
		The prescribed tangent subsequence identifies the limit with
		$\int_\C\psi\rho\dd\mu_\C$. Positivity and smoothness of $\rho$ on
		compact sets give local Radon convergence in the original coordinates.
		For any parabolic scales $\lambda_i\downarrow0$, let
		$s_i=-2\log\lambda_i$. At every fixed $\tau<0$, we have
		\[
		\lambda_i^{-1}(M_{T+\lambda_i^2\tau}-x_0)
		=\sqrt{-\tau}\,N_{s_i-\log(-\tau)}
		\longrightarrow\sqrt{-\tau}\,\C=\C.
		\]
		Thus it follows that every tangent flow is the same static cone. By applying interior regularity to the corresponding mean curvature
		flows, we obtain smooth convergence on compact subsets of
		$\operatorname{Reg}(\C)$.
		Finally, smooth convergence on the fitting annulus and continuity
		of the fitting map in Lemma~\ref{lem:gauge} imply $a(s)\to0$.
		By the parameter-speed estimate in Lemma~\ref{lem:gauge} and
		\eqref{eq:continuation-tail-length}, we have
		\[
		\begin{aligned}
			|a(s)|&=\left|\lim_{t\to\infty}(a(s)-a(t))\right|
			\le\int_s^\infty|a'(\sigma)|\dd\sigma\\
			&\le C\int_s^\infty G(N_\sigma)\dd\sigma
			\le CE(s)^\beta
			\le C(1+s-s_0)^{-\beta/(1-2\beta)}.
		\end{aligned}
		\]
		This proves \eqref{eq:fitted-parameter-rate} and completes the proof.
	\end{proof}
	
	\begin{proof}[Proof of Theorems~\ref{thm:main} and \ref{thm:energy-loja}]
		It follows from Proposition~\ref{prop:continuation} directly.
	\end{proof}

	\subsection{Rigidity among smooth self-shrinkers}

	\begin{proof}[Proof of Corollary~\ref{cor:shrinker-rigidity}]
		We argue by contradiction. Suppose that smooth normalized shrinkers
		$S_i$ converge locally with
		multiplicity one to a nonflat cone $\C$ as stated. We first obtain
		uniform Gaussian bounds from this local convergence. Properness implies
		completeness, and the volume-growth theorem of Ding--Xin
		\cite[Theorem~1.1]{DingXin} gives polynomial volume growth for each
		noncompact $S_i$; for a compact component the required integrability
		is immediate. At a minimum of $|X|^2$ on each component of $S_i$, we have
		$X=X^\perp$ and $\Delta_{S_i}|X|^2\ge0$, so the identity
		$\Delta_{S_i}|X|^2=2m-|X^\perp|^2$ shows that every component
		meets $\overline B_{\sqrt{2m}}$.
		By properness, only finitely many components meet this ball, so the
		componentwise integrability implies that every Gaussian polynomial
		moment of $S_i$ is finite.
		The constants at this stage may depend on $i$.
		
		The first-variation calculation in Lemma~\ref{lem:scale-comparison},
		with $\Phi=0$, gives
		\begin{equation}\label{eq:shrinker-variance-monotonicity}
			\frac{\mathrm d}{\mathrm d\upsilon}F_\upsilon(S_i)
			=\frac{1-\upsilon}{4\upsilon^2}
			\int_{S_i}|X^\perp|^2\rho_\upsilon\dd\mu
			\le0\qquad(\upsilon\ge1).
		\end{equation}
		Polynomial volume growth justifies differentiation and removal of the
		outer cutoffs in this calculation without a uniform area-ratio bound.
		In particular, we have $F_2(S_i)\le F(S_i)$.
		Since $\rho/\rho_2=2^{m/2}e^{-r^2/8}$, we obtain
		\begin{equation}\label{eq:shrinker-gaussian-tail}
			\begin{aligned}
				\int_{S_i\setminus B_R}\rho\dd\mu
				\le2^{m/2}e^{-R^2/8}F_2(S_i) \le2^{m/2}e^{-R^2/8}F(S_i).
			\end{aligned}
		\end{equation}
		Choose a fixed $R_0$ so large that
		$2^{m/2}e^{-R_0^2/8}\le1/2$. Splitting the Gaussian integral at $R_0$
		and applying \eqref{eq:shrinker-gaussian-tail}, we have
		\begin{equation}\label{eq:shrinker-uniform-gaussian-bound}
			F_2(S_i)\le F(S_i)
			\le2\int_{S_i\cap B_{R_0}}\rho\dd\mu\le C.
		\end{equation}
		The last bound follows from local varifold convergence. Moreover,
		$\dist(X,\C)\le r$, so \eqref{eq:shrinker-uniform-gaussian-bound}
		gives
		\begin{equation}\label{eq:shrinker-distance-tail}
			\begin{aligned}
				\int_{S_i\setminus B_R}\dist(X,\C)^2\rho\dd\mu
				\le2^{m/2}\sup_{r\ge R}(r^2e^{-r^2/8})F_2(S_i) \le C\sup_{r\ge R}(r^2e^{-r^2/8})\longrightarrow0
				\qquad(R\to\infty),
			\end{aligned}
		\end{equation}
		uniformly in $i$. Local convergence, together with
		\eqref{eq:shrinker-gaussian-tail}--\eqref{eq:shrinker-distance-tail},
		therefore yields
		\begin{equation}\label{eq:shrinker-global-gaussian-convergence}
			F(S_i)\longrightarrow F(\C),\qquad D_0(S_i)\longrightarrow0.
		\end{equation}
		
		We next derive the all-center area bound. For $x\in\R^{m+1}$ and
		$\upsilon>0$, put
		\[
		F_{x,\upsilon}(S_i):=(4\pi\upsilon)^{-m/2}
		\int_{S_i}e^{-|X-x|^2/(4\upsilon)}\dd\mu.
		\]
		Apply Huisken's monotonicity formula \cite{Huisken} to the ancient
		flow $\sqrt{-t}\,S_i$, with spacetime center $(x,-1+R^2)$, between
		times $-L$ and $-1$. After dilation of the initial slice, we obtain
		\[
		F_{x,R^2}(S_i)
		\le F_{x/\sqrt L,\,(L+R^2-1)/L}(S_i)\qquad(L>1).
		\]
		For fixed $i,x,R$, the Gaussian on the right converges pointwise to
		$\rho$. For all sufficiently large $L$ it is bounded by a constant
		multiple of $\rho_2$, so \eqref{eq:shrinker-uniform-gaussian-bound}
		and dominated convergence permit $L\to\infty$. Thus we obtain
		\begin{equation}\label{eq:shrinker-all-center-area-bound}
			\begin{aligned}
				\mathcal H^m(S_i\cap B_R(x))
				\le e^{1/4}(4\pi)^{m/2}R^mF_{x,R^2}(S_i) \le e^{1/4}(4\pi)^{m/2}R^mF(S_i)\le CR^m.
			\end{aligned}
		\end{equation}
		The constant in \eqref{eq:shrinker-all-center-area-bound} is independent
		of $i,x,R$. By \eqref{eq:shrinker-global-gaussian-convergence}, the
		Gaussian areas lie in $[F(\C)/2,3F(\C)/2]$ for all sufficiently large $i$.
		On each fixed regular annulus their mean curvature is bounded by $|X|/2$.
		Multiplicity-one varifold convergence and Allard regularity give complete
		multiplicity-one $C^{1,\alpha}$ convergence there; the shrinker graph equation
		and interior elliptic estimates make it smooth. Additional points away
		from the cone are excluded by the local area lower bound for bounded
		mean curvature, since such points would give positive measure outside the
		limiting support.
		
		Lemma~\ref{lem:gauge} gives reduced fitted parameters $a_i\to0$.
		Parameter comparison implies $d_i=D_{a_i}(S_i)\to0$. Regard each
		$S_i$ as a constant rescaled flow on a fixed larger time interval.
		Its shrinker quantity is zero and its fitted distance is $d_i$ at
		both endpoints. Proposition~\ref{prop:residual-contraction}, with a fixed
		$0<\eta<1$, gives
		\[
		d_i^2+|q(a_i)|^2\le\eta d_i^2.
		\]
		The Gaussian area interval, area bound, annular graph, and fitting
		hypotheses have all been verified above. The contraction does not
		require a prescribed future tangent or nonnegative relative Gaussian
		area. Thus we obtain $d_i=0$ and $q(a_i)=0$ for large $i$.
		
		Positivity of $\rho$ and continuity of distance imply
		$S_i\subset\C_{a_i}$. At a nonzero point this inclusion of hypersurfaces
		of the same dimension is open in the regular cone. Properness makes its
		image closed. Since the regular cone is connected and $S_i$ has a
		nonzero point, it contains the whole regular cone and its vertex.
		The nearby cones remain nonflat: otherwise a sequence of their totally
		geodesic links would converge smoothly to the nonflat link of $\C$.
		Therefore we obtain $S_i=\C_{a_i}$, contradicting smoothness at the vertex.
		This proves the sequential assertion; metrizability of local varifold
		convergence gives the neighborhood formulation.
		This completes the proof.
	\end{proof}

	\section{Examples}\label{sec:examples}
	In this section, we verify the spectral and stationary rigidity conditions
	for the cones in Corollary~\ref{cor:examples}.
	\subsection{Strictly minimizing cones}\label{sec:strict-minimization}
	\paragraph{Definition of strict minimization.}
	Let $\C$ be a nonflat regular area-minimizing hypersurface cone
	satisfying (S1). Write $E_+$ and $E_-$ for the two components of its
	complement, and choose $n_\Sigma$ pointing into $E_+$. The
	Hardt--Simon theorem \cite[Theorem~2.1]{HS} gives smooth embedded
	area-minimizing hypersurfaces $L_\pm\subset E_\pm$, normalized by
	$\dist(0,L_\pm)=1$, whose positive dilations foliate $E_\pm$.
	Outside a compact set they are normal graphs over $\C$ of signed
	heights $h_\pm$. With $\phi_0>0$ normalized in $L^2(\Sigma)$, we
	call $\C$ \emph{strictly minimizing} if there are constants
	$c_+,c_->0$ such that
	\begin{equation}\label{eq:strict-minimizing-leaves}
		h_\pm(r,\omega)=\pm c_\pm r^{\gamma_0^+}\phi_0(\omega)
		+o_2(r^{\gamma_0^+})\qquad(r\to\infty).
	\end{equation}
	Here $a=o_2(r^\gamma)$ means
	\[
	\lim_{R\to\infty}\sum_{j=0}^2R^{j-\gamma}
	\|\nabla^ja\|_{C^0(\C\cap\{R<r<2R\})}=0.
	\]
	The normal and signs in \eqref{eq:strict-minimizing-leaves} are fixed
	on the cone, so the requirement is imposed on both complementary
	sides. Since (S1) gives $\gamma_0^-<\gamma_0^+$, this requires the
	larger principal branch on each side. This is the usual
	Hardt--Simon strict-minimizing condition; see
	\cite[Theorem~3.2]{HS} and \cite[Section~2.2]{Wang} for its
	equivalent formulations.
	
	\begin{prop}\label{e:strict-s3}
		A strictly minimizing cone satisfying (S1) satisfies (S3).
	\end{prop}
	\begin{proof}
		The flat case is immediate. We assume that $\C$ is nonflat, and let $V$
		satisfy the hypotheses in (S3). By the definition in Section~\ref{sec:strict-minimization}, $\C$ is
		area minimizing and satisfies \eqref{eq:strict-minimizing-leaves}.
		Since the blowdown is the smooth multiplicity-one cone away from its vertex,
		Allard regularity \cite{Allard} and the minimal graph equation imply that,
		outside a compact set, $V$ is a single normal graph over $\C$. We write $h$ for this normal graph height, defined on $\C\cap\{r>R_2\}$
		for some $R_2>0$. Thus, for all sufficiently large $R_1$, every point of
		$\operatorname{spt}\|V\|\setminus B_{R_1}$ has the form
		$Y+h(Y)n_\C(Y)$ with $Y\in\C\cap\{r>R_2\}$.
		To quantify its decay,  put
		$\epsilon_R^2:=\int_{\{R/4<|X|<8R\}}\dist(X,\C)^2(1+|X|^2)^{-1}\dd\|V\|$.
		Then we have $\epsilon_R\to0$. Setting $\widetilde h_R(Y)=R^{-1}h(RY)$, we obtain
		\begin{equation}\label{eq:strict-minimizer-height-decay}
			\begin{aligned}
				\|\widetilde h_R\|_{L^2(\C\cap\{1/2<r<4\})}^2
				&\le CR^{-m}\epsilon_R^2,\\
				\sup_{\{R<r<2R\}}\sum_{k=0}^2R^k|\nabla^kh|
				&\le CR\|\widetilde h_R\|_{L^2(\C\cap\{1/2<r<4\})}
				\le CR^{1-m/2}\epsilon_R=o(R^{-(m-2)/2}).
			\end{aligned}
		\end{equation}
		For the second line, we apply the interior elliptic estimate to the
		small minimal graph equation. Thus we obtain $h=o(r^{-(m-2)/2})$. If
		$V\ne|\C|$, write $\sigma$ for its asymptotic rate in
		\cite[Section~4.1]{Wang}. \eqref{eq:strict-minimizer-height-decay} gives
		$\sigma\le-(m-2)/2$, whereas the strictly minimizing case of
		\cite[Theorem~4.5(1)]{Wang} and (S1) give
		\[
		\gamma_0^+\le\sigma\le-\frac{m-2}{2}
		<-\frac{m-2}{2}+\nu_0=\gamma_0^+.
		\]
		This contradiction proves (S3).
		This finishes the proof.
	\end{proof}
	
	\subsection{The spectrum of quadratic cones}
	In this subsection, we compute the link spectrum of $\C_{p,q}$ in \eqref{e:quadratic},
	verifying the two spectral conditions in Proposition~\ref{e:spectrum}.
	
	Throughout this subsection $p,q$ denote the dimensions of the sphere. Let
	\begin{equation}\label{e:quadratic}
		m=p+q+1,\qquad
		\mathcal C_{p,q}
		=\{(x,y)\in\mathbb R^{p+1}\oplus\mathbb R^{q+1}:
		q|x|^2=p|y|^2\},\qquad p,q\ge1.
	\end{equation}
	Its link is
	\[
	\Sigma_{p,q}:=S^p\left(\sqrt{\frac p{m-1}}\right)
	\times S^q\left(\sqrt{\frac q{m-1}}\right).
	\]
	\begin{prop}\label{e:spectrum}
		Let $\C_{p,q}$ be the cone in \eqref{e:quadratic}, with $m\ge7$.
		Then $\C_{p,q}$ satisfies \textup{(S1) and (S2)}, and its link
		spectrum is given by
		\begin{equation}\label{e:product-spectrum}
			\mu_{k,l}=(m-1)\left[
			\frac{k(k+p-1)}p+\frac{l(l+q-1)}q-1\right],
			\qquad k,l\in\mathbb Z_{\ge0}.
		\end{equation}
	\end{prop}
	\begin{proof}
		The principal curvatures of the link are $\sqrt{q/p}$, with multiplicity
		$p$, and $-\sqrt{p/q}$, with multiplicity $q$. In particular their sum is
		zero and the squared norm of the second fundamental form is $m-1$.
		The eigenvalues of the negative Laplacian on $S^d(R_{\mathrm{sph}})$ are
		$k(k+d-1)/R_{\mathrm{sph}}^2$. Using the tensor-product basis of spherical harmonics, we compute the
		Laplacian eigenvalues. After subtracting $m-1$ from these eigenvalues,
		we obtain \eqref{e:product-spectrum}.
		
		The bottom eigenvalue is $\mu_0=-(m-1)$. The strict Hardy inequality
		\textup{(S1)} follows from
		\[
		\frac{(m-2)^2}{4}-(m-1)=\frac{m^2-8m+8}{4}>0
		\qquad(m\ge7).
		\]
		To verify (S2), we put $T_d(k):=k(k+d-1)/d$. Then we have
		\[
		T_d(0)=0,\qquad T_d(1)=1,\qquad
		T_d(k)\ge T_d(2)>2\quad(k\ge2).
		\]
		If a forbidden value in \textup{(S2)} were an eigenvalue, we would have
		\begin{equation}
			T_p(k)+T_q(l)=
			\frac{2[(m-1)-sm+2s^2]}{m-1},\qquad 1\le s<m/4.
		\end{equation}
		For $s=1$ the right side is $2/(m-1)\in(0,1)$. For $s\ge2$, since
		$m>4s$, its numerator divided by two is strictly less than
		$-2s^2+4s-1<0$. The left side is either zero or at least one, a
		contradiction in both cases. This proves \textup{(S2)}.
		
		Thus we conclude the proof.
	\end{proof}
	
	\subsection{Weighted rigidity for quadratic cones}
	For a quadratic cone $\C_{p,q}$ of dimension $m=p+q+1\ge7$, the proof of
	(S3) uses the quantitative regularity theorem in Lemma~\ref{e:es}.
	We treat the minimizing and exceptional nonminimizing cases in
	Propositions~\ref{e:minimizing-s3} and \ref{e:exceptional-s3}, respectively.
	
	\begin{lem}\label{e:es}
		Let $\mathcal C=\mathcal C_{p,q}$ be an area-minimizing quadratic cone
		of dimension $m$, and put
		\[
		\gamma_+=\frac{-(m-2)+\sqrt{(m-2)^2-4(m-1)}}2.
		\]
		Let $S_\lambda$ be its signed Hardt--Simon dilation family, normalized so
		that $S_0=\mathcal C$ and $R S_\lambda=S_{R\lambda}$.
		There are constants $\delta,\Lambda,c,\beta>0$, depending on the cone,
		with the following property. Suppose $M$ is a stationary integral
		$m$-varifold in $B_1$, $|\lambda|\le\Lambda$, and
		\[
		\int_{B_1}\operatorname{dist}(X,S_\lambda)^2\dd \|M\|(X)
		\le E\le\delta^2,
		\]
		\[
		\|M\|(B_1)\le\tfrac32\|\mathcal C\|(B_1),\qquad
		\|M\|(B_{1/10})\ge\tfrac12\|\mathcal C\|(B_{1/10}).
		\]
		Then there are $z_*\in\mathbb R^{m+1}$, $Q\in SO(m+1)$, and
		$\lambda'\in\mathbb R$ such that we have
		\begin{equation}\label{e:es-parameters}
			|z_*|+|Q-\mathrm{Id}|+
			|(\lambda')^{1-\gamma_+}-\lambda^{1-\gamma_+}|
			\le cE^{1/2}.
		\end{equation}
		Here signed powers mean $\lambda^{1-\gamma_+}=\operatorname{sign}(\lambda)|\lambda|^{1-\gamma_+}$.
		In $B_{1/2}$ the support of $M$ is a single normal graph $u$ over the
		corresponding portion of $z_*+Q S_{\lambda'}$. The graph function satisfies
		\begin{equation}\label{e:es-graph}
			r^{-1}\|u\|_{C^0(B_r(z_*))}+\|\nabla u\|_{C^0(B_r(z_*))}
			+r^\beta[\nabla u]_{\beta,B_r(z_*)}
			\le cr^\beta E^{1/2}\quad(0<r\le1/2).
		\end{equation}
		The graph and its norms are understood on the regular reference
		hypersurface; when $\lambda'=0$ this is the punctured cone.
	\end{lem}
	The statement above is the form of \cite[Theorem~3.1]{ES} that we use.
	By \cite[Sections~2.3 and~2.5]{ES}, the leading term of the normal graph
	of a foliation leaf over the cone has radial dependence $r^{\gamma_+}$,
	where $\gamma_+$ is the larger
	root of the indicial equation displayed above. The theorem applies to
	any stationary integral varifold $M$ satisfying the bounds stated above.
	
	\begin{lem}\label{e:subgrowth}
		Let $V$ be a stationary integral $m$-varifold with multiplicity-one
		blowdown $\mathcal C_{p,q}$, where $\mathcal C_{p,q}$ is area-minimizing.
		Assume that
		\begin{equation}\label{e:subgrowth-condition}
			I(R):=\int_{B_R}\operatorname{dist}(X,\mathcal C_{p,q})^2
			\dd \|V\|(X)=o(R^{m+2\gamma_+}).
		\end{equation}
		Then we have $V=|\mathcal C_{p,q}|$.
	\end{lem}
	\begin{proof}
		Let $V_R$ be the dilation of $V$ by $R^{-1}$. Its unit-ball excess is
		\[
		E_R=R^{-m-2}I(R)=o(R^{2\gamma_+-2}).
		\]
		For every area-minimizing nonflat quadratic cone $m\ge7$ and
		$-m/2<\gamma_+<0$. The blowdown convergence supplies both area conditions
		in Lemma~\ref{e:es}; the relevant spheres have zero cone measure.
		We apply Lemma~\ref{e:es} with $\lambda=0$ and denote the resulting parameters
		by $z_R,Q_R,\lambda_R$. Their parameters in the original coordinates are
		$A_R:=Rz_R$ and $\ell_R=R\lambda_R$. Equation~\eqref{e:es-parameters}
		gives
		\begin{align*}
			|A_R|&\le cR^{-m/2}I(R)^{1/2}=o(R^{\gamma_+})\longrightarrow0,\\
			|\ell_R|^{1-\gamma_+}
			&\le cR^{-(m+2\gamma_+)/2}I(R)^{1/2}\longrightarrow0,
			\qquad Q_R\longrightarrow\mathrm{Id}.
		\end{align*}
		Thus $A_R+Q_RS_{\ell_R}$ tends smoothly to $\mathcal C_{p,q}$ on each
		fixed compact annulus.
		
		Fix an annulus $s\le |X|\le S$, with $s>0$.
		Its image at scale $R^{-1}$ lies in $B_{c_{s,S}R^{-1}}(z_R)$ for large $R$.
		Equation~\eqref{e:es-graph}, multiplied by $R$, bounds the ambient
		height error there by
		\[
		C_{s,S}R^{-\beta}E_R^{1/2}\longrightarrow0;
		\]
		the corresponding gradient error tends to zero as well. The full graph
		conclusion gives convergence of supports in both directions on a
		slightly smaller annulus. Since the varifold in the original coordinates $V$ is fixed,
		its support away from zero is exactly the regular cone. The constancy
		theorem makes its multiplicity a constant integer there. Rectifiability
		excludes point measure at zero and the unit blowdown fixes the integer to
		be one. This proves the assertion.
		The proof is complete.
	\end{proof}
	
	\begin{prop}\label{e:minimizing-s3}
		Let $\C$ be an area-minimizing quadratic cone. Then $\C$ satisfies
		\textup{(S3)}.
	\end{prop}
	\begin{proof}
		Let $V$ satisfy the hypotheses in \textup{(S3)}. Recall that
		\[
		\mathcal E_H(V)=\int_{\R^{m+1}}\frac{\operatorname{dist}(X,\mathcal C)^2}{1+|X|^2}
		\dd \|V\|(X)<\infty.
		\]
		For fixed $R_0$ and $R>R_0$, we have
		\[
		I(R)\le I(R_0)+(1+R^2)
		\int_{|X|>R_0}\frac{\operatorname{dist}(X,\mathcal C)^2}{1+|X|^2}
		\dd \|V\|(X).
		\]
		We divide by $R^2$, let $R\to\infty$, and then let $R_0\to\infty$.
		Thus we obtain $I(R)=o(R^2)$. On the other hand, we know
		\[
		m+2\gamma_+=2+\sqrt{m^2-8m+8}>2\quad(m\ge7).
		\]
		Thus we obtain \eqref{e:subgrowth-condition}, and
		Lemma~\ref{e:subgrowth} proves $V=|\mathcal C|$.
		This concludes the proof.
	\end{proof}
	In the normalization \eqref{e:quadratic}, all quadratic cones with $m>7$
	are area-minimizing, and those with $m=7$ are area-minimizing exactly when
	$|p-q|\le2$ \cite[Theorem~1.1]{Davini04}. Thus, up to exchanging the factors,
	it remains to prove weighted rigidity directly for $\C_{1,5}$.
	
	\begin{prop}\label{e:exceptional-s3}
		Let
		\[
		\mathcal C=\mathcal C_{1,5}
		=\{(x,y)\in\mathbb R^2\oplus\mathbb R^6:|y|^2=5|x|^2\}
		\]
		Then $\C$ satisfies \textup{(S3)}.
	\end{prop}
	\begin{proof}
		Its complementary components will be denoted by
		\[
		E_s:=\{|y|<\sqrt5|x|\},\qquad E_c:=\{|y|>\sqrt5|x|\}.
		\]
		We first record the foliations and maximum principles used below.
		
		By \cite[p.~210]{Lin}, $\C_{1,5}$ is one-sided strictly area-minimizing
		in $\overline{E_s}$. The side $E_s$ is foliated by the dilations
		$H_\lambda=\lambda H_1$, $\lambda>0$, of a connected, complete, smooth
		embedded minimal hypersurface $H_1$ \cite[pp.~210 and~213]{Lin}. Lin's asymptotic expansion for its
		radial generating curve has leading difference $u-v=u^{-2}+o(u^{-2})$.
		We convert the fixed radial orbit coordinates to normal height over
		$\mathcal C$ to obtain, for each fixed $\lambda>0$,
		\begin{equation}\label{e:slow-leaf}
			h_\lambda(r,\theta)\ge c_\lambda r^{-2}\quad(r\ge R_\lambda),
			\qquad c_\lambda>0.
		\end{equation}
		The positive sign denotes displacement into $E_s$.
		
		Since the cone is regular, stable, and area-minimizing on exactly one side,
		we may apply \cite[Lemma~8.2]{Wang}; the cone over $S^5\times S^1$
		is explicitly among its examples. In the present notation it provides
		a compact connected smooth minimal hypersurface with boundary
		\begin{equation}
			\mathscr K\subset\overline{E_c},\qquad
			\partial \mathscr K=\mathcal C\cap S^7,\qquad
			\mathscr K\cap\mathcal C=\partial \mathscr K,\qquad \operatorname{dist}(0,\mathscr K)>0,
		\end{equation}
		such that the interiors of $\mathscr K_\tau=\tau \mathscr K$, $\tau>0$, foliate $E_c$.
		
		We also use two varifold maximum principles. White's boundary principle
		\cite[Theorem~4]{WhiteMP} implies that a stationary integral hypersurface
		varifold supported on one side of a connected smooth minimal
		hypersurface, and touching that hypersurface, contains the entire
		hypersurface as a positive integer stationary summand. Its nonnegative
		stationary remainder is disjoint from the boundary. The assertion is
		applied in the relevant open ambient manifold, so a cap interior is an
		admissible boundary. Wickramasekera's disjointness theorem
		\cite[Theorem~3.2]{Wick} says that two stationary integral hypersurface
		$m$-varifolds whose supports intersect in a set of zero
		$\mathcal H^{m-1}$ measure have disjoint supports. In particular, for
		$m=7$ their only contact cannot be an isolated point.
		
		We now take $V$ to be as in \textup{(S3)}. Multiplicity-one blowdown convergence
		to the regular cone implies that the whole support outside a compact
		set is a single small normal graph $h$ over $\mathcal C$. To see the
		uniform nature of this assertion, we apply Allard regularity to the
		rescalings on a fixed larger annulus; the stationary minimal graph
		equation then gives smooth convergence. Failure at any sequence of
		radii would contradict the prescribed blowdown. Local monotonicity
		excludes additional components away from the graph. Thus the scaled
		$C^2$ norm of $h$ tends to zero.
		
		Finite Hardy distance gives more precise decay. Define
		\[
		\varepsilon_R^2:=
		\int_{R/4<|X|<4R}
		\frac{\operatorname{dist}(X,\mathcal C)^2}{1+|X|^2}
		\dd \|V\|(X).
		\]
		Then we have $\varepsilon_R\to0$. For the rescaled graph height
		$\widetilde h_R(Y)=R^{-1}h(RY)$, graph Jacobians and the normal tube comparison
		imply
		\[
		\|\widetilde h_R\|_{L^2(A_{1/2,2}\cap\mathcal C)}^2
		\le CR^{-7}\varepsilon_R^2.
		\]
		The small minimal-graph equation is a uniformly elliptic linear equation
		for $\widetilde h_R$ with uniformly bounded coefficients. We subtract the
		zero graph equation and integrate its linearization along the segment
		of graphs. The local Caccioppoli estimate and Sobolev iteration give
		\[
		\|\widetilde h_R\|_{L^\infty(A_{3/4,3/2})}
		\le C\|\widetilde h_R\|_{L^2(A_{1/2,2})}.
		\]
		Thus we obtain
		\begin{equation}\label{e:hardy-pointwise}
			\sup_{\mathcal C\cap A_{3R/4,3R/2}}|h|
			\le CR^{-5/2}\varepsilon_R=o(R^{-5/2}).
		\end{equation}
		We obtain the same conclusion if we shrink or enlarge the annuli by fixed factors.
		
		We first consider the leaf coordinate $\lambda:E_s\to(0,\infty)$ of the
		foliation $H_\lambda$. It extends continuously by zero to the cone.
		For clarity, this extension follows by writing $H_1$ as a radial graph
		$r=R(\omega)$ over $E_s\cap S^7$. Its asymptotically conical end has
		$R(\omega)\to\infty$ as $\omega$ approaches the boundary link, while
		$\lambda(r\omega)=r/R(\omega)$. Also $\lambda(X)\le C|X|$, because
		$H_1$ has positive distance from the origin. This proves continuity
		also at the vertex. The foliation order, \eqref{e:slow-leaf}, and
		\eqref{e:hardy-pointwise} show that
		\[
		\lambda(X)\longrightarrow0\quad\text{along }
		\operatorname{spt}V\cap E_s\text{ as }|X|\to\infty.
		\]
		Indeed, for every fixed positive $\lambda_0$, the graph of $V$ eventually
		lies between the cone and $H_{\lambda_0}$ whenever it is on this side.
		
		If the support meets $E_s$, its leaf coordinate therefore attains a
		positive maximum $\lambda_0$. Compactness of a maximizing sequence
		follows from decay at infinity, and its limit cannot lie on the cone
		because the extended coordinate vanishes there. In the open ambient
		manifold $E_s$, the restricted varifold is stationary, supported in
		$\{0<\lambda\le\lambda_0\}$, and touches its connected minimal
		boundary $H_{\lambda_0}$. White's principle gives a positive integer
		summand on that entire leaf. But its Hardy distance diverges:
		\[
		\int_{H_{\lambda_0}}
		\frac{\operatorname{dist}(X,\mathcal C)^2}{1+|X|^2}\dd \mathcal H^7
		\ge c\int_{R_0}^\infty r^6r^{-4}r^{-2}\dd r=\infty.
		\]
		This contradicts the Hardy bound for $V$. Therefore we have
		\begin{equation}\label{e:no-slow-side}
			\operatorname{spt}V\cap E_s=\varnothing.
		\end{equation}
		
		We next suppose that the support touches the punctured cone. In
		$\mathbb R^8\setminus\{0\}$, equation~\eqref{e:no-slow-side} puts it on
		one side of the connected smooth minimal boundary
		$\mathcal C\setminus\{0\}$. White's principle gives
		\[
		V=d|\mathcal C|+V',\qquad d\in\mathbb Z_{>0},
		\]
		where the remainder is nonnegative and stationary away from zero.
		A rectifiable seven-varifold has no atom at zero, so the equality holds
		on all of $\mathbb R^8$ as a varifold identity. Since $V$ and the cone
		are globally stationary, subtraction makes $V'$ globally stationary
		as well. The unit blowdown implies $d=1$. If $V'$ were nonzero,
		stationary monotonicity at a positive-density point would give
		$\|V'\|(B_R(0))\ge cR^7$ for large $R$, whereas its blowdowns converge
		to zero. Thus we obtain $V'=0$ and $V=|\mathcal C|$.
		
		The remaining possibility of contact only at the vertex is excluded
		by Wickramasekera's theorem: the singleton has zero $\mathcal H^6$
		measure, so such an intersection of the stationary supports is empty.
		Unless the conclusion already holds, we have therefore proved
		\begin{equation}\label{e:only-cap-side}
			\operatorname{spt}V\subset E_c,\qquad
			\operatorname{spt}V\cap\mathcal C=\varnothing.
		\end{equation}
		
		Finally, let $\tau:E_c\to(0,\infty)$ be the coordinate of the compact-cap
		foliation. Compactness of $\mathscr K$ and its positive distance from zero give
		\[
		c_0\tau(X)\le |X|\le C_0\tau(X)\quad(X\in E_c)
		\]
		for constants $0<c_0<C_0<\infty$. Thus it follows that $\tau$ tends to infinity along
		an unbounded sequence. A sequence in the support with $\tau\to0$
		would converge to zero, contradicting closedness and
		\eqref{e:only-cap-side}. A minimizing sequence for $\tau$ is therefore
		bounded and bounded away from zero; its support limit cannot lie on
		the cone by the same closedness. Thus it follows that $\tau$ attains a
		positive minimum $\tau_0$ on the support. In $E_c$, the varifold is
		supported in $\{\tau\ge\tau_0\}$ and touches the connected minimal
		boundary $\operatorname{int}\mathscr K_{\tau_0}$. White's principle implies that its
		support contains the whole cap interior. Closedness in the full
		ambient space then implies that it contains
		$\partial \mathscr K_{\tau_0}=\tau_0(\mathcal C\cap S^7)$, contradicting
		\eqref{e:only-cap-side}. This completes the proof.
	\end{proof}
	
	\begin{cor}\label{e:all-quadratics}
		Let $\mathcal C_{p,q}$ be a quadratic cone with $p,q\ge1$ and
		$m=p+q+1\ge7$. Then $\mathcal C_{p,q}$ satisfies
		\textup{(S1), (S2), and (S3)}. Thus it follows that all these
		cones fall under Theorem~\ref{thm:main}, including the nonminimizing
		cone $\mathcal C_{1,5}$ and its factor exchange.
	\end{cor}
	\begin{proof}
		By Proposition~\ref{e:spectrum}, we obtain (S1) and (S2).
		For an area-minimizing pair, we apply Proposition~\ref{e:minimizing-s3}.
		For the remaining pair, we apply Proposition~\ref{e:exceptional-s3}. Weighted rigidity is invariant under
		ambient orthogonal changes of coordinates, so factor exchange preserves
		it.
		This finishes the proof.
	\end{proof}
	
	\begin{proof}[Proof of Corollary~\ref{cor:examples}]
		Proposition~\ref{e:strict-s3} verifies (S3) for a strictly minimizing
		cone. Theorem~\ref{thm:main} proves the first assertion. For quadratic
		cones, Proposition~\ref{e:spectrum} proves (S1) and (S2), and
		Corollary~\ref{e:all-quadratics} proves (S3). Thus the quantitative
		conclusions of Theorem~\ref{thm:energy-loja} also apply.
		Thus we conclude the proof.
	\end{proof}

	Hanbing Fang, Department of Mathematics, The Chinese University of Hong
	Kong, Shatin, N.T., Hong Kong; E-mail: hanbingfang@cuhk.edu.hk\\

	Yongheng Han, Institute of Mathematics, Academy of Mathematics and Systems Science,
	Chinese Academy of Sciences, Beijing 100190, China; E-mail: hyh98@amss.ac.cn

\end{document}